\documentclass[a4paper,11pt]{article}

\usepackage{ucs}
\usepackage[utf8x]{inputenc}
\usepackage{graphicx}
\usepackage{amsfonts}
\usepackage{dsfont}
\usepackage{amssymb}
\usepackage{amsmath}
\usepackage{amsthm}
\usepackage{enumerate}
\usepackage{stmaryrd}
\usepackage{fullpage}
\usepackage{ifthen}
\usepackage{subfigure}
\usepackage{epic}
\usepackage{authblk}
\usepackage{textcomp}
\usepackage{wrapfig}
\usepackage[scriptsize]{caption}
\usepackage{url}
\usepackage{eurosym}
\usepackage{float}
\usepackage{mathrsfs}

\usepackage{tikz,tkz-tab}
\usetikzlibrary{matrix}

\usepackage{xcolor}
\usepackage{titlesec}

\graphicspath{{eps/}{pdf/}}

\newcommand{\K}{\mathcal{K}}

\newcommand{\R}{\mathbb{R}}
\newcommand{\N}{\mathbb{N}}
\newcommand{\Z}{\mathbb{Z}}
\newcommand{\C}{\mathbb{C}}

\newcommand{\ds}{\mathrm{d} s}
\renewcommand{\d}{\mathrm{d}}
\newcommand{\ee}{\mathrm{e}}
\renewcommand{\Re}{\mathrm{Re}}
\renewcommand{\Im}{\mathrm{Im}}
\newcommand{\Arg}{\mathrm{Arg}}

\newcommand{\bqq}{\begin{equation}}
\newcommand{\eqq}{\end{equation}}
\newcommand{\bqs}{\begin{equation*}}
\newcommand{\eqs}{\end{equation*}}

\newtheorem{proposition}{Proposition}[section]
\newtheorem{theorem}{Theorem}

\makeatletter
\newcommand{\leqnomode}{\tagsleft@true\let\veqno\@@leqno}
\makeatother

\numberwithin{equation}{section}

\newenvironment{Proof}[1][.]%
 {\begin{trivlist}\item[]\textbf{Proof#1 }}%
 {\hspace*{\fill}$\rule{0.3\baselineskip}{0.35\baselineskip}$\end{trivlist}}

\title{Spreading Properties of a City-Road Reaction-Diffusion Model on One-Dimensional Lattice}

\author[1]{Gr\'egory Faye\footnote{email: \texttt{gregory.faye@math.univ-toulouse.fr}}}
\author[1]{Jean-Michel Roquejoffre}
\author[2]{Min Zhao}
\affil[1]{\small Univ Toulouse, CNRS, Institut de Math\'ematiques de Toulouse, Toulouse, France}
\affil[2]{\small Aix Marseille Univ, CNRS, Institut de Math\'ematiques de Marseille, Marseille, France}

\begin{document}
\maketitle

\begin{abstract}
We propose and study a new model to describe biological invasions constrained on infinite homogeneous one dimensional metric graphs. Our model consists of an infinite PDE-ODE system where, at each vertex of the one-dimensional lattice $\Z$, we have a logistic equation, and connections between vertices are given by diffusion equations on the edges supplemented with Robin like boundary conditions at the vertices. We establish the main properties of the system and study the long time behavior of the solutions, especially by characterizing an asymptotic spreading speed for the system. In the fast diffusion regime, we derive a novel asymptotic model which exhibits similar propagation properties as the classical discrete Fisher-KPP on the one-dimensional lattice $\Z$.
\end{abstract}



\section{Introduction}

Traveling waves in biology are ubiquitous and have been found in many contexts, such as the spread of cancer cells in healthy tissue, traveling bands of
bacteria, the diffusion of genes within a population, or the spread of an epidemic, to name a few. One common feature of these biological spreading phenomena is that they are highly  complex, network-driven dynamic processes. In many applications, the intrinsic heterogeneity of the underlying networks makes it very challenging to analyze these processes and assess the relevant factors that effectively drive the propagation. From a modeling perspective, it is also quite difficult, given the multiscale nature of the considered biological processes, to adopt a formalism that could combine intricate network structures and complex dynamics, and still provide comprehensive and valuable feedbacks to the biological community. 

Our focus here will be on biological processes that can be well approximated by macroscopic models set on metric graphs. More precisely, given a metric graph, that is, a collection of interconnected vertices and edges with prescribed lengths, we will consider  non classical reaction-diffusion models where, schematically, diffusion processes take place along the edges of the graph while reaction kinetics occur at the vertices with prescribed rules of exchanges between vertices and adjacent edges. Such a formalism has typically been proposed to study coupled membrane-bulk diffusion systems \cite{Gou-Ward-2016,Gou-Li-Nagata-Ward-2015,Paquin-Lefebvre-Nagata-Ward-2020} and to analyze the effects of transportation networks such as roads, railways or waterways on the spread of epidemics among cities \cite{BF21a,KDB} as reported for example for the spread of COVID-19, Chikungunya virus, Zika virus and HIV virus \cite{Gatto-Bertuzzo-2020,Hale-Thomas-2019,Faria-Rambaut-2014}. Other types of reaction-diffusion models on metric graphs have been proposed in the past decades. In population dynamics, the so-called river network models \cite{JPS,DLPZ} describe the dynamics of organisms living in a river system subject to a forced flow in the downstream direction. It typically consists of reaction-diffusion equations set on the vertices of a given prescribed network with a continuity condition at the edges, together with a Kirchoff law that translates the continuity of fluxes through the edges. In cellular physiology, models of cells coupled by gap junctions \cite{KS09,RB90,PCB22} are typically set on networks, and concentrations of diffusing particles follow a diffusion equation within each cell, idealized by an edge, and at the junction between two cells, that is, at each vertex of the network, specific boundary conditions are prescribed to account for the permeability properties of the cells membrane.

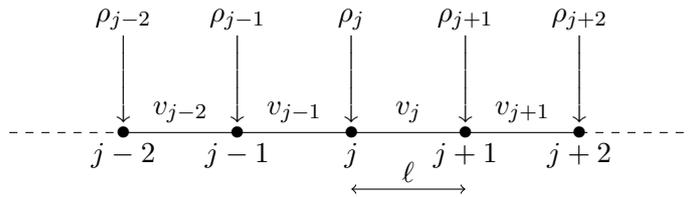
\begin{figure}
\begin{center}
\begin{tikzpicture}[scale=1.5]
    \coordinate (A1) at (0,0);
    \coordinate (A2) at (1,0);
    \coordinate (A3) at (2,0);
    \coordinate (A4) at (3,0);
    \coordinate (A5) at (4,0);
    \coordinate (A6) at (5,0);
    \coordinate (A7) at (6,0);

   \coordinate (B1) at (1.5,0);
    \coordinate (B2) at (2.5,0);
    \coordinate (B3) at (3.5,0);
    \coordinate (B4) at (4.5,0);
    
    \coordinate (C1) at (1,1);
    \node at (C1) {$\rho_{j-2}$};
    
    \coordinate (C2) at (2,1);
    \node at (C2) {$\rho_{j-1}$};
    
    \coordinate (C3) at (3,1);
    \node at (C3) {$\rho_{j}$};
    
    \coordinate (C4) at (4,1);
    \node at (C4) {$\rho_{j+1}$};

 \coordinate (C5) at (5,1);
    \node at (C5) {$\rho_{j+2}$};

  \node[above] at (B1) {$v_{j-2}$};
    \node[above] at (B2) {$v_{j-1}$};
    \node[above] at (B3) {$v_{j}$};
    \node[above] at (B4) {$v_{j+1}$};

    \node[below] at (A2) {$j-2$};
    \node[below] at (A3) {$j-1$};
    \node[below] at (A4) {$j$};
    \node[below] at (A5) {$j+1$};
    \node[below] at (A6) {$j+2$};

  \node[above] at (A2) {$\Bigg\downarrow$};
    \node[above] at (A3) {$\Bigg\downarrow$};
    \node[above] at (A4) {$\Bigg\downarrow$};
    \node[above] at (A5) {$\Bigg\downarrow$};
    \node[above] at (A6) {$\Bigg\downarrow$};

\draw[<->] (3,-0.5) -- (4,-0.5);
 \node at (3.5,-0.35) {$\ell$};

    \node at (A2) {$\bullet$};
    \node at (A3) {$\bullet$};
    \node at (A4) {$\bullet$};
    \node at (A5) {$\bullet$};
    \node at (A6) {$\bullet$};
   
    \draw[dashed] (A1) -- (A2);
    \draw  (A2) -- (A3) -- (A4) -- (A5) -- (A6);
    \draw[dashed] (A6) -- (A7);

\end{tikzpicture}
\end{center}
\caption{Schematic spatial configuration of the system \eqref{eq1}-\eqref{bd} where the unknowns $\rho_j$ are indexed on the lattice $\Z$ while each $v_j$ is locally defined on $(0,\ell)$.}\label{fig:model}
  \end{figure}

In the present work, we propose a model where the underlying metric graph is indexed by the one-dimensional lattice $\Z$, such that each vertex of the graph is exactly connected to two incident edges, and all edges have exactly the same length, denoted by $\ell>0$. As previously emphasized, our framework is very general and relevant in a wide array of situations in biology. Nevertheless, for convenience, we will adopt a population dynamics point of view and for the sake of simplicity and illustration, we shall from now on refer to the vertices as ``the cities'' and the edges ``the roads''. For $j\in\Z$, we denote by $\rho_j(t)$ the density of individuals which reside in city $j$ while we denote by $v_j(t,x)$ the density of individuals diffusing along the road connecting city $j$ to city $j+1$ where $t$ refers to time and $x\in[0,\ell]$ represents the local position on the road. We refer to Figure~\ref{fig:model} for a schematic illustration of the spatial configuration of our so-called ``city-road'' model. Exchanges of populations take place between cities and adjacent roads. Namely, given a city, indexed by $j\in\Z$, a fraction $\alpha>0$ of individuals from the two adjacent roads at the city, that are $v_{j-1}(t,\ell)$ and $v_{j}(t,0)$, joins the city, while a fraction $\beta>0$ of individuals from the city $j$ transfers to each of the two adjacent roads. It is further  assumed that the population in each city is subject to a logistic-type growth, resulting in a nonlinear reaction term $f(\rho)$ that models effective birth rate and intrinsic competition. On the other hand, we assume that no such reaction is relevant on the roads and consider solely a diffusion process, with diffusion coefficient $d>0$, to describe the dynamics of each $v_j$. Transposing the above principles into equations, we are thus led to consider the following system of equations:
\bqq\label{eq1}
\forall t>0,~j\in\Z, \quad \left\{
	\begin{split}
		\partial_{t} v_{j}(t, x)&=d\partial_{x}^2v_{j}(t,x),\quad x\in(0,\ell),\\
		\rho_{j}^{\prime}(t)&=f(\rho_{j}(t))+\alpha(v_{j}(t,0)+v_{j-1}(t,\ell))-2\beta\rho_{j}(t),
	\end{split}\right.
\eqq
with inhomogeneous Robin boundary conditions
\begin{equation}\label{bd}
\forall t>0,~j\in\Z, \quad	\begin{cases}
		-d\partial_{x}v_{j}(t,0)+\alpha v_{j}(t,0)=\beta \rho_{j}(t),\\
		d\partial_{x}v_{j}(t,\ell)+\alpha v_{j}(t,\ell)=\beta \rho_{j+1}(t).
	\end{cases}
\end{equation}
Regarding the nonlinearity $f$, we shall assume the following running assumptions.
\bqs\leqnomode
\begin{array}{l}
f\in\mathscr{C}^1([0,1]), \,  f(0)=f(1)=0, \, f'(0)>0, \, f'(1)<0, \text{ and } \forall u\in(0,1), 0<f(u)\leq f'(0) u, \\
\text{ furthermore } f \text{ is extended linearly on } (-\infty,0)\cup(1,+\infty).
\end{array}\tag{H}
\eqs
Let us already remark that by performing the following rescaling
\bqs
x'\longleftrightarrow \frac{x}{\ell}, \quad \widetilde{v}_{j}(t,x') \longleftrightarrow \ell v_{j}(t,x), \quad d' \longleftrightarrow \frac{d}{\ell^2}\quad\text{ and }\quad \alpha' \longleftrightarrow \frac{\alpha}{\ell},
\eqs
we may assume, for the rest of the paper, and without loss of generality, that $\ell=1$.

Model \eqref{eq1}--\eqref{bd} is largely inspired by the SIR model proposed by the first author and Besse in \cite{BF21a}. There are nevertheless three important differences between the two models. First of all, the intrinsic dynamics at each city are different, in our case it is given by a single logistic equation, while in \cite{BF21a} it was given by an SIR compartment model resulting in a system of equations. Second, the study \cite{BF21a} considered compact connected graphs, meaning that the number of cities and roads was finite, while here model \eqref{eq1}--\eqref{bd} is indexed by the one-dimensional lattice $\Z$ and thus infinite. Finally, the model in \cite{BF21a} allowed a fraction of individuals to pass from one road to another one. This is not taken into account in the boundary conditions \eqref{bd} and we refer to the last section of the present manuscript for a longer discussion about this possible extension into the model. One of the key feature of system \eqref{eq1}--\eqref{bd} is the preservation of the total population in the absence of reaction kinetics at the cities. Indeed, assume that $(\mathbf{v},\boldsymbol{\rho})$ with $\mathbf{v}=(v_j)_{j\in\Z}$ and $\boldsymbol{\rho}=(\rho_j)_{j\in\Z}$ is a solution of \eqref{eq1}--\eqref{bd} with $f=0$ and such that the following quantity is well defined for all time $t\geq0$ for which the solution exists:
\bqs
M(t):=\sum_{j\in\Z} \left[\rho_j(t)+\int_0^1 v_j(t,x)\d x\right].
\eqs 
Then, formally, integrating by parts in the first equation and using the boundary conditions, we obtain
\begin{align*}
\int_0^1 v_j(t,x)\d x - \int_0^1 v_j(0,x)\d x & = d \int_0^t \left( \partial_{x}v_{j}(s,1)-\partial_{x}v_{j}(s,0)\right)\ds \\
&=\beta  \int_0^t \left(\rho_j(s)+\rho_{j+1}(s)\right)\ds- \alpha \int_0^t \left(v_{j}(s,1)+v_{j}(s,0)\right)\ds,
\end{align*}
while the second equation gives
\bqs
\rho_j(t)-\rho_j(0)= \alpha \int_0^t(v_{j}(s,0)+v_{j-1}(s,1))\ds-2\beta\int_0^t\rho_{j}(s)\ds.
\eqs
Summing over $\Z$, we deduce that $M(t)=M(0)$ for all $t\geq0$. As a consequence, we see that, in the absence of reaction kinetics, the exchanges between the cities and the roads exactly compensate each other.

Our aim here is to study the long time behavior of the solutions \eqref{eq1}--\eqref{bd} as a function of the various parameters of the model: $\alpha$, $\beta$, and $d$, the nonlinearity $f$ and the chosen initial condition. We are especially interested in characterizing the spreading properties of the system. More precisely, given a compactly supported initial condition, that is, given an initial condition for which only finitely many cities and/or roads have a nonzero initial population, does the corresponding solution of the Cauchy problem converge to a unique positive stationary configuration? And if yes, at which speed does the convergence towards this eventual steady state take place? In a nutshell, our main results regarding our model \eqref{eq1}--\eqref{bd} are as follows. At this stage of the presentation, we remain formal and refer to the following sections for precise statements and assumptions.

\paragraph{Existence and uniqueness of classical solutions.} We prove in Theorem~\ref{thmcauchy} below that for each well-prepared initial condition our model \eqref{eq1}--\eqref{bd} admits a unique positive bounded classical solution which is global in time. The structure of \eqref{eq1}--\eqref{bd} is non standard, and since the graph considered here is infinite, we cannot readily rely on the existing results of \cite{BF21a}, which only apply for compact graphs. We adopt a similar approach and construct solutions via an iterative scheme. To obtain compactness and extract converging subsequences, we combine a priori estimates via comparison principle techniques and standard parabolic estimates for the heat equation with inhomogeneous Robin boundary conditions. This analysis is conducted in Section~\ref{CP}.

\paragraph{Long time behavior of the solutions.} We fully characterize the long time behavior of the unique solution of our model. More precisely, in Theorem~\ref{thm-sta}, we first prove that the only positive, bounded, stationary solution of \eqref{eq1}--\eqref{bd} is the constant sequence $\left(\frac{\beta}{\alpha},1\right)_{j\in\Z}$. Interestingly enough, the proof relies on the fact that stationary solutions of \eqref{eq1}--\eqref{bd} are in one-to-one correspondence with the stationary solutions of the discrete Fisher-KPP equation given by
\bqs
\lambda\left( \rho_{j-1}-2\rho_j+\rho_{j+1}\right)+f(\rho_j)=0, \quad j\in\Z,
\eqs
for some $\lambda>0$ depending explicitly on $\alpha$, $\beta$ and $d$. Finally, we demonstrate that the positive bounded stationary solution $\left(\frac{\beta}{\alpha},1\right)_{j\in\Z}$ is the global attractor of the system \eqref{eq1}--\eqref{bd} when initialized with nontrivial nonnegative bounded initial condition. We refer to Theorem~\ref{long} for a precise statement but we already emphasize that the convergence is  locally uniform in $j\in\Z$ and uniform in $x\in[0,1]$. The aforementioned results are proved in Section~\ref{secLTB} and rely on comparison principle techniques and the construction of adequate sub and super-solutions for the system \eqref{eq1}--\eqref{bd}.

\paragraph{Linear spreading speed.} In Section~\ref{secEXP}, we analyze the linearized problem around the trivial constant state $\left(0,0\right)_{j\in\Z}$ and derive a theoretical formula for the linear spreading speed, denoted by $c_*$, and defined as the smallest possible speed $c>0$ for which there exists an exponential solution of the form
\bqs	(v_{j}(t,x),\rho_{j}(t))=\left(\ee^{-\mu(j-ct)}V(x),\ee^{-\mu(j-ct)}\right),
\eqs
for some prescribed positive profile $V$. The formula for $c_*$ is given in equation \eqref{c_*} below and we refer to Figure~\ref{fig:Speed} and Figure~\ref{fig:SpeedAB} in Subsection~\ref{subsecnumspeed}  for illustrations of the dependence of $c_*$ as a function of the other parameters $\alpha$, $\beta$, $d$ and $f'(0)$. The characterization leading to the definition of $c_*$ is quite intricate. Although we manage to prove that $c_*$ is well-defined, it is yet a problem to prove that there exists a unique corresponding $\mu_*>0$ at which the spreading speed is attained, as it is usually the case for reaction-diffusion systems having a monotone structure. We conjecture that it is indeed the case based on our numerical computation of the linear spreading speed via its formula \eqref{c_*}.

\paragraph{Asymptotic spreading.}  It turns out that the linear spreading speed $c_*$ defined in formula \eqref{c_*} is precisely the asymptotic spreading speed of the nonlinear system \eqref{eq1}-\eqref{bd} as proved in Theorem~\ref{thmSpSp} in Section~\ref{secAS}. More precisely, we show that solutions of system \eqref{eq1}-\eqref{bd} starting from compactly supported initial conditions spread at speed $c_*$. Translated mathematically, if $(\mathbf{v},\boldsymbol{\rho})$ is a corresponding solution, then we have the following dichotomy:
\begin{itemize}
\item[(i)] for all $c>c_*$, we have
\bqs
\underset{t\to+\infty}{\lim}~\underset{\substack{|j| \geq ct\\ x\in[0,1]}}{\sup}(v_{j}(t,x),\rho_{j}(t))=\left(0,0\right);
\eqs
\item[(ii)] for all $c\in(0,c_*)$, we have
\bqs
\underset{t\to+\infty}{\lim}~\underset{\substack{|j| \leq ct  \\ x\in[0,1]}}{\inf}(v_{j}(t,x),\rho_{j}(t))=\left(\frac{\beta}{\alpha},1\right).
\eqs
\end{itemize}
A key element of the proof is the construction of compactly supported generalized subsolutions for the nonlinear system \eqref{eq1}-\eqref{bd}.
 
\paragraph{Large diffusion limit.} We finally investigate the large diffusion limit $d\rightarrow+\infty$ of the system in Section~\ref{secLDL}. Our first result, see Theorem~\ref{thmLD} for a precise statement, ensures that for well-prepared initial conditions, the solution of system \eqref{eq1}-\eqref{bd} converges as $d\rightarrow+\infty$ towards $(\mathbf{V},\mathbf{P})$ locally uniformly in $(t,j)\in(0,+\infty)\times\Z$ and uniformly in $x\in[0,1]$, and $(\mathbf{V},\mathbf{P})$ is the solution of the asymptotic system
\bqs
\forall t>0,~j\in\Z, \quad \left\{
	\begin{split}
		 V_j'(t)&=-2\alpha V_j(t)+\beta(P_j(t)+P_{j+1}(t)),\\
		P_{j}^{\prime}(t)&=f(P_{j}(t))+\alpha(V_{j}(t)+V_{j-1}(t))-2\beta P_{j}(t).
	\end{split}\right.
\eqs
For this asymptotic system, we also prove that $\left(\frac{\beta}{\alpha},1\right)_{j\in\Z}$ is the only positive bounded stationary solution and that solutions to the corresponding Cauchy problem starting from bounded nonnegative initial conditions asymptotically converge towards it, locally uniformly in $j\in\Z$. We further prove in Theorem~\ref{thmASdinf} the existence of an asymptotic spreading speed, denoted by $c_*^\infty$ (see formula~\ref{c_11}), for the asymptotic system. We also conjecture that
\bqs
c_* \underset{d\rightarrow+\infty}{\longrightarrow}c_*^\infty,
\eqs
 where $c_*$ is the spreading speed of the full system \eqref{eq1}-\eqref{bd}, and leave it to future work to rigorously demonstrate this asymptotic limit. Let us further mention that this conjecture is verified numerically in Figure~\ref{fig:Speed}(c) below.
 
 Our asymptotic spreading result echoes the ones obtained for standard reaction-diffusion equations set on graphs such as, for instance, the Fisher-KPP equation set on the lattice \cite{Weinberger-1982} or homogeneous trees \cite{Hoffman-Holzer-2019}, and where the linear spreading speed characterizes the long time behavior of the solutions of the Cauchy problem starting from compactly supported initial data. We also refer to \cite{Besse-Faye-Roquejoffre-Zhang-2023} for the most recent results in the direction of the so-called logarithmic Bramson correction for the level sets of the solutions for the Fisher-KPP equations on the lattice. In our setting, as expected, the characterization of the spreading speed is less explicit and more intricate. Let us also remark that our framework is at the crossroad of the aforementioned standard discrete reaction-diffusion models and continuous models that take into account lines of transportation such as the so-called ``field-road" model of Berestycki, Roquejoffre and Rossi~\cite{Berestycki-Roquejoffre-Rossi-2013}. Indeed, on a formal level, our proposed model can be thought of as being a one-dimensional version of the planar reaction-diffusion system of \cite{Berestycki-Roquejoffre-Rossi-2013}, if we consider only one city and one semi-infinite road.

\section{The Cauchy problem}\label{CP}

In this section, we focus on the well-posedness of the problem \eqref{eq1}--\eqref{bd}. As a consequence, we supplement the system with the initial condition
\begin{equation}\label{inv}
\forall j\in\Z,\quad \begin{cases}
		v_{j}(0,x)=h_{j}(x),\quad x\in(0,1),\\
		\rho_{j}(0)=\Lambda_{j}.
	\end{cases}
\end{equation}
We shall always assume that the initial sequences $\mathbf{h}=(h_j)_{j\in\Z}$ and $\boldsymbol{\Lambda}=(\Lambda_j)_{j\in\Z}$ satisfy the following compatibility condition
\bqq\label{compatibility}
\forall j\in\Z, \quad	\begin{cases}
		-dh_j'(0)+\alpha h_j(0)=\beta \Lambda_{j},\\
		dh_{j}'(1)+\alpha h_{j}(1)=\beta \Lambda_{j+1}.
	\end{cases}
\eqq
Throughout the paper, we let $\mathcal{\ell}^{\infty}(\Z)$ denote the Banach space of bounded valued sequences indexed by $\Z$ and equipped with the norm:
$$
\|\mathbf{u}\|_{\mathcal{\ell}^{\infty}(\Z)}:=\max_{j\in\Z}|u_{j}|, \text{ for } \mathbf{u}=(u_j)_{j\in\Z},
$$
and also define
$$
\mathcal{X}^0:=\left\{\mathbf{u}=(u_j)_{j\in\Z} ~|~ \forall j\in\Z,~ u_{j}\in\mathscr{C}^0([0,1],\R) \text{ and } \|\mathbf{u}\|_{\infty}<+\infty\right\},
$$
with norm
$$
\|\mathbf{u}\|_{\infty}:=\sup_{j\in\Z}\max_{x\in[0,1]}|u_{j}(x)|.
$$

The main result of this section is the following.

\begin{theorem}\label{thmcauchy}
The Cauchy problem \eqref{eq1}-\eqref{bd}-\eqref{inv} with nontrivial nonnegative bounded initial sequences $\mathbf{h}=(h_j)_{j\in\Z} \in \mathcal{X}^0$ and $\boldsymbol{\Lambda}=(\Lambda_j)_{j\in\Z} \in \ell^\infty(\Z)$ satisfying the compatibility condition \eqref{compatibility} admits a unique bounded positive global classical solution $(\mathbf{v},\boldsymbol{\rho})=(v_j,\rho_j)_{j\in\Z}$ with $\rho_j\in\mathscr{C}^1([0,+\infty),\R)$ and
\bqs
v_j\in \mathscr{C}^{0}([0,+\infty)\times[0,1],\R),~ \partial_t v_j, \partial_x^2v_j\in\mathscr{C}^{0}((0,+\infty)\times(0,1),\R), \text{ and } \partial_x v_j\in \mathscr{C}^{0}((0,+\infty)\times[0,1],\R),
\eqs
for all $j\in\Z$. Furthermore, for all $t>0$ and $j\in\Z$, one has 
\bqs
 0<v_{j}(t,x)\leq \max\left\{\dfrac{\beta}{\alpha},\|\mathbf{h}\|_{\infty} , \dfrac{\beta}{\alpha}\|\boldsymbol{\Lambda}\|_{\ell^{\infty}(\Z)}\right\},~ x\in[0,1], 
\eqs
and
\bqs
 0< \rho_{j}(t)\leq \max\left\{\|\boldsymbol{\Lambda}\|_{\ell^{\infty}(\Z)},1,\dfrac{\alpha}{\beta}\|\mathbf{h}\|_{\infty} \right\}.
\eqs
\end{theorem}

\begin{figure}[t!]
\centering
\includegraphics[width=.45\textwidth]{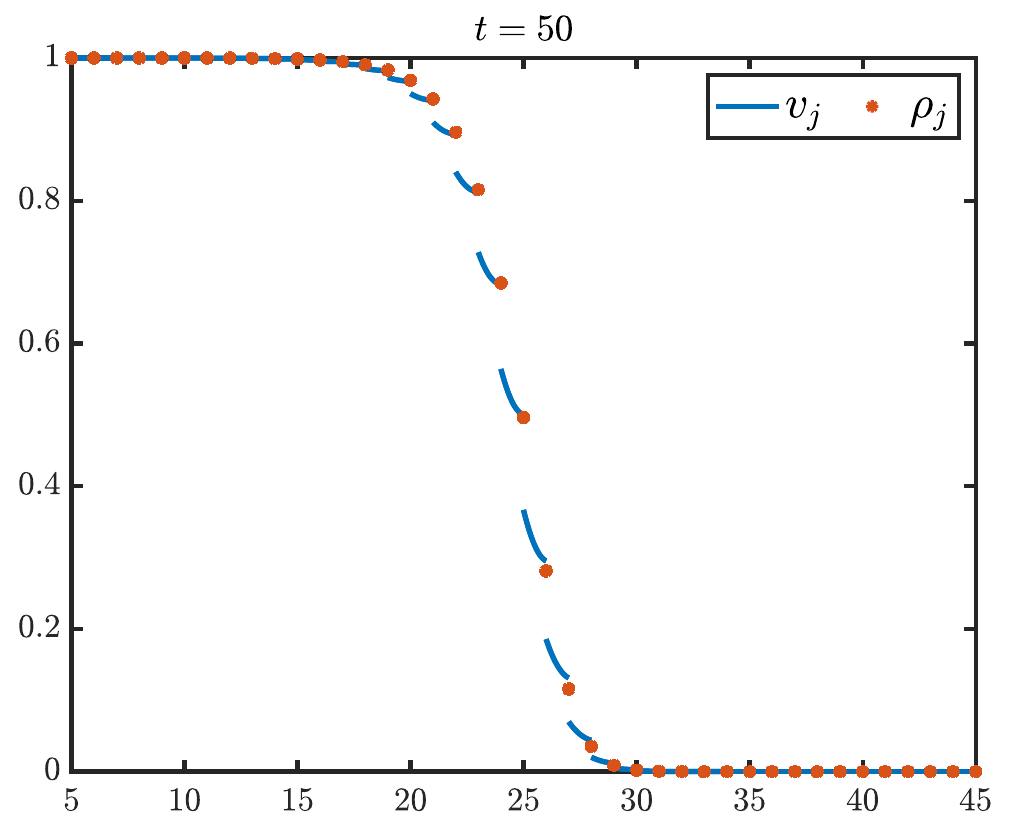}
\caption{Numerically computed solution of system \eqref{eq1}--\eqref{bd} at time $t=50$ for $f(u)=u(1-u)$ and $(\alpha,\beta,d)=(1,1,1)$ starting from an initial condition where $h_j\equiv 0$ for all $j\in\Z$ and $\Lambda_j=1$ for $j\leq0$ and $\Lambda_j=0$ for $j\geq1$. The red dots represent $\rho_j$ located at position $j$ while each blue curve represents $v_j$ located on the interval $[j,j+1]$.}
 \label{fig:sol}
\end{figure}

\subsection{Uniqueness}

In order to establish the uniqueness of the solution of the Cauchy problem \eqref{eq1}-\eqref{bd}-\eqref{inv}, we shall rely on a comparison principle for \eqref{eq1}-\eqref{bd}. We first define the notion of super and subsolutions to \eqref{eq1}-\eqref{bd}. Let $(\overline{\mathbf{v}},\overline{\boldsymbol{\rho}})=(\overline{v}_j,\overline{\rho}_j)_{j\in\Z}$ with $\overline{\rho}_j\in\mathscr{C}^1([0,+\infty),\R)$ and
\bqs
\overline{v}_j\in \mathscr{C}^{0}([0,+\infty)\times[0,1],\R),~ \partial_t \overline{v}_j, \partial_x^2\overline{v}_j\in\mathscr{C}^{0}((0,+\infty)\times(0,1),\R), \text{ and } \partial_x \overline{v}_j\in \mathscr{C}^{0}((0,+\infty)\times[0,1],\R),
\eqs
for all $j\in\Z$. We say that $(\overline{\mathbf{v}},\overline{\boldsymbol{\rho}})$ is a supersolution to \eqref{eq1}-\eqref{bd} if it has the above regularity and satisfies
\bqs
\begin{cases}
\partial_{t} \overline{v}_{j}(t, x) \geq  d\partial_x^2\overline{v}_{j}(t,x),\quad  x\in(0,1),\\
		{\overline{\rho}_{j}}'(t)
		\geq f(\overline{\rho}_{j}(t))+\alpha(\overline{v}_{j}(t,0)+\overline{v}_{j-1}(t,1))-2\beta\overline{\rho}_{j}(t),\\
		-d\partial_{x}\overline{v}_{j}(t,0)+\alpha \overline{v}_{j}(t,0)\geq \beta \overline{\rho}_{j}(t),\\
		d\partial_{x}\overline{v}_{j}(t,1)+\alpha \overline{v}_{j}(t,1)\geq\beta \overline{\rho}_{j+1}(t),
	\end{cases}
\eqs
for all $t>0$. We define similarly a subsolution $(\underline{\mathbf{v}},\underline{\boldsymbol{\rho}})$ to \eqref{eq1}--\eqref{bd} with the same regularity and all above inequalities being reversed.

\begin{proposition}\label{propCPcauchy}
Let $(\underline{\mathbf{v}},\underline{\boldsymbol{\rho}})$ and $(\overline{\mathbf{v}},\overline{\boldsymbol{\rho}})$ be respectively a subsolution and supersolution to \eqref{eq1}--\eqref{bd}. If we assume that $(\underline{\mathbf{v}},\underline{\boldsymbol{\rho}})$ and $(\overline{\mathbf{v}},\overline{\boldsymbol{\rho}})$ are locally bounded in time and satisfy for all $j\in\Z$ that $\underline{v}_j(0,x)\leq\overline{v}_j(0,x)$ for all $x\in[0,1]$ and $\underline{\rho}_j(0)\leq\overline{\rho}_j(0)$, then we have $\underline{v}_j(t,x)\leq\overline{v}_j(t,x)$ and $\underline{\rho}_j(t)\leq\overline{\rho}_j(t)$ for all $t>0$, $x\in[0,1]$ and  $j\in\Z$. Furthermore, if $\underline{\mathbf{v}}(0)\not\equiv\overline{\mathbf{v}}(0)$ or $\underline{\boldsymbol{\rho}}(0)\not\equiv\overline{\boldsymbol{\rho}}(0)$, then we have $\underline{v}_j(t,x)<\overline{v}_j(t,x)$ and $\underline{\rho}_j(t)<\overline{\rho}_j(t)$ for all $t>0$, $x\in[0,1]$ and  $j\in\Z$.
\end{proposition}

For future reference, we shall often write $(\underline{\mathbf{v}},\underline{\boldsymbol{\rho}})\leq(\overline{\mathbf{v}},\overline{\boldsymbol{\rho}})$ as a short-hand notation for the property that $\underline{v}_j(t,x)\leq\overline{v}_j(t,x)$ and $\underline{\rho}_j(t)\leq\overline{\rho}_j(t)$ for all $t>0$, $x\in[0,1]$ and  $j\in\Z$. Let us finally note that the above comparison principle immediately extends to generalized sub and supersolutions, given by the supremum of subsolutions and the infimum of supersolutions respectively.

\begin{Proof}
We start by defining the following sequences $\mathbf{w}=\overline{\mathbf{v}}-\underline{\mathbf{v}}$ and $\mathbf{z}:=\overline{\boldsymbol{\rho}}-\underline{\boldsymbol{\rho}}$ which share the same regularity as the super and subsolutions, and satisfy the following system of equations
\bqs
\forall j\in\Z,\quad \begin{cases}
\partial_{t} w_{j}(t, x) \geq  d\partial_x^2w_{j}(t,x),\quad  x\in(0,1),\\
		{z_{j}}'(t)
		\geq \left(g_j(t)-2\beta\right)z_{j}(t)+\alpha(w_{j}(t,0)+w_{j-1}(t,1)),\\
		-d\partial_{x}w_{j}(t,0)+\alpha w_{j}(t,0)\geq \beta z_{j}(t),\\
		d\partial_{x}w_{j}(t,1)+\alpha w_{j}(t,1)\geq\beta z_{j+1}(t),
	\end{cases}
\eqs
for all $t>0$, together with
\bqs
\forall j\in\Z,\quad \begin{cases}
w_j(0,x) \geq0, \quad x\in[0,1],\\
z_j(0) \geq 0.
\end{cases}
\eqs
In the above system, we have also defined 
\bqs
\forall t>0,~j\in\Z, \quad g_j(t) := \left\{
\begin{split}
\frac{f(\overline{\rho}_{j}(t))-f(\underline{\rho}_{j}(t))}{\overline{\rho}_{j}(t)-\underline{\rho}_{j}(t)},& \quad \overline{\rho}_{j}(t) \neq \underline{\rho}_{j}(t),\\
f'(\overline{\rho}_{j}(t)), &\quad  \overline{\rho}_{j}(t) = \underline{\rho}_{j}(t).
\end{split}
\right.
\eqs
Since $(\underline{\mathbf{v}},\underline{\boldsymbol{\rho}})$ and $(\overline{\mathbf{v}},\overline{\boldsymbol{\rho}})$ are locally bounded in time, we have that for all $T>0$, there exists a constant $C>0$ such that $|g_j(t)|\leq C$ for all $t\in(0,T]$ and $j\in\Z$. As a consequence, we can rely on Proposition~\ref{propCPapp1} to infer that $w_j(t,x)\geq0$ and $z_j(t)\geq0$ for all $t>0$, $x\in[0,1]$ and $j\in\Z$. The same Proposition~\ref{propCPapp1} also ensures that if furthermore $\mathbf{w}(0)\not\equiv0$ or $\mathbf{z}(0)\not\equiv0$, then $w_j(t,x)>0$ and $z_j(t)> 0$ for all $t>0$, $x\in[0,1]$ and $j\in\Z$. This concludes the proof.
\end{Proof}

\begin{Proof}[ of uniqueness of Theorem~\ref{thmcauchy}.]
Assume that $(\mathbf{v}_1,\boldsymbol{\rho}_1)$ and $(\mathbf{v}_2,\boldsymbol{\rho}_2)$ are two bounded positive global classical solutions to \eqref{eq1}-\eqref{bd}-\eqref{inv} starting from the same initial condition $(\mathbf{h},\boldsymbol{\Lambda})$. Applying twice the comparison principle of Proposition~\ref{propCPcauchy}, we readily obtain that $(\mathbf{v}_1,\boldsymbol{\rho}_1)\equiv(\mathbf{v}_2,\boldsymbol{\rho}_2)$.
\end{Proof}

\subsection{Existence}

Throughout this section, we shall always assume that the nontrivial nonnegative bounded initial sequences $\mathbf{h}=(h_j)_{j\in\Z} \in \mathcal{X}^0$ and $\boldsymbol{\Lambda}=(\Lambda_j)_{j\in\Z} \in \ell^\infty(\Z)$ satisfy the compatibility condition \eqref{compatibility}.

To establish the existence of a solution to system \eqref{eq1}-\eqref{bd}-\eqref{inv}, we shall construct a monotone iterative sequence.  More precisely, we obtain a solution to \eqref{eq1}-\eqref{bd}-\eqref{inv} as the limit of the sequence of solutions 
$(\mathbf{v}^{n},\boldsymbol{\rho}^{n})_{n\in\N}$ starting from some well chosen initial sequence $(\mathbf{v}^{0},\boldsymbol{\rho}^{0})$, and where for each $n\geq1$, the sequences $\mathbf{v}^n=(v_j^n)_{j\in\Z}$ and $\boldsymbol{\rho}^n=(\rho_j^n)_{j\in\Z}$ are solutions to the following problem:
\begin{equation}\label{eqn}
\forall t>0,\quad \left\{
	\begin{split}
	\partial_{t} v_{j}^{n}(t, x)&=d\partial_{x}^2v_{j}^{n}(t,x),\quad x\in(0,1), \\
		(\rho_{j}^{ n})'(t)
		&=f(\rho_{j}^{n}(t))+\alpha(v_{j}^{n-1}(t,0)+v_{j-1}^{n-1}(t,1))-2\beta\rho_{j}^{n}(t),
	\end{split}	
	\right.
\end{equation}
with Robin boundary conditions
\begin{equation}\label{bdn}
\forall t>0,\quad \left\{
	\begin{split}
		-d\partial_{x}v_{j}^{n}(t,0)+\alpha v_{j}^{n}(t,0)&=\beta \rho_{j}^{n}(t), \\
		d\partial_{x}v_{j}^{n}(t,1)+\alpha v_{j}^{n}(t,1)&=\beta \rho_{j+1}^{n}(t),
	\end{split}\right.
\end{equation}
and initial datum
\begin{equation}\label{Invn}
	\begin{cases}
		v_{j}^{n}(0,x)=h_{j}(x),\quad x\in[0,1], \\
		\rho_{j}^{n}(0)=\Lambda_{j},
	\end{cases}
\end{equation}
for all $j\in\Z$.  We divide the proof into several steps.

\paragraph{Step 1: initialization.} We first set
\bqs
0< \Theta :=  \max\left\{\dfrac{\beta}{\alpha},\|\mathbf{h}\|_{\infty} , \dfrac{\beta}{\alpha}\boldsymbol{\|\Lambda}\|_{\ell^{\infty}(\Z)} \right\}.
\eqs
For each $j\in\Z$, we consider the following Cauchy problem
\begin{equation}\label{rhon=0}\left\{
\begin{split}
\rho_{j}'(t)
	&=f(\rho_{j}(t))+2\alpha \Theta -2\beta\rho_{j}(t), \quad t>0, \\
	\rho_{j}(0)&=\Lambda_{j}.
\end{split}\right.
\end{equation}
Since $f$ is globally Lipschitz continuous on $\R$ thanks to our assumption (H), for each $j\in\Z$, there exists a unique global solution, denoted by $\rho^0_j$ defined on $\R_+$ which satisfies $\rho^0_j\in\mathscr{C}^1(\R_+,\R)$. By definition of $\Theta$, we further have that $0\leq \rho_j^0(t) \leq \alpha\Theta/\beta$ for all $t\geq0$ and $j\in\Z$. For each $j\in\Z$, given $\rho^0_j$, we can then consider
 \begin{equation}\label{vn=0}
\forall t>0, \quad\left\{
	\begin{split}
	\partial_{t} v_{j}(t, x)&=d\partial_{x}^2v_{j}(t,x),\quad x\in(0,1), \\
		-d\partial_{x}v_{j}(t,0)+\alpha v_{j}(t,0)&=\beta \rho_{j}^0(t), \\
		d\partial_{x}v_{j}(t,1)+\alpha v_{j}(t,1)&=\beta \rho_{j+1}^0(t),
	\end{split}\right. 
\end{equation}
with initial data
\bqs
v_{j}(0,x)=h_{j}(x),\quad x\in[0,1].
\eqs
Since $\rho_{j}^0\in \mathscr{C}^1(\R_+,\R)$ for all $j\in\Z$, there exists a unique classical solution $v_j^0$, that is
\bqs
v_j^0\in \mathscr{C}^{0}(\R_+\times[0,1],\R),~ \partial_t v_j^0, \partial_x^2v_j^0\in\mathscr{C}^{0}((0,+\infty)\times(0,1),\R), \text{ and } \partial_x v_j^0\in \mathscr{C}^{0}((0,+\infty)\times[0,1],\R).
\eqs
From the comparison principle applied to \eqref{vn=0}, we further obtain that
\bqs
\forall t\geq0, \, x\in[0,1], \quad 0\leq v_j^0(t,x) \leq \Theta.
\eqs

\paragraph{Step 2: solvability of \eqref{eqn}-\eqref{bdn}-\eqref{Invn} for each $n\geq1$.} We use induction to show that \eqref{eqn}-\eqref{bdn}-\eqref{Invn} has a unique solution. For $n=1$ and each $j\in\Z$, we consider the following Cauchy problem
\begin{equation}\label{rhon=1}\left\{
\begin{split}
\rho_{j}'(t)
	&=f(\rho_{j}(t))+\alpha(v_{j}^{0}(t,0)+v_{j-1}^{0}(t,1))-2\beta\rho_{j}(t), \quad t>0, \\
	\rho_{j}(0)&=\Lambda_{j}.
\end{split}\right.
\end{equation}
Since $f$ is globally Lipschitz continuous, and since $t\longmapsto v_{j}^{0}(t,0)$ and $t\longmapsto v_{j}^{0}(t,1)$ are continuous and uniformly bounded on $\R_+$, the Cauchy-Lipschitz theorem ensures the existence of a unique global solution to \eqref{rhon=1} denoted by $\rho_{j}^1$ and which satisfies $\rho_{j}^1\in \mathscr{C}^1(\R_+,\R)$. Since for all $t\geq 0$ and $j\in\Z$, $0\leq v_{j}^{0}(t,0)+v_{j}^{0}(t,1)\leq 2\Theta$, by comparison, $\rho_j^1$ is a subsolution to \eqref{rhon=0}, and we deduce that $0\leq \rho_j^1(t)\leq \rho_j^0(t) \leq \alpha \Theta/\beta$ for all $t\geq0$. Next, we for each $j\in\Z$, we look at the following evolutionary problem 
\begin{equation}\label{vn=1}
\forall t>0, \quad \left\{
	\begin{split}
	\partial_{t} v_{j}(t, x)&=d\partial_{x}^2v_{j}(t,x),\quad x\in(0,1), \\
		-d\partial_{x}v_{j}(t,0)+\alpha v_{j}(t,0)&=\beta \rho_{j}^1(t), \\
		d\partial_{x}v_{j}(t,1)+\alpha v_{j}(t,1)&=\beta \rho_{j+1}^1(t),
	\end{split}\right. 
\end{equation}
with initial data
\bqs
		v_{j}(0,x)=h_{j}(x),\quad x\in[0,1].
\eqs
Since $\rho_{j}^1\in \mathscr{C}^1(\R_+,\R)$ for all $j\in\Z$, there exists a unique classical solution $v_j^1$, that is
\bqs
v_j^1\in \mathscr{C}^{0}(\R_+\times[0,1],\R),~ \partial_t v_j^1, \partial_x^2v_j^1\in\mathscr{C}^{0}((0,+\infty)\times(0,1),\R), \text{ and } \partial_x v_j^1\in \mathscr{C}^{0}((0,+\infty)\times[0,1],\R).
\eqs
Finally, since $\rho_j^1(t)\leq \rho_j^0(t)$, the comparison principle shows that $v_j^1$ is a subsolution to \eqref{vn=0}, and thus we get that $ v_j^1(t,x)\leq v_j^0(t,x)$ for all $t>0$ and $x\in[0,\ell]$. On other hand, given that $\rho_j^1(t)\geq0$ and $h_{j}(x)\geq0$ for all $x\in[0,1]$ the comparison principle applied to \eqref{vn=1} shows that $v_j^1(t,x)\geq0$ for all $t>0$ and $x\in[0,\ell]$. As a consequence, we have proved that
\bqs
\forall t\geq0, \, x\in[0,1], 0\leq v_j^1(t,x) \leq v_j^0(t,x)\leq \Theta.
\eqs
We can then apply an induction argument to obtain the existence of a nonincreasing sequence of unique classical global solution to system \eqref{eqn}-\eqref{bdn}-\eqref{Invn}. More precisely, for each $n\geq1$, there is exists a unique couple of solution $\rho_{j}^n\in \mathcal{C}^1(\R_+,\R)$ and $v_j^n$ having the following regularity 
\bqs
v_j^n\in \mathscr{C}^{0}(\R_+\times[0,1],\R),~ \partial_t v_j^n, \partial_x^2v_j^n\in\mathscr{C}^{0}((0,+\infty)\times(0,1),\R), \text{ and } \partial_x v_j^n\in \mathscr{C}^{0}((0,\infty)\times[0,1],\R),
\eqs
for each $j\in\Z$. For each $n\geq1$, we also have the following uniform bounds
\bqq\label{bounded1}
0 \leq  \rho_j^{n}(t) \leq \frac{\alpha \Theta}{\beta}= \max\left\{\|\boldsymbol{\Lambda}\|_{\ell^{\infty}(\Z)},1,\dfrac{\alpha}{\beta}\|\mathbf{h}\|_{\infty} \right\},
\eqq
and
\bqq\label{bounded2}
 0 \leq  v_j^{n}(t,x)\leq \Theta =  \max\left\{\dfrac{\beta}{\alpha},\|\mathbf{h}\|_{\infty} , \dfrac{\beta}{\alpha}\|\boldsymbol{\Lambda}\|_{\ell^{\infty}(\Z)} \right\}, \quad x\in[0,1],
\eqq
for all $j\in\Z$ and $t\geq0$. Let us further remark that $(\mathbf{v}^n)_{n\in\N}$ and $(\boldsymbol{\rho}^n)_{n\in\N}$ are decreasing sequences, that is
\bqs
\rho_j^{n+1}(t)\leq \rho_j^{n}(t) \quad \text{ and } \quad v_j^{n+1}(t,x) \leq v_j^{n}(t,x), \quad x\in[0,1],
\eqs
for all $n\in\N$, $j\in\Z$ and $t\geq0$. Finally, let us observe that thanks to our assumption on $f$ and the above uniform bound \eqref{bounded1}, one gets the following uniform bound for the time derivative of $\rho_{j}^{ n}$, namely
\bqs
\forall t>0, \quad \left|(\rho_{j}^{ n})'(t)\right|\leq \left(\frac{f'(0)}{\beta}+3\right)\max\left\{\beta, \alpha \|\mathbf{h}\|_{\infty} , \beta\|\boldsymbol{\Lambda}\|_{\ell^{\infty}(\Z)} \right\},
\eqs
for any $n\geq1$ and $j\in\Z$.

\paragraph{Step 3: existence of a solution.} Let $T>0$ be fixed and $(\mathbf{v}^n,\boldsymbol{\rho}^n)_{n\in\N}$ be the monotone sequence of solutions of \eqref{eqn}-\eqref{bdn}-\eqref{Invn} constructed in the previous step. We already know that for each $n\geq1$ and $j\in\Z$ the function $\rho_j^n$ is globally Lipschitz continuous. As a consequence, since each $v_j^n$ is a solution of
\bqs
\partial_t v_j^n(t,x)=d\partial_x^2 v_j^n(t,x),\quad t>0,~x\in(0,1),
\eqs
with Robin boundary conditions
\bqs\left\{
	\begin{split}
		-d\partial_{x}v_{j}^{n}(t,0)+\alpha v_{j}^{n}(t,0)&=\beta \rho_{j}^{n}(t), \\
		d\partial_{x}v_{j}^{n}(t,1)+\alpha v_{j}^{n}(t,1)&=\beta \rho_{j+1}^{n}(t),
	\end{split}\right.\quad	\forall t>0,
\eqs
and initial datum
\bqs
v_{j}^{n}(0,x)=h_{j}(x),\quad x\in[0,1],
\eqs
we have, by standard parabolic estimates for the heat equation on bounded domain with Robin boundary conditions \cite{LSU}, that there exists $0<\nu<1$ such that for any $\tau\in(0,T)$ and 
\begin{align*}
\forall n\geq1,~j\in\Z, \quad \| v_j^n\|_{\mathscr{C}^{0,\nu}([\tau,T]\times[0,1])}&+\| \partial_x v_j^n\|_{\mathscr{C}^{0,\nu}([\tau,T]\times[0,1])}\\
 &\leq C \left( \|\rho_j^n\|_{L^\infty([0,T+1])} + \|\rho_{j+1}^n\|_{L^\infty([0,T+1])} + \|v_j^n\|_{L^\infty([0,T+1]\times[0,1])} \right) \\
&\leq C\left( \boldsymbol{\|\Lambda}\|_{\ell^{\infty}(\Z)}+\|\mathbf{h}\|_{\infty} \right),
\end{align*}
where the constant $C>0$ only depends on $\nu,\tau,T,d,\alpha,\beta$ . Then, by Schauder estimates \cite{LSU}, we also have for all $n\geq1$ and $j\in\Z$
\begin{align*}
\| \partial_t v_j^n\|_{\mathscr{C}^{0,\nu}([\tau,T]\times[0,1])}&+\| \partial_x^2 v_j^n\|_{\mathscr{C}^{0,\nu}([\tau,T]\times[0,1])}\\
 &\leq C' \left( \|\rho_j^n\|_{\mathscr{C}^{0,\nu}([\tau/2,T])} + \|\rho_{j+1}^n\|_{\mathscr{C}^{0,\nu}([\tau/2,T])} + \|v_j^n\|_{L^\infty([0,T+1]\times[0,1])} \right) \\
&\leq C'\left( \boldsymbol{\|\Lambda}\|_{\ell^{\infty}(\Z)}+\|\mathbf{h}\|_{\infty} \right),
\end{align*}
for some constant $C'>0$ independent of $n$ and $j$. Finally, returning to the equation satisfied by $\rho_j^n$, we also deduce that $(\rho_{j}^{ n})' \in \mathscr{C}^{0,\nu}([\tau,T])$ with the following uniform estimate
\bqs
\|(\rho_{j}^{ n})'\|_{\mathscr{C}^{0,\nu}([\tau,T])} \leq (f'(0)+2\beta)\|{\rho_{j}^{ n}}\|_{\mathscr{C}^{0,\nu}([\tau,T])}+2\alpha \| v_j^n\|_{\mathscr{C}^{0,\nu}([\tau,T]\times[0,1])} \leq C''(\boldsymbol{\|\Lambda}\|_{\ell^{\infty}(\Z)}+\|\mathbf{h}\|_{\infty}),
\eqs
for some $C''>0$ independent of $n$ and $j$.

As a consequence, for any $(N,M)\in\Z^2$ such that $N<M$, the sequence 
\bqs
((v_j^n)_{j=N,\dots,M},(\rho_j^n)_{j=N,\dots,M})_{n\in\N},
\eqs
 together with its respective time derivatives and space derivatives for $v_j^n$ up to order 2, is uniformly bounded  in $\mathscr{C}^{0,\nu}$ norm on the compact set $[\tau,T]\times[0,1]$. By Arzela-Ascoli's theorem, up to a subsequence, there exists a limit sequence $(\mathbf{v},\boldsymbol{\rho})=(v_j,\rho_j)_{j\in\Z}$ such that $(\mathbf{v}^n,\boldsymbol{\rho}^n)$ converges to $(\mathbf{v},\boldsymbol{\rho})$ as $n\rightarrow+\infty$ on any compact of $(0,+\infty)\times[0,1]\times\Z$, but also its respective time derivative and space derivatives (up to order 2). Thanks to the monotonicity of $(\mathbf{v}^n)_{n\in\N}$, we also get that $\mathbf{v}^{n-1}$ converges to $\mathbf{v}$ as $n\rightarrow+\infty$ on any compact of $(0,+\infty)\times[0,1]\times\Z$ along the same subsequence such that one can pass into the limit in \eqref{eqn} and recover \eqref{eq1}.

From Proposition~\ref{proprep}, one has
\bqs
\begin{split}
v_j^n(t,x)&=\int_{0}^1 \K(t,x-y)h(y)\d y+\int_0^t \left[ \K(t-s,x-1)\rho_{j+1}^n(s)+\K(t-s,x)\rho_{j}^n(s)\right]\d s \\
&~~~+\int_0^t \left[ -\alpha \K(t-s,x-1) +d \partial_x\K(t-s,x-1)\right] v_j^n(s,1)\d s\\
&~~~-\int_0^t \left[ \alpha \K(t-s,x) +d \partial_x\K(t-s,x)\right] v_j^n(s,0)\d s,
\end{split}
\eqs
where $\K(t,x):=\frac{1}{\sqrt{4\pi dt}}\exp\left(-\frac{x^2}{4dt}\right)$, and passing to the limit as $n\rightarrow+\infty$ for $t\in[\tau,T]$, $x\in[0,1]$ and $j\in\llbracket N,M\rrbracket$ for $N<M$, we end up with
\bqs
\begin{split}
v_j(t,x)&=\int_{0}^1 \K(t,x-y)h_j(y)\d y+\int_0^t \left[ \K(t-s,x-1)\rho_{j+1}(s)+\K(t-s,x)\rho_{j}(s)\right]\d s \\
&~~~+\int_0^t \left[ -\alpha \K(t-s,x-1) +d \partial_x\K(t-s,x-1)\right] v_j(s,1)\d s\\
&~~~-\int_0^t \left[ \alpha \K(t-s,x) +d \partial_x\K(t-s,x)\right] v_j(s,0)\d s.
\end{split}
\eqs
Taking $t=\tau\rightarrow0$, we recover
\bqs
v_j(t,x)\underset{t\rightarrow0}{\longrightarrow}h_j(x), \quad x\in[0,1].
\eqs
Integrating the second equation in \eqref{eqn} from $0$ to $t$, we get that
\bqs
\rho_{j}^n(t)=\Lambda_{j}+\int_{0}^{t}(f(\rho_{j}^n(s))-2\beta\rho_{j}^n(s))\ds+\alpha\int_{0}^{t}(v_{j}^{n-1}(s,0)+v_{j-1}^{n-1}(s,1))\ds,
\eqs
and passing to the limit as $n\rightarrow+\infty$, we get
\bqs
\rho_{j}(t)=\Lambda_{j}+\int_{0}^{t}(f(\rho_{j}(s))-2\beta\rho_{j}(s))\ds+\alpha\int_{0}^{t}(v_{j}(s,0)+v_{j-1}(s,1))\ds,
\eqs
from which we also recover that
\bqs
\rho_{j}(t)\underset{t\rightarrow0}{\longrightarrow}\Lambda_{j}.
\eqs
Finally, thanks to the regularity of $\partial_xv_j^n$ up to the boundary at $x=0$ and $x=1$, we can also pass to the limit as $n\rightarrow+\infty$ in the boundary condition. As a consequence, $(\mathbf{v},\boldsymbol{\rho})$ is a classical solution to \eqref{eq1}-\eqref{bd}-\eqref{inv}. By uniqueness of the problem \eqref{eq1}-\eqref{bd}-\eqref{inv}, we remark that the convergence of $(\mathbf{v}^n,\boldsymbol{\rho}^n)$ towards $(\mathbf{v},\boldsymbol{\rho})$ holds for all $n$, and not only up to a subsequence. Finally, the a priori bounds \eqref{bounded1} and \eqref{bounded2} gives that the solution $(\mathbf{v},\boldsymbol{\rho})$ satisfies 
\bqs
0 \leq  \rho_j(t) \leq \max\left\{\|\boldsymbol{\Lambda}\|_{\ell^{\infty}(\Z)},1,\dfrac{\alpha}{\beta}\|\mathbf{h}\|_{\infty} \right\},
\eqs
and
\bqs
 0 \leq  v_j(t,x)\leq  \max\left\{\dfrac{\beta}{\alpha},\|\mathbf{h}\|_{\infty} , \dfrac{\beta}{\alpha}\|\boldsymbol{\Lambda}\|_{\ell^{\infty}(\Z)} \right\}, \quad x\in[0,1],
\eqs
for all $j\in\Z$ and $t\geq0$. Since $(\mathbf{h},\boldsymbol{\Lambda})\not\equiv(0,0)$, the comparison principle from Proposition~\ref{propCPcauchy} ensures that
\bqs
\forall t>0,~j\in\Z,\quad 0< v_{j}(t,x),~ x\in[0,1], \text{ and } 0< \rho_{j}(t).
\eqs
This concludes the proof of Theorem~\ref{thmcauchy}.

\section{Long time behavior}\label{secLTB}

We now turn to the study of the long time behavior of \eqref{eq1}--\eqref{bd}. 
\begin{theorem}\label{thm-sta}
The unique non-negative, bounded stationary solutions for equation \eqref{eq1}-\eqref{bd} are $(v_{j}^{\infty},\rho_{j}^{\infty})_{j\in\Z}\equiv(0,0)$ and $(v_{j}^{\infty},\rho_{j}^{\infty})_{j\in\Z}\equiv	(\frac{\beta}{\alpha},1)$.
\end{theorem}

\begin{Proof}[ of Theorem \ref{thm-sta}.]Bounded nonnegative stationary solutions of system \eqref{eq1}-\eqref{bd} are solutions to
\begin{equation}\label{seq1}
	\begin{cases}
		0=dv^{\prime\prime}_{j}(x),\quad x\in(0,1),j\in\Z,\\
		0=f(\rho_{j})+\alpha(v_{j}(0)+v_{j-1}(1))-2\beta\rho_{j},\quad j\in\Z,
	\end{cases}
\end{equation}
together with the boundary conditions
\begin{equation}\label{sbd}
	\begin{cases}
		-dv^{\prime}_{j}(0)+\alpha v_{j}(0)=\beta \rho_{j},\quad j\in\Z,\\
		dv^{\prime}_{j}(1)+\alpha v_{j}(1)=\beta \rho_{j+1},\quad j\in\Z.
	\end{cases}
\end{equation}

	It follows from the $v_{j}$-equation of \eqref{seq1} that there exist two sequences $(a_j)_{j\in\Z}$ and $(b_j)_{j\in\Z}$ of real numbers such that 
	\begin{equation}\label{v_j}
		v_{j}(x)=a_jx+b_j, \quad x\in(0,1),\quad j\in\Z.
	\end{equation}
	Then using \eqref{v_j}, we have that
	\begin{equation}\label{v'}
		v^{\prime}_{j}(0)=v^{\prime}_{j}(1)=v_{j}(1)-v_{j}(0).
	\end{equation}
	Substituting \eqref{v'} into \eqref{sbd} we get that
	\begin{equation*}
		\begin{cases}
			(d+\alpha)v_{j}(0)-d v_{j}(1)=\beta\ \rho_{j},\quad j\in\Z,\\
			-dv_{j}(0)+(d+\alpha) v_{j}(1)=\beta \rho_{j+1},\quad j\in\Z.
		\end{cases}
	\end{equation*}
	Solving the above system, we obtain that
	\begin{equation}\label{vj0}
		v_{j}(0)=\dfrac{\beta d}{\alpha(2d+\alpha )}\left(\dfrac{d+\alpha }{d}\rho_{j}+\rho_{j+1}\right),\quad \forall j\in\Z,
	\end{equation}
	and
	\begin{equation}\label{vjl}
		v_{j}(1)=\dfrac{\beta d}{\alpha(2d+\alpha )}\left(\rho_{j}+\dfrac{d+\alpha }{d}\rho_{j+1}\right),\quad \forall j\in\Z.
	\end{equation}
By applying the second equation of the system \eqref{seq1}, and combining it with \eqref{vj0} and \eqref{vjl}, we derive that
	\begin{equation}\label{srho}
		\frac{\beta d}{2d+\alpha}(\rho_{j-1}-2\rho_{j}+\rho_{j+1})+f(\rho_{j})=0,\quad\forall j\in\Z.
	\end{equation}
	
As a consequence, the existence and uniqueness of bounded nonnegative stationary solutions to system \eqref{seq1}-\eqref{sbd} is equivalent to the existence and uniqueness of bounded nonnegative stationary solutions to \eqref{srho}. 

Since $f(0)=f(1)=0$, it is clear that $\rho_{j}\equiv0,\forall j\in\Z$ and $\rho_{j} \equiv1,\forall j\in\Z$ are always  solutions of \eqref{srho}. Let us prove that they are actually the only bounded nonnegative stationary solutions of \eqref{srho}.  

Let $\boldsymbol{\rho}=(\rho_j)_{j\in\Z} \neq 0$ be a non zero, bounded, nonnegative stationary solution to \eqref{srho}. Then, necessarily, one has $\rho_j>0$ for all $j\in\Z$. Indeed, if $\rho_{j_0}=0$ for some $j_0\in\Z$, then equation \eqref{srho} implies that $\rho_{j_0+1}=\rho_{j_0-1}=0$, and by induction, $\rho_j=0$ for all $j\in\Z$. So from now on, we assume that $\boldsymbol{\rho}$ satisfies $\rho_j>0$ for all $j\in\Z$. 

We let $N\geq1$ be an integer which satisfies 
\bqq\label{eqcondN}
\dfrac{2\beta d}{2d+\alpha}\left(1-\cos\left(\frac{\pi}{2(N+1)}\right)\right) < f'(0).
\eqq
Consider the eigenvalue problem
\bqs
-(\rho_{j-1}-2\rho_{j}+\rho_{j+1})=\mu \rho_j, \quad j=-N,\dots,N,
\eqs
and $\rho_j=0$ for all $j\leq-N-1$ and $j\geq N+1$. One easily finds that the eigenvalues are given by
\bqs
\mu_p=2\left(1-\cos\left(\frac{p\pi}{2(N+1)}\right)\right), \quad p=1,\dots,2N+1,
\eqs
with corresponding eigenfunctions $\boldsymbol{\phi}^p=(\phi_j^p)_{j\in\Z}$ defined as
\bqs
\phi_j^p=\left\{\begin{split}
\cos\left(\frac{pj\pi}{2(N+1)}\right), & \quad  j=-N,\dots,N,\\
0, & \quad \text{ otherwise.}
\end{split}\right.
\eqs
We consider the principal eigenfunction $\boldsymbol{\phi}^1$ with eigenvalue $\mu_1>0$. Thanks to condition \eqref{eqcondN} on $N$ which ensures that 
\bqs
\dfrac{\beta d}{2d+\alpha}\mu_1 < f'(0),
\eqs
and thanks to the regularity assumption on the function $f$, there exists $\epsilon_0>0$ such that for all $\epsilon\in(0,\epsilon_0]$ one has
\bqs
- \mathscr{L}(\epsilon \boldsymbol{\phi}^1)_j < f(\epsilon \phi^1_j), \quad j=-N,\dots,N
\eqs
where the operator $\mathscr{L}:\ell^{\infty}(\Z)\rightarrow\ell^{\infty}(\Z)$ is defined
as
\bqs
\mathscr{L}(\boldsymbol{\rho})_j = \dfrac{\beta d}{2d+\alpha}(\rho_{j-1}-2\rho_{j}+\rho_{j+1}), \quad j\in\Z,
\eqs
for any $\boldsymbol{\rho}\in\ell^{\infty}(\Z)$.

Given $\boldsymbol{\rho}=(\rho_j)_{j\in\Z} \neq 0$ a non zero, bounded, positive stationary solution to \eqref{srho}, for any $k\in\Z$, we define the shifted sequence $\widehat{\boldsymbol{\rho}}=(\widehat{\rho}_j)_{j\in\Z}$ as $\widehat{\rho}_j=\rho_{j+k}$ for all $j\in\Z$. By translation invariance of equation \eqref{srho}, we readily have that $\widehat{\boldsymbol{\rho}}$ is also a non zero, bounded, positive stationary solution to \eqref{srho}. With $\epsilon_0>0$ defined above, we claim that
\bqq\label{eqClaimRhoHat}
\widehat{\rho}_j>\epsilon_0 \phi^1_j, \quad j=-N-1,\dots,N+1.
\eqq
If it were not the case, then we can find $\epsilon_*\in(0,\epsilon_0]$ such that $\widehat{\rho}_j \geq \epsilon_* \phi^1_j$ for all $j=-N-1,\dots,N+1$ and $\widehat{\rho}_{j_*} = \epsilon_* \phi^1_{j_*}$ for some $j_*\in \llbracket-N-1,N+1\rrbracket$. Let us first remark that since $\phi^1_{-N-1}=\phi^1_{N+1}=0$ and $\widehat{\rho}_j>0$ for all $j\in\Z$, we necessarily have that $j_*\in \llbracket-N,N\rrbracket$. We next note that $\epsilon_*\boldsymbol{\phi}^1$ is still a sub-solution of \eqref{srho}. We can then apply the strong maximum principle of Proposition~\ref{propSMPdiscret} which shows that $\widehat{\rho}_{j} = \epsilon_* \phi^1_{j}$ for all $j\in \llbracket-N-1,N+1\rrbracket$ which leads to a contradiction since it implies that $0<\widehat{\rho}_{N+1} = \epsilon_* \phi^1_{N+1}=0$. As a conclusion, from the inequality \eqref{eqClaimRhoHat} and using the definition of $\widehat{\boldsymbol{\rho}}$, we have proved that
\bqs
\forall k\in\Z, \quad \rho_{j+k}>\epsilon_0 \phi^1_j, \quad j=-N-1,\dots,N+1,
\eqs
which, evaluated at $j=0$, gives the uniform lower bound:
\bqs
\forall k\in\Z, \quad \rho_{k}>\epsilon_0.
\eqs
As a consequence, we can define
\bqs
0<m:=\inf_{j\in\Z}~\rho_j .
\eqs
Assume that $m<1$. We let $(j_k)_{k\in\N}$ such that
\bqs
\rho_{j_k}\underset{k\rightarrow+\infty}{\longrightarrow}m.
\eqs
For each $j\in\Z$, we denote
\bqs
\widetilde{\rho}_j:=\underset{k\rightarrow+\infty}{\lim}\rho_{j+j_k},
\eqs
and we remark that $\widetilde{\rho}_j$ also satisfies \eqref{srho}. We also note that $\widetilde{\rho}_0=m$ and by construction
\bqs
\widetilde{\rho}_0=m=\inf_{j\in\Z}~\widetilde{\rho}_j.
\eqs
It also satisfies 
\bqs
\dfrac{\beta d}{2d+\alpha }\left(\underbrace{\widetilde{\rho}_{-1}-\hat{\rho}_0}_{\geq0}+\underbrace{\widetilde{\rho}_{1}-\widetilde{\rho}_0}_{\geq0}\right)=-f(\widetilde{\rho}_0)=-f(m)<0,
\eqs
which is impossible. As a consequence, one has $m\geq1$. By a similar argument, this time with $M=\sup_{j\in\Z}~\rho_j>0$, one gets that necessarily $M\leq1$. This implies that $\rho_j=1$ for all $j\in\Z$ and concludes the proof of the theorem.
\end{Proof}

Next, we demonstrate that the positive stationary solution $(v_{j}^{\infty},\rho_{j}^{\infty})_{j\in\Z}\equiv	(\frac{\beta}{\alpha},1)$
is the global attractor of the system \eqref{eq1}--\eqref{bd} starting from nontrivial nonnegative bounded initial condition. More precisely, we shall prove the following result.

\begin{theorem}\label{long}
	Let $(\mathbf{v},\boldsymbol{\rho})$ be the unique global classical solution of \eqref{eq1}-\eqref{bd}-\eqref{inv} starting from a nontrivial nonnegative bounded initial sequence $(\mathbf{h},\boldsymbol{\Lambda}) \in \mathcal{X}^0\times \ell^\infty(\Z)$ satisfying the compatibility condition \eqref{compatibility}. Then,
	$$
	\lim_{t\to+\infty}(v_{j}(t,x),\rho_{j}(t))=\left(\frac{\beta}{\alpha},1\right), \quad \forall x \in[0,1],
	$$
	locally uniformly $j\in\Z$.
\end{theorem}

\begin{Proof}[ of Theorem \ref{long}.] The first part of the proof consists of constructing a nonnegative, compactly supported, stationary subsolution to \eqref{eq1}-\eqref{bd}. Actually, following the proof of Theorem \ref{thm-sta}, there exists $N_0>1$ large enough such that condition \eqref{eqcondN} is satisfied for all $N\geq N_0$. Next, with $N\geq N_0$, we define $\underline{\boldsymbol{\rho}}=\left(\underline{\rho}_j\right)_{j\in\Z}$ as
\bqs
\underline{\rho}_j:=\left\{\begin{split}
\cos\left(\frac{j\pi}{2(N+1)}\right), & \quad  j=-N,\dots,N,\\
0, & \quad \text{ otherwise,}
\end{split}\right.
\eqs
and set
\bqs
\underline{v}_j(x):=a_j x + b_j, \quad x\in[0,1], \quad j\in\Z,
\eqs
with
\bqs
a_j:=\frac{\beta}{2d+\alpha}\left(\underline{\rho}_{j+1}-\underline{\rho}_j\right)\quad \text{ and }\quad b_j:=\dfrac{\beta d}{\alpha(2d+\alpha)}\left(\dfrac{d+\alpha}{d}\underline{\rho}_{j}+\underline{\rho}_{j+1}\right),\quad \forall j\in\Z.
\eqs
By construction, $\underline{\mathbf{v}}=(\underline{v}_j)_{j\in\Z }$ is compactly supported, and we have
\bqs
d \partial_x^2\underline{v}_j(x)=0, \quad x\in(0,1),
\eqs
together with
\bqs
\begin{cases}
		-d\partial_x \underline{v}_{j}(0)+\alpha \underline{v}_{j}(0)=\beta \underline{\rho}_{j},\\
		d\partial_x\underline{v}_{j}(1)+\alpha \underline{v}_{j}(1)=\beta \underline{\rho}_{j+1},
	\end{cases}
\eqs
for all $j\in\Z$. Finally, there exists $\epsilon_0>0$ such that for all $\epsilon\in(0,\epsilon_0]$ one has:
\bqs
\dfrac{\epsilon\beta d}{2d+\alpha}\left(\underline{\rho}_{j-1}-2\underline{\rho}_{j}+\underline{\rho}_{j+1}\right)+f(\epsilon \underline{\rho}_j)>0, \quad, j=1,\dots,N.
\eqs
As a consequence $\left( \epsilon \underline{\mathbf{v}},\epsilon\underline{\boldsymbol{\rho}}\right)$ is a stationary, compactly supported, subsolution for all $N\geq N_0$ and $\epsilon\in(0,\epsilon_0]$.

We now use the method of super- and subsolutions to prove the theorem. First, we consider $\left(\overline{\mathbf{v}},\overline{\boldsymbol{\rho}}\right)$ defined by
	$$
\overline{v}_j(x):=\max\left\{\dfrac{\beta}{\alpha},\|\mathbf{h}\|_{\infty} , \dfrac{\beta}{\alpha}\|\boldsymbol{\Lambda}\|_{\ell^{\infty}(\Z)}\right\}, \quad \overline{\rho}_j=\max\left\{1 , \|\boldsymbol{\Lambda}\|_{\ell^{\infty}(\Z)},\dfrac{\alpha}{\beta}\|\mathbf{h}\|_{\infty}\right\},\quad x\in[0,1],~j\in\Z.
	$$
Let $t\mapsto \left(\widehat{\mathbf{v}}(t),\widehat{\boldsymbol{\rho}}(t)\right)$ be the global solution of \eqref{eq1}-\eqref{bd} with initial condition $ \left(\widehat{\mathbf{v}}(0),\widehat{\boldsymbol{\rho}}(0)\right)=\left(\overline{\mathbf{v}},\overline{\boldsymbol{\rho}}\right)$. It follows from the comparison principle from Proposition~\ref{propCPcauchy} that $t\mapsto\left(\widehat{\mathbf{v}}(t),\widehat{\boldsymbol{\rho}}(t)\right)$ is non increasing in time $t$ and satisfies $\left(\frac{\beta}{\alpha} , 1\right) \leq \left(\widehat{\mathbf{v}}(t),\widehat{\boldsymbol{\rho}}(t)\right)$ for all $t>0$. Thus, owing to Theorem \ref{thm-sta}, as $t \rightarrow+\infty$, it converges locally uniformly in $j$ to the unique positive solution of \eqref{seq1}-\eqref{sbd}, namely $\left(v_{j}^{\infty}, \rho_j^{\infty}\right) \equiv\left(\frac{\beta}{\alpha} , 1\right)$, that is
\bqs
\forall x\in[0,1], \quad \widehat{v}_j(t,x)\underset{t\rightarrow+\infty}{\longrightarrow} \frac{\beta}{\alpha} \quad \text{ and } \quad \widehat{\rho}_j(t)\underset{t\rightarrow+\infty}{\longrightarrow} 1,
\eqs
locally uniformly in $j\in\Z$. Now, let $t\mapsto (\mathbf{v}(t),\boldsymbol{\rho}(t))$ be the solution of \eqref{eq1}-\eqref{bd}-\eqref{inv} starting from the nonnegative, not identically equal to zero, bounded initial datum $(\mathbf{h},\boldsymbol{\Lambda})$. Since $(\mathbf{h},\boldsymbol{\Lambda}) \leq \left(\overline{\mathbf{v}},\overline{\boldsymbol{\rho}}\right)$, we have $(\mathbf{v}(t),\boldsymbol{\rho}(t))\leq \left(\widehat{\mathbf{v}}(t),\widehat{\boldsymbol{\rho}}(t)\right)$ for all $t>0$ and thus
$$
\forall x\in[0,1], \quad \limsup_{t\to+\infty}~(v_j(t, x), \rho_j(t)) \leq\left(\frac{\beta}{\alpha} ,1\right),
$$
locally uniformly in $j \in \Z$. Furthermore, since $0\lneqq (\mathbf{h},\boldsymbol{\Lambda})$, by the comparison principle from Proposition~\ref{propCPcauchy}, we have that $0<(\mathbf{v}(t),\boldsymbol{\rho}(t))$ for all $t>0$. As a consequence, upon  reducing further the size of $\epsilon$, we can always ensure that $\left( \epsilon \underline{\mathbf{v}},\epsilon\underline{\boldsymbol{\rho}}\right) \leq (\mathbf{v}(1),\boldsymbol{\rho}(1))$. We now let $t\mapsto \left(\underline{\mathbf{v}}(t),\underline{\boldsymbol{\rho}}(t)\right)$ be the global solution of \eqref{eq1}-\eqref{bd} with initial condition $\left( \epsilon \underline{\mathbf{v}},\epsilon\underline{\boldsymbol{\rho}}\right)$, which by the comparison principle, is nondecreasing in $t$. As a consequence, it also converges locally uniformly $j$ to $\left(v_{j}^{\infty}, \rho_j^{\infty}\right) \equiv\left(\frac{\beta}{\alpha} , 1\right)$, the unique positive solution of \eqref{seq1}-\eqref{sbd}. Thus, we have that
	$$
\forall x\in[0,1], \quad \left(\frac{\beta}{\alpha},1\right) = \lim_{t\to+\infty} (\underline{v}_j(t, x), \underline{\rho}_j(t))	\leq \liminf_{t\to+\infty} (v_j(t+1, x), \rho_j(t+1)) ,
	$$
locally uniformly in $j$. This completes the proof of Theorem~\ref{long}.
\end{Proof}

\section{Exponential solutions and linear spreading speed}\label{secEXP}

In this section, we fully characterize the linear spreading speed of system \eqref{eq1}-\eqref{bd} by studying pure exponential solutions to the linearized system around the trivial state. We also numerically explore in a second time the dependence of the spreading speed as function of the various parameters of the model.

\subsection{Derivation of the linear spreading speed}

In order to study the spreading properties of system \eqref{eq1}-\eqref{bd}, we  consider the existence of exponential solutions for the linearized problem around the trivial state which writes:
\begin{equation}\label{leq1}
\forall t>0,~j\in\Z,\quad	\begin{cases}
		\partial_{t} v_{j}(t, x)=d\partial_x^2v_{j}(t,x),\quad x\in(0,1),\\
		\rho_{j}^{\prime}(t)=f^{\prime}(0)\rho_{j}(t)+\alpha(v_{j}(t,0)+v_{j-1}(t,1))-2\beta\rho_{j}(t),
	\end{cases}
\end{equation}
with the boundary condition
\begin{equation}\label{lbd}
\forall t>0,~j\in\Z,\quad	\begin{cases}
		-d\partial_{x}v_{j}(t,0)+\alpha v_{j}(t,0)=\beta \rho_{j}(t),\\
		d\partial_{x}v_{j}(t,1)+\alpha v_{j}(t,1)=\beta \rho_{j+1}(t).
	\end{cases}
\end{equation}

We will be looking for solutions of the form
\begin{equation}\label{exp}
	(v_{j}(t,x),\rho_{j}(t))=\left(\ee^{-\mu(j-ct)}V(x),\ee^{-\mu(j-ct)}\right),
\end{equation}
where 
\begin{equation*}\label{V}
	V(x)=a\cosh\left(\sqrt{\frac{\lambda}{d}} x\right)+b\sinh\left(\sqrt{\frac{\lambda}{d}} x\right), \quad\forall x\in[0,1],
\end{equation*}
for some $\lambda>0$, $\mu>0$, $c>0$ and $(a,b)\in\R^2$ that will be determined later. We substitute ansatz \eqref{exp} into \eqref{leq1}, \eqref{lbd} and get that
\begin{equation}\label{expeq}
\begin{cases}
\mu c=\lambda,\\
\mu c=f^{\prime}(0)+\alpha(V(0)+\ee^{\mu}V(1))-2\beta,\\
-dV^{\prime}(0)+\alpha V(0)=\beta,\\
dV^{\prime}(1)+\alpha V(1)=\beta \ee^{-\mu}.
\end{cases}
\end{equation}
We first express $(a,b)$ as a function of $(V(0), V(1))$, that is 
$$
a=V(0),\quad b=\frac{V(1)}{\sinh\left(\sqrt{\frac{\lambda}{d}}\right)}-\frac{V(0)}{\tanh\left(\sqrt{\frac{\lambda}{d}}\right)}.
$$
We also deduce that
\begin{align*}
V^{\prime}(0)&=\sqrt{\frac{\lambda}{d}}\left[-\frac{1}{\tanh\left(\sqrt{\frac{\lambda}{d}}\right)}V(0)+\frac{1}{\sinh\left(\sqrt{\frac{\lambda}{d}}\right)}V(1)\right],\\
V^{\prime}(1)&=\sqrt{\frac{\lambda}{d}}\left[-\frac{1}{\sinh\left(\sqrt{\frac{\lambda}{d}}\right)}V(0)+\dfrac{1}{\tanh\left(\sqrt{\frac{\lambda}{d}}\right)}V(1)\right].
\end{align*}
As a consequence, using the Robin type boundary in \eqref{expeq}, we deduce that
\begin{equation}\label{matr}
\left(\begin{matrix}
		\alpha+\dfrac{\sqrt{d\lambda}}{\tanh\left(\sqrt{\frac{\lambda}{d}}\right)} & -\dfrac{\sqrt{d\lambda}}{\sinh\left(\sqrt{\frac{\lambda}{d}}\right)}\\
		-\dfrac{\sqrt{d\lambda}}{\sinh\left(\sqrt{\frac{\lambda}{d}}\right)}   & \alpha+\dfrac{\sqrt{d\lambda}}{\tanh\left(\sqrt{\frac{\lambda}{d}}\right)}\\
	\end{matrix}\right)
\begin{pmatrix}
	V(0)\\
	V(1)
\end{pmatrix}
=\beta
\begin{pmatrix}
	1\\
	\ee^{-\mu}
\end{pmatrix}
.
\end{equation}
Define
$$
\Delta(\lambda):=\alpha^2+d\lambda+\dfrac{2\sqrt{d\lambda}\alpha}{\tanh\left(\sqrt{\frac{\lambda}{d}}\right)},
$$
and let us remark that when $\lambda>0$ we have $\Delta(\lambda)>0$ and $\Delta$ is well-defined up to $\lambda=0$ with $\Delta(0)=\alpha^2+2d\alpha>0$. Thus, we can invert the above system \eqref{matr} and deduce that
$$
V(0)=\dfrac{\beta}{\Delta(\lambda)}\left(\frac{\sqrt{d\lambda}}{\tanh\left(\sqrt{\frac{\lambda}{d}}\right)}+\alpha+\frac{\sqrt{d\lambda}}{\sinh\left(\sqrt{\frac{\lambda}{d}}\right)}\ee^{-\mu}\right),
$$
and
$$
V(1)=\dfrac{\beta}{\Delta(\lambda)}\left(\frac{\sqrt{d\lambda}}{\tanh\left(\sqrt{\frac{\lambda}{d}}\right)}\ee^{-\mu}+\alpha\ee^{-\mu}+\frac{\sqrt{d\lambda}}{\sinh\left(\sqrt{\frac{\lambda}{d}}\right)}\right).
$$

Therefore, we substitute the above two formulas into \eqref{expeq} to obtain 
\begin{equation}\label{muc}
	\begin{cases}
		\mu c=\lambda,\\
		\mu c=f^{\prime}(0)-2\beta+\frac{2\alpha\beta}{\Delta(\lambda)}\left[\alpha+\dfrac{\sqrt{d\lambda}}{\tanh\left(\sqrt{\frac{\lambda}{d}}\right)}+\dfrac{\sqrt{d\lambda}}{\sinh\left(\sqrt{\frac{\lambda}{d}}\right)}\cosh(\mu)\right].
	\end{cases}
\end{equation}

\begin{figure}[t!]
\centering
\includegraphics[width=.45\textwidth]{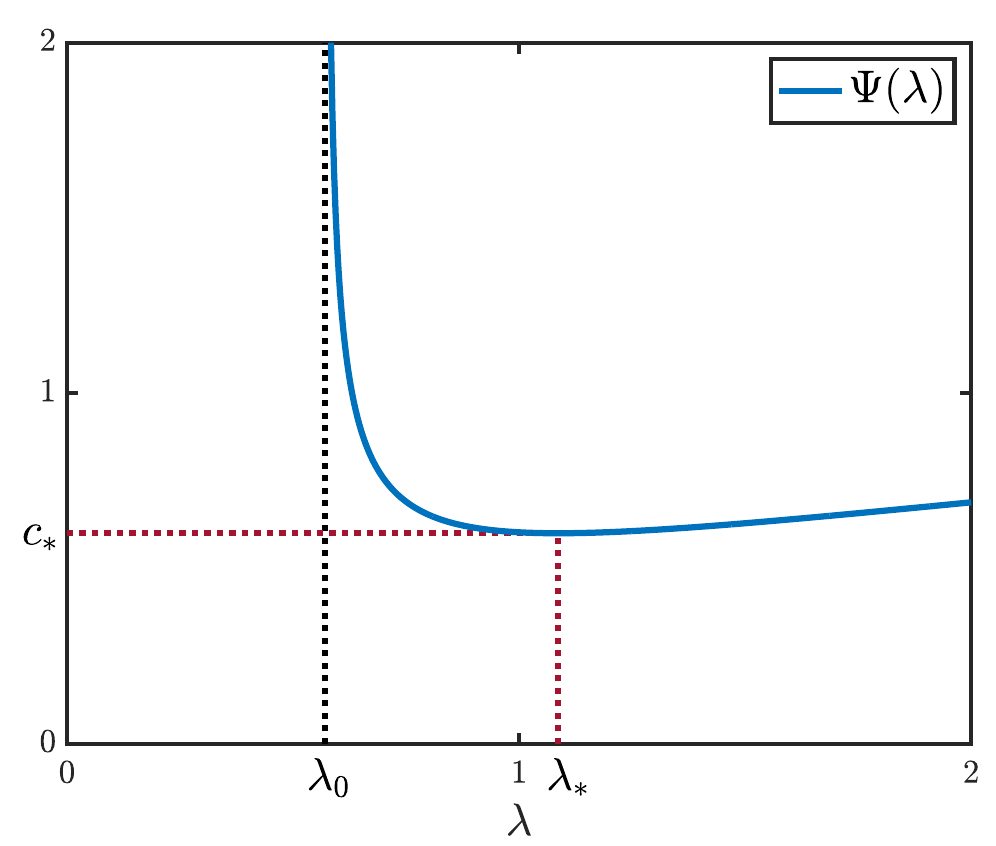}
\caption{Typical representation of the map $\Psi$ on $(\lambda_0,+\infty)$ with a unique global minimum at $\lambda=\lambda_*$. Here, parameters values are set to $(\alpha,\beta,d,f'(0))=(1,1,1,1)$.}
 \label{fig:SpreadingSpeed}
\end{figure}

As a consequence, we have that
\begin{equation*}\label{cosh}
	\cosh\left(\dfrac{\lambda}{c}\right)=y(\lambda),
\end{equation*}
where
\begin{equation}\label{y}
	y(\lambda):= \left[\dfrac{\Delta(\lambda)}{2\alpha\beta}(\lambda+2\beta-f^{\prime}(0))-\left(\alpha+\dfrac{\sqrt{d\lambda}}{\tanh\left(\sqrt{\frac{\lambda}{d}}\right)}\right)\right]\dfrac{\sinh\left(\sqrt{\frac{\lambda}{d}}\right)}{\sqrt{d\lambda}}.
\end{equation}
Solving the above equation, one gets that
\begin{equation}\label{c}
	c=\dfrac{\lambda}{\mu(\lambda)},
\end{equation}
where
\begin{equation*}
	\mu(\lambda):=\ln\left(y(\lambda)+\sqrt{y^2(\lambda)-1}\right),
\end{equation*}
and $y(\lambda)$ is defined in \eqref{y}. Now, to justify all the above computations, one needs to show that $\lambda \mapsto\mu(\lambda)$ is well-defined which is equivalent to proving that $y(\lambda)>1$ for some range of $\lambda$. First of all, one can actually check that $y(\lambda)\to+\infty$ as $\lambda\to+\infty$.  By a direct computation, we also have that
\begin{equation*}
	y(0^+)=1-\dfrac{f^{\prime}(0)}{\beta}\left(\dfrac{\alpha}{2d}+1\right)<1.
\end{equation*}
Thus, there exists $\lambda_0>0$ such that $y(\lambda_0)=1$ and
\begin{equation*}\label{y_gamma}
	y(\lambda)>1,\quad\forall\lambda>\lambda_0.
\end{equation*}
We claim that such a $\lambda_0>0$ is unique. First, we observe that the equality $y(\lambda_0)=1$ is equivalent to
\bqs
\frac{2\alpha\beta}{\Delta(\lambda_0)}\left[\alpha+\dfrac{\sqrt{d\lambda_0}}{\tanh\left(\sqrt{\frac{\lambda_0}{d}}\right)}+\dfrac{\sqrt{d\lambda_0}}{\sinh\left(\sqrt{\frac{\lambda_0}{d}}\right)}\right]=\lambda_0+2\beta-f^{\prime}(0).
\eqs

On the one hand, the map $g:\lambda\mapsto g(\lambda)=\lambda+2\beta-f^{\prime}(0)$ is strictly increasing on $\R_+$ with $g(0)=2\beta-f^{\prime}(0)$ and $g(\lambda)\sim \lambda$ as $\lambda\rightarrow+\infty$. On the other hand, the map $G$ defined as 
\bqs
G:\lambda\mapsto G(\lambda)=\frac{2\alpha\beta}{\Delta(\lambda)}\left[\alpha+\dfrac{\sqrt{d\lambda}}{\tanh\left(\sqrt{\frac{\lambda}{d}}\right)}+\dfrac{\sqrt{d\lambda}}{\sinh\left(\sqrt{\frac{\lambda}{d}}\right)}\right],
\eqs 
is decreasing on $\R_+$ with $G(0)=2\beta$ and $G(\lambda)\sim \frac{2\alpha\beta}{\sqrt{d\lambda}}$ as $\lambda\rightarrow+\infty$. The fact that $G$ is decreasing on $\R_+$ comes from the direct computation:
\bqs
\forall \lambda >0, \quad G'(\lambda)= - \frac{\alpha \beta\left(\lambda+\sqrt{\frac{\lambda}{d}}\sinh\left(\sqrt{\frac{\lambda}{d}}\right) \right)}{\lambda \sinh\left(\sqrt{\frac{\lambda}{d}}\right)^2 \Delta(\lambda)^2}\Theta(\lambda)<0,
\eqs
where
\bqs
\Theta(\lambda):=d\lambda\left(\cosh\left(\sqrt{\frac{\lambda}{d}}\right)-1\right)+\alpha^2\left(\cosh\left(\sqrt{\frac{\lambda}{d}}\right)+1\right)+2\sqrt{d\lambda}\sinh\left(\sqrt{\frac{\lambda}{d}}\right)>0.
\eqs
As a consequence since $g(0)=2\beta-f'(0)<2\beta=G(0)$ and $g$ is strictly increasing with $g(\lambda)\rightarrow+\infty$ as $\lambda\rightarrow+\infty$ while $G$ is strictly decreasing  $G(\lambda)\rightarrow0$ as $\lambda\rightarrow+\infty$, we obtain the existence of a unique $\lambda_0>0$ such that $g(\lambda_0)=G(\lambda_0)$ and $g(\lambda)<G(\lambda)$ for all $\lambda\in[0,\lambda_0)$ together with $g(\lambda)>G(\lambda)$ for all $\lambda>\lambda_0$. This proves the claim.

\begin{figure}[t!]
\centering
\subfigure[$(\beta,d,f'(0))=(1,1,1)$.]{\includegraphics[width=.31\textwidth]{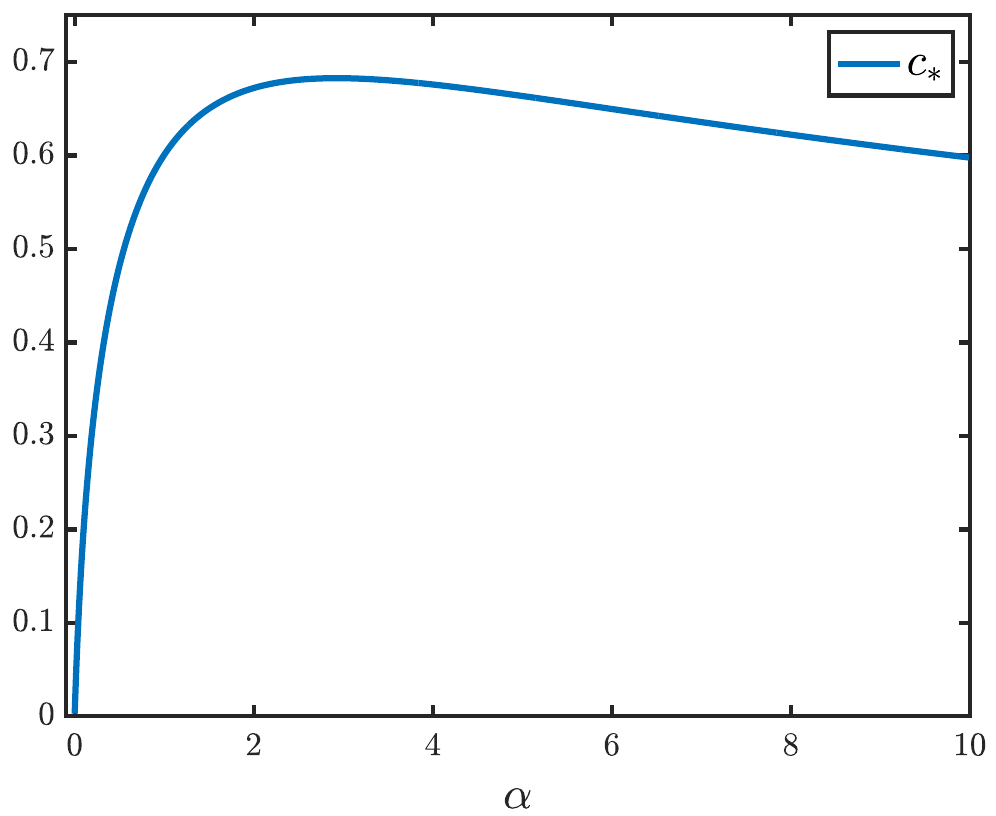}}
\subfigure[$(\alpha,d,f'(0))=(1,1,1)$.]{\includegraphics[width=.31\textwidth]{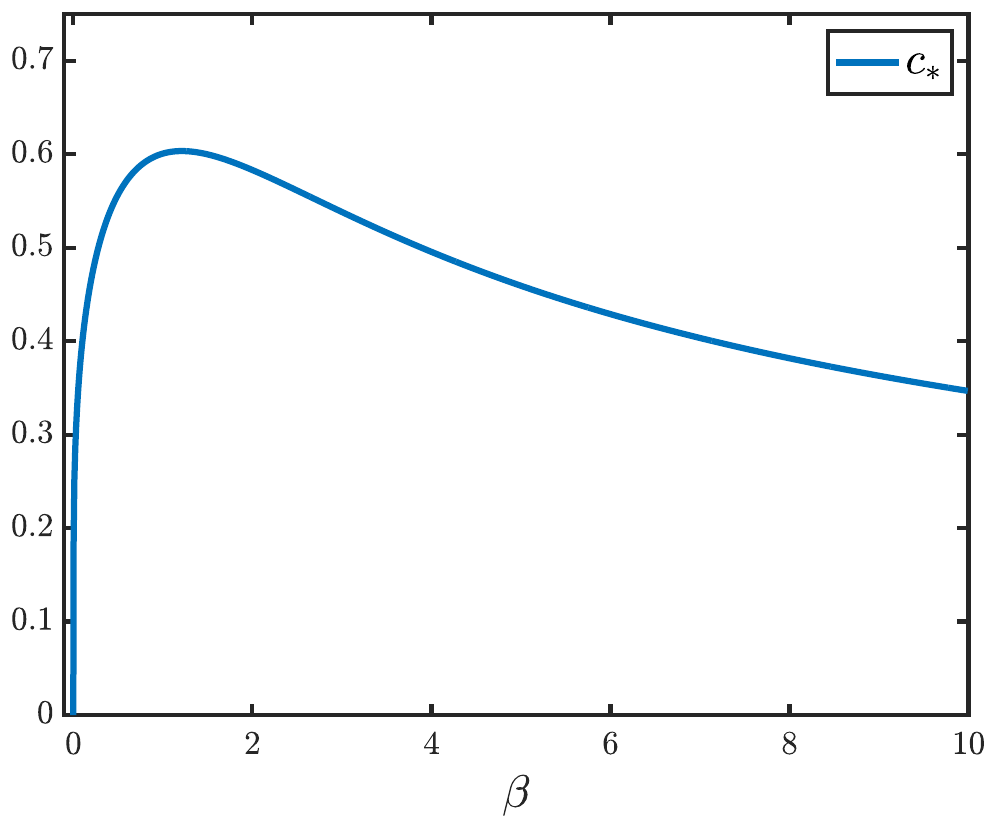}}
\subfigure[$(\alpha,\beta,f'(0))=(1,1,1)$.]{\includegraphics[width=.32\textwidth]{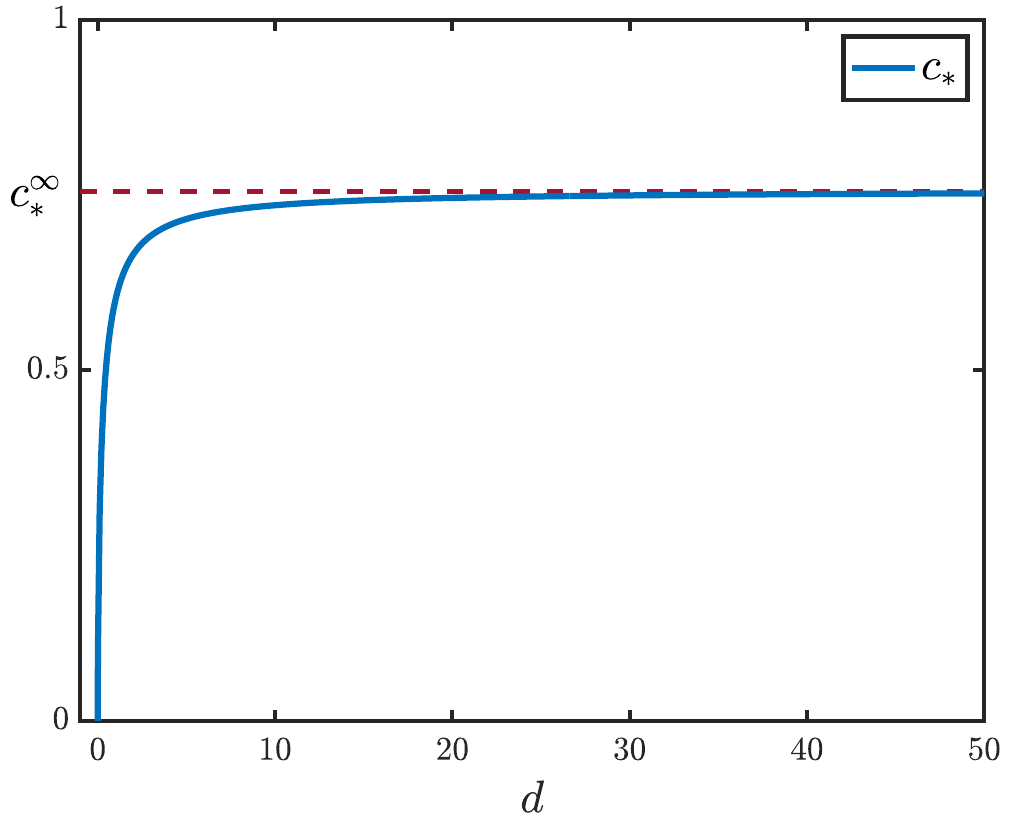}}
\subfigure[$(\alpha,\beta,d)=(1,1,1)$.]{\includegraphics[width=.31\textwidth]{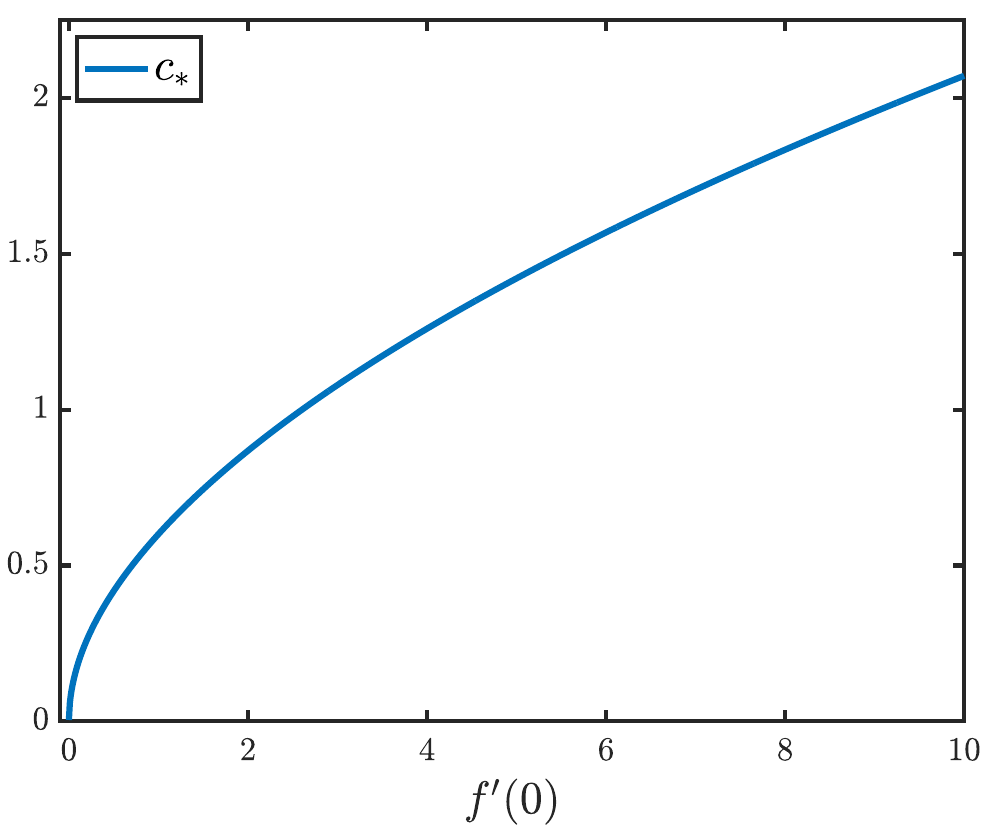}}\hspace{0.5cm}
\subfigure[$(\alpha,\beta,d)=(1,1,1)$.]{\includegraphics[width=.31\textwidth]{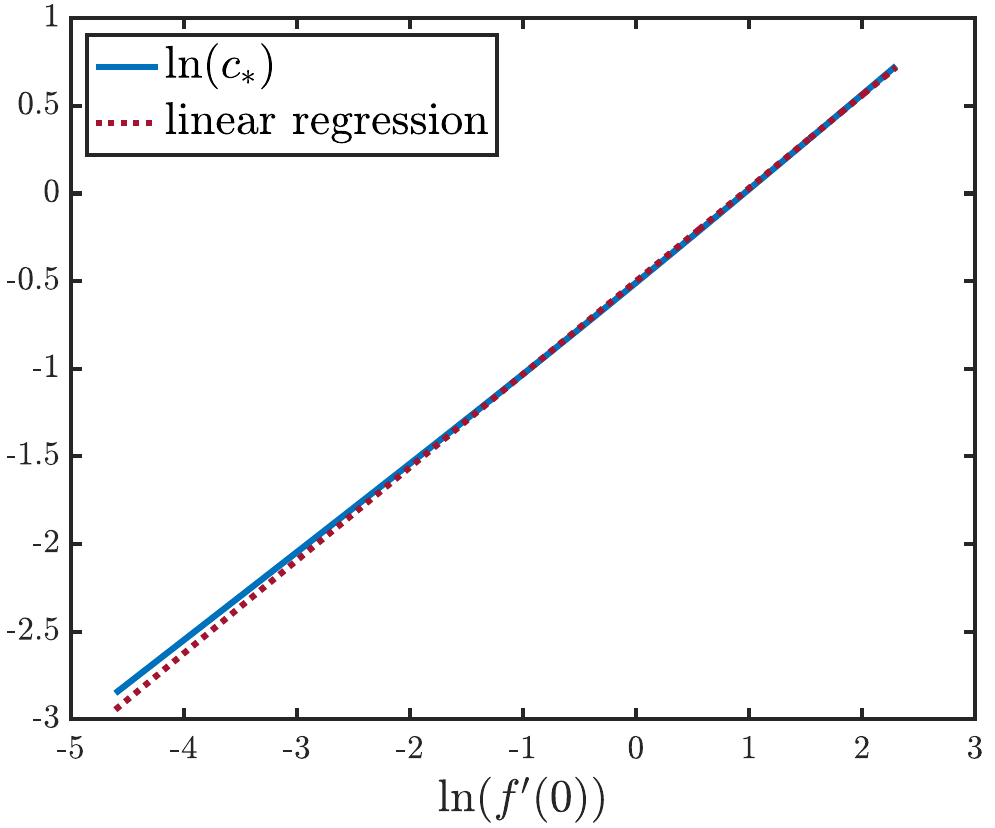}}
\caption{Plots of the linear spreading speed \eqref{eq1} as parameters are varied.}
 \label{fig:Speed}
\end{figure}

The uniqueness of $\lambda_0$ implies that $y(\lambda)<1$ for all $\lambda\in(0,\lambda_0)$ such that $\mu(\lambda)$ is only well-defined for all $\lambda>\lambda_0$. Coming back to \eqref{c}, we define $\Psi:(\lambda_0,+\infty)\to \R_+$ as 
\bqs
\Psi(\lambda):=\frac{\lambda}{\mu(\lambda)}.
\eqs
We readily note that $\Psi(\lambda)\rightarrow+\infty$ as $\lambda\rightarrow\lambda_0^+$ and as $\lambda\rightarrow+\infty$. Since $\Psi$ is smooth on $(\lambda_0,+\infty)$, it achieves a minimum on $(\lambda_0,+\infty)$, and we can define
\begin{equation}\label{c_*}
	c_*:=\min_{\lambda>\lambda_0}~\dfrac{\lambda}{\mu(\lambda)}.
\end{equation}
We present in Figure~\ref{fig:SpreadingSpeed} a typical representation of the map $\Psi$ on $(\lambda_0,+\infty)$. It exhibits a unique global minimum at some $\lambda_*>\lambda_0$, values at which one has
\bqs
c_*=\frac{\lambda_*}{\mu(\lambda_*)}.
\eqs

\subsection{Numerical computations of the linear spreading speed $c_*$}\label{subsecnumspeed}

We numerically computed the linear spreading speed $c_*$ by systematically evaluating the global minimum of the function $\Psi$ as given by formula \eqref{c_*} as a function of the various parameters of the system. We reported the corresponding results in Figure~\ref{fig:Speed}. For the chosen parameter values, variations of the linear spreading speed $c_*$ as a function of $\alpha$ and $\beta$ show a similar pattern with, in both cases, the existence of a maximal spreading speed (see panels (a) and (b) of Figure~\ref{fig:Speed}) at some optimal value of the parameters $\alpha$ or $\beta$. More precisely, as either $\alpha$ or $\beta$ is varied, while all other parameters are kept fixed, the linear spreading speed is first increasing from zero towards a global maximum value and then decreasing. On the other hand, when varying the parameter $d$, we clearly observe a monotone convergence towards a limiting asymptotic value. The limiting value exactly matches the asymptotic spreading speed $c_*^\infty$ defined in formula \eqref{c_11} of Section~\ref{secLDL} below. Finally, as it is the case for spreading speeds for scalar continuous Fisher-KPP equations, we see that the spreading speed $c_*$ is a strictly monotone function of the parameter $f'(0)$, and we conjecture that $c_*$ is proportional to $\sqrt{f'(0)}$. This is numerically confirmed (see panel (e) of Figure~\ref{fig:Speed}) by performing a linear regression of $\ln(c_*)$ as a function of $\ln(f'(0))$. We find that $\ln(c_*) \sim a_1 \ln(f'(0))+a_0$ with $(a_1,a_0)\simeq(0.5306,-0.5012)$, where the relative error of the coefficient $a_1$ compared to the predicted value of $1/2$ is approximately $0.0613$.

\begin{figure}[t!]
\centering
\includegraphics[width=.43\textwidth]{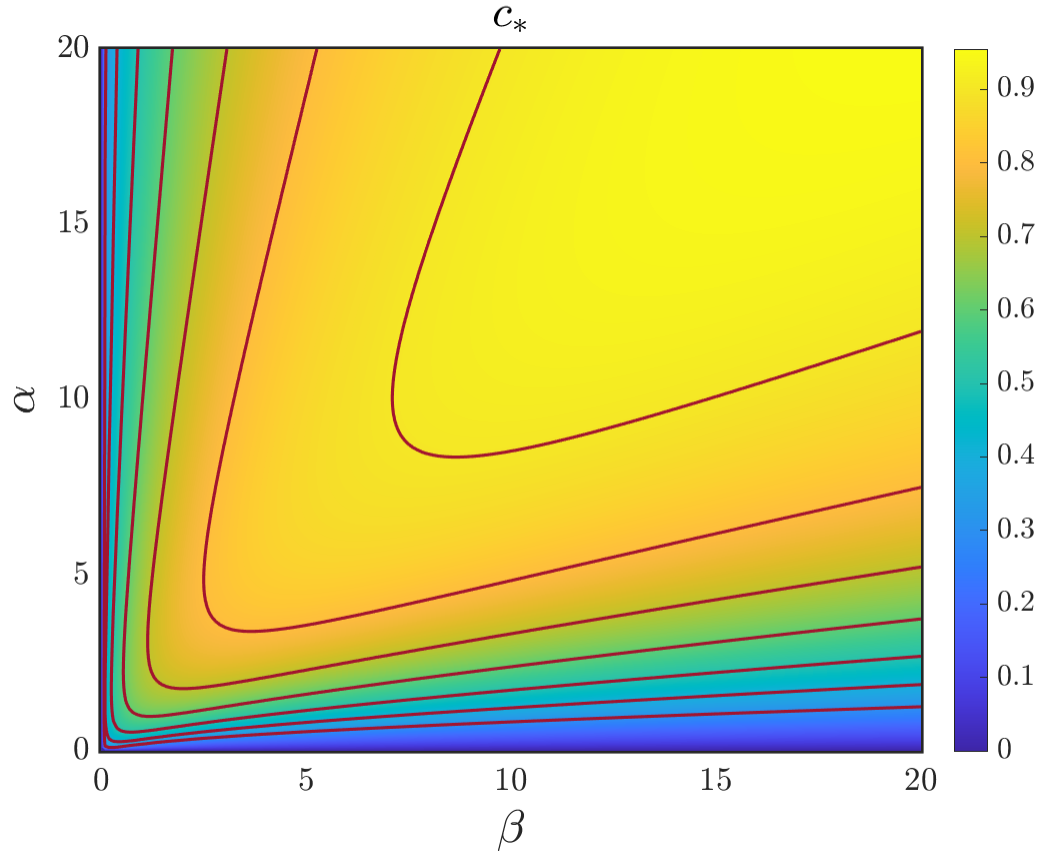}\hspace{0.6cm}
\includegraphics[width=.4\textwidth]{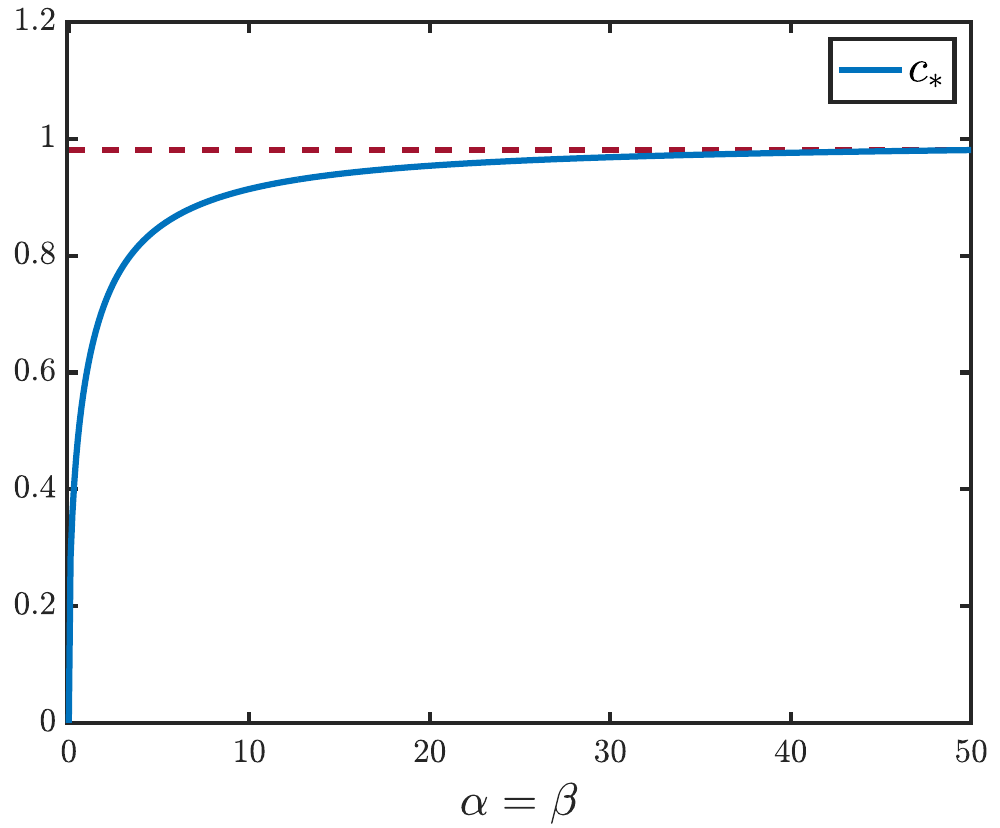}
\caption{Left: Amplitude of the spreading speed $c_*$ as a function of $(\alpha,\beta)$ in the square $[0,20]\times[0,20]$. Several isolines (red curves) are also reported. Right: Linear spreading speed $c_*$ as a function of $\alpha=\beta$. Other values of the parameters are set to $(d,f'(0))=(1,1)$.}
 \label{fig:SpeedAB}
\end{figure}

We have also further explored the dependence of the spreading speed as a function of $\alpha$ and $\beta$ by showing in Figure~\ref{fig:SpeedAB} the color plot of the map $(\alpha,\beta)\mapsto c_*(\alpha,\beta)$ and several of its isolines (red curves). It shows that the spreading speed seems to converge towards a limiting value as $\alpha=\beta\rightarrow+\infty$. We numerically confirmed this behavior by plotting the linear spreading speed $c_*$ as a function of $\alpha=\beta$ where we observe a monotone convergence towards a limiting asymptotic value, see the right panel of Figure~\ref{fig:SpeedAB}.  We leave it as  future work to theoretically investigate this asymptotic limit.

\section{Asymptotic spreading}\label{secAS}

In this section, we investigate the asymptotic spreading properties of system \eqref{eq1}-\eqref{bd} starting from compactly supported initial conditions. We anticipate that the linear spreading speed $c_*$ defined in the previous section via formula \eqref{c_*} is precisely the asymptotic spreading speed of the nonlinear system \eqref{eq1}-\eqref{bd} as stated in the following theorem.

\begin{theorem}\label{thmSpSp}
Let $(\mathbf{v},\boldsymbol{\rho})$ be the unique bounded classical solution of the Cauchy problem \eqref{eq1}-\eqref{bd}-\eqref{inv} starting from a nontrivial bounded compactly supported initial datum $(0,0)\not\equiv(\mathbf{h},\boldsymbol{\Lambda})\leq \left(\frac{\beta}{\alpha},1\right)$. Let $c_*>0$ be defined in \eqref{c_*}. Then:
\begin{itemize}
\item[(i)] for all $c>c_*$, we have
\bqs
\underset{t\to+\infty}{\lim}~\underset{\substack{|j| \geq ct\\ x\in[0,1]}}{\sup}(v_{j}(t,x),\rho_{j}(t))=\left(0,0\right),
\eqs
\item[(ii)] for all $c\in(0,c_*)$, we have
\bqs
\underset{t\to+\infty}{\lim}~\underset{\substack{|j| \leq ct  \\ x\in[0,1]}}{\inf}(v_{j}(t,x),\rho_{j}(t))=\left(\frac{\beta}{\alpha},1\right).
\eqs
\end{itemize}
\end{theorem}

\begin{figure}[t!]
\centering
\includegraphics[width=.4\textwidth]{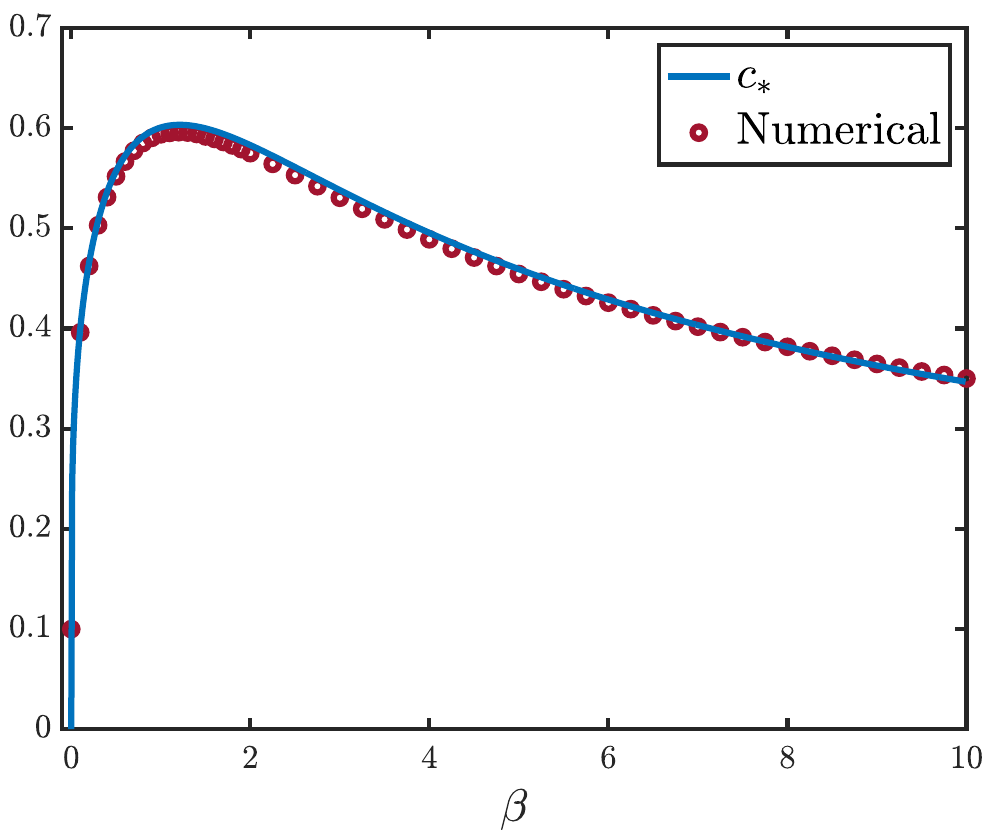}
\caption{Comparison between the theoretical spreading speed $c_*$ given by formula \eqref{c_*} and numerically computed spreading speed (dark red circles) obtained by numerically solving system \eqref{eq1}-\eqref{bd} from compactly supported initial conditions.}
 \label{fig:SpeedNumTh}
\end{figure}

We illustrate the above result in Figure~\ref{fig:SpeedNumTh} by directly comparing the theoretical spreading speed $c_*$ given by formula \eqref{c_*} and numerically computed spreading speed (dark red circles) obtained by numerically solving system \eqref{eq1}-\eqref{bd} from compactly supported initial conditions using the numerical scheme proposed in \cite{BF21a}.

\subsection{Upper estimate}

We first prove item (i) of Theorem~\ref{thmSpSp}, which is a direct consequence of the analysis conducted in the previous section. Let $c_*>0$ be given by formula \eqref{c_*} and let $\lambda_*>\lambda_0$ be such that 
\bqs
c_* = \underset{\lambda>\lambda_0}{\min}~\frac{\lambda}{\mu(\lambda)}=  \frac{\lambda_*}{\mu_*},
\eqs
where we have set $\mu_*:=\mu(\lambda_*)>0$. Then the following sequence
\begin{equation}\label{exp}
\forall t\geq0,\quad j\in\Z, \quad x\in[0,1], \quad (v_{j}(t,x),\rho_{j}(t))=\left(\ee^{-\mu_*(j-c_*t)}V_*(x),\ee^{-\mu_*(j-c_*t)}\right),
\end{equation}
where 
\begin{align*}
\forall x\in[0,1], \quad V_*(x)&=\frac{\beta}{\sinh\left(\sqrt{\frac{\lambda_*}{d}}\right)\Delta(\lambda_*)}\left[\sqrt{\lambda_* d}\cosh\left(\sqrt{\frac{\lambda_*}{d}} (1-x)\right)+\alpha\sinh\left(\sqrt{\frac{\lambda_*}{d}}(1- x) \right)\right]\\
	&~~~+\frac{\beta\ee^{-\mu_*}}{\sinh\left(\sqrt{\frac{\lambda_*}{d}}\right)\Delta(\lambda_*)}\left[\sqrt{\lambda_* d}\cosh\left(\sqrt{\frac{\lambda_*}{d}} x\right)+\alpha\sinh\left(\sqrt{\frac{\lambda_*}{d}}x \right)\right], 
	\end{align*}
is a solution of the linearized problem \eqref{leq1}. We readily remark that $V_*(x)>0$ for all $x\in[0,1]$. We can then introduce the sequences
\bqs
\forall t\geq0,\quad j\in\Z, \quad x\in[0,1], \quad \overline{v}_j(t,x)=\min\left( \vartheta \ee^{-\mu_*(j-c_*t)}V_*(x) ,\frac{\beta}{\alpha}\right),
\eqs
and
\bqs
\forall t\geq0,\quad j\in\Z, \quad  \overline{\rho}_j(t)=\min\left(\vartheta \ee^{-\mu_*(j-c_*t)},1\right),
\eqs
for some $\vartheta>0$ to be fixed. Since the initial datum $(0,0)\not\equiv(\mathbf{h},\boldsymbol{\Lambda})\leq \left(\frac{\beta}{\alpha},1\right)$ is assumed to be compactly supported, we can always find $\vartheta>0$ sufficiently large such that
\bqs
\forall j\in\Z, \quad x\in[0,1], \quad h_j(x)\leq \overline{v}_j(0,x) \quad \text{ and } \quad \Lambda_j\leq \overline{\rho}_j(0).
\eqs
From the comparison principle of Proposition~\ref{propCPcauchy}, we deduce that
\bqs
\forall t\geq0,\quad j\in\Z, \quad x\in[0,1], \quad   v_j(t,x)\leq \overline{v}_j(t,x) \quad \text{ and } \quad \rho_j(t) \leq \overline{\rho}_j(t),
\eqs
and thus for all $c>c_*$ one has
\bqs
\underset{t\to+\infty}{\lim}~\underset{\substack{j \geq ct\\ x\in[0,1]}}{\sup}( v_{j}(t,x),\rho_{j}(t))\leq \underset{t\to+\infty}{\lim}~\underset{\substack{j \geq ct\\ x\in[0,1]}}{\sup}( \overline{v}_{j}(t,x),\overline{\rho}_{j}(t)) =  \left(0,0\right).
\eqs
By symmetry, that is using $\overline{v}_{-j}(t,x)$ and $\overline{\rho}_{-j}(t)$ instead, we deduce item (i) of Theorem~\ref{thmSpSp}.

\subsection{Lower estimate}

Our aim is now to prove the lower estimate of item (ii) of Theorem~\ref{thmSpSp}. For that purpose, we shall construct compactly supported subsolutions of the linear system penalized by $\delta>0$ which reads
\bqq\label{lindelta}
\forall t>0,~j\in\Z,\quad	\begin{cases}
		\partial_{t} v_{j}(t, x)=d\partial_x^2v_{j}(t,x),\quad x\in(0,1),\\
		\rho_{j}^{\prime}(t)=\left(f^{\prime}(0)-\delta\right)\rho_{j}(t)+\alpha(v_{j}(t,0)+v_{j-1}(t,1))-2\beta\rho_{j}(t),
	\end{cases}
\eqq
together with the usual boundary conditions
\bqq\label{bddelta}
\forall t>0,~j\in\Z,\quad	\begin{cases}
		-d\partial_{x}v_{j}(t,0)+\alpha v_{j}(t,0)=\beta \rho_{j}(t),\\
		d\partial_{x}v_{j}(t,1)+\alpha v_{j}(t,1)=\beta \rho_{j+1}(t).
	\end{cases}
\eqq

The main result of this section is the following.

\begin{proposition}\label{propsubsoldelta}
Let $c_*$ be given by formula \eqref{c_*}. For all $c\in(0,c_*)$ close enough to $c_*$ there exists $\delta>0$ such that the penalized linear system \eqref{lindelta}-\eqref{bddelta} admits a nonnegative, compactly supported, generalized subsolution $(\underline{\mathbf{v}},\underline{\boldsymbol{\rho}})\not\equiv(0,0)$.
\end{proposition}

\begin{Proof}
In order to keep the presentation as light as possible, we will proceed with $f'(0)$ instead of $f^{\prime}(0)-\delta$ in \eqref{lindelta} since our arguments naturally perturb for $\delta>0$ small enough. 

We set $c\in(0,c_*)$ and consider once again exponential solutions of \eqref{lindelta}-\eqref{bddelta} of the form
\bqs
	(v_{j}(t,x),\rho_{j}(t))=\left(\ee^{-\mu(j-ct)}V(x),\ee^{-\mu(j-ct)}\right),
\eqs
where 
\bqs
	V(x)=a\cosh\left(\sqrt{\frac{\lambda}{d}} x\right)+b\sinh\left(\sqrt{\frac{\lambda}{d}} x\right), \quad\forall x\in[0,1],
\eqs
this time with eventual complex parameters $(\lambda,\mu,a,b)\in\C^4$ that will be fixed along the proof. Performing similar computations as in the previous section, we readily obtain that, given $c\in(0,c_*)$, the couple $(\lambda,\mu)\in\C^2$ is a solution of the system \eqref{muc} from which one obtains equation \eqref{c} which we rewrite as
\bqs
\Phi(c,\lambda)=0,
\eqs
with $\Phi(c,\lambda):=\lambda -c \mu(\lambda)$ where $\mu(\lambda)=\ln\left(y(\lambda)+\sqrt{y(\lambda)^2-1}\right)$ and $y(\lambda)$ is given in \eqref{y}. By definition of $c_*$ and analyticity of the map $\Phi$ on its domain of definition we have that there exists a positive integer $p\geq1$ such that
\bqs
\Phi(c_*,\lambda_*)=0,\quad \partial_\lambda^{k} \Phi(c_*,\lambda_*)=0 \quad\text{ for }\quad k=1,\dots,2p-1\quad\text{ and }\quad \partial_\lambda^{2p} \Phi(c_*,\lambda_*)>0.
\eqs
Next, introducing the auxiliary variables   
\bqs
\xi := c_*-c>0 \quad\text{ and }\quad z := \lambda-\lambda_*\in\C,
\eqs
we see that $\Phi(c,\lambda)=0$ is equivalent, in a neighborhood of $(c,\lambda)=(c_*,\lambda_*)$, to
\bqs
\mu_*\xi+ \mu'(\lambda_*)\xi z+a_* z^{2p} = \phi(z,\xi), \quad a_*:=\frac{\partial_\lambda^{2p} \Phi(c_*,\lambda_*)}{(2p)!}=-\frac{c_*}{(2p)!}\mu^{(2p)}(\lambda_*)>0,
\eqs
where $\phi$ is analytic in a neighborhood of $(0,0)$ and $\phi(z,\xi)=O\left(|z|^{2p+1}+\xi|z|^2\right)$ as $(z,\xi)\rightarrow(0,0)$. For small $\xi>0$, the polynomial equation $\mu_*\xi+ \mu'(\lambda_*)\xi z+a_* z^{2p} =0$ has $2p$ complex conjugate roots which writes
\bqs
z_\pm^k(\xi)= \left( \frac{\mu_*}{a_*}\xi\right)^{\frac{1}{2p}}\ee^{\mathbf{i}\left[\pm \frac{\pi}{2p}+\frac{2k\pi}{p} \right]}+O\left(\xi^{\frac{1}{p}}\right), \text{ for } k=0,\dots,p-1.
\eqs
Applying Rouch\'e's theorem, we get that the algebraic equation $\mu_*\xi+ \mu'(\lambda_*)\xi z+a_* z^{2p} = \phi(z,\xi)$ has also $2p$ complex roots which we denote by $\widetilde{z}_\pm^k(\xi)$ and these roots still satisfy
\bqs
\widetilde{z}_\pm^k(\xi)= \left( \frac{\mu_*}{a_*}\xi\right)^{\frac{1}{2p}}\ee^{\mathbf{i}\left[\pm \frac{\pi}{2p}+\frac{2k\pi}{p} \right]}+O\left(\xi^{\frac{1}{p}}\right), \text{ for } k=0,\dots,p-1.
\eqs
As a consequence, reverting to the full notation, we observe that for $c$ strictly less than and sufficiently close to $c_*$, the equation $\Phi(c,\lambda)=0$ admits a solution of the form
\bqs
\lambda = \lambda_*+\widetilde{z}_+^0(\xi),
\eqs
with the following properties:
\bqs
\Re(\lambda)=\lambda_*+O\left(\xi^{\frac{1}{2p}}\right)>0, \quad \Im(\lambda) = \left( \frac{\mu_*}{a_*}\xi\right)^{\frac{1}{2p}} \sin\left(\frac{\pi}{2p}\right)+O\left(\xi^{\frac{1}{p}}\right)>0.
\eqs
The corresponding profile $V$, given by
\begin{align*}
\forall x\in[0,1], \quad V(x)&=\frac{\beta}{\sinh\left(\sqrt{\frac{\lambda}{d}}\right)\Delta(\lambda)}\left[\sqrt{\lambda d}\cosh\left(\sqrt{\frac{\lambda}{d}} (1-x)\right)+\alpha\sinh\left(\sqrt{\frac{\lambda}{d}}(1- x) \right)\right]\\
	&~~~+\frac{\beta\ee^{-\mu(\lambda)}}{\sinh\left(\sqrt{\frac{\lambda}{d}}\right)\Delta(\lambda)}\left[\sqrt{\lambda d}\cosh\left(\sqrt{\frac{\lambda}{d}} x\right)+\alpha\sinh\left(\sqrt{\frac{\lambda}{d}}x \right)\right], 
	\end{align*}
satisfies 
\bqs
\forall x\in[0,1], \quad \Re\left(V(x)\right)=\Re\left(V_*(x)\right)+O\left(\xi^{\frac{1}{2p}}\right)>0, \quad \Im\left(V(x)\right)=O\left(\xi^{\frac{1}{2p}}\right)\neq 0
\eqs
and
\bqs
\forall x\in[0,1], \quad \mathrm{Arg}(V(x))=O\left(\xi^{\frac{1}{2p}}\right),
\eqs
where we denoted by $\mathrm{Arg}(V(x))\in(-\pi,\pi]$ the principal argument of $V(x)$. Taking the real parts of the just constructed exponential solutions, we set
\bqs
\widetilde{v}_j(t,x)=|V(x)|\ee^{-\frac{\lambda_*+\Re(\widetilde{z}_+^0(\xi))}{c}(j-ct)}\cos\left(\frac{\Im(\widetilde{z}_+^0(\xi))}{c}(j-ct)-\mathrm{Arg}(V(x))\right),
\eqs
and
\bqs
\widetilde{\rho}_j(t)=\ee^{-\frac{\lambda_*+\Re(\widetilde{z}_+^0(\xi))}{c}(j-ct)}\cos\left(\frac{\Im(\widetilde{z}_+^0(\xi))}{c}(j-ct)\right),
\eqs
for all $t\geq0$, $j\in\Z$ and $x\in[0,1]$. Our aim is now to truncate $\widetilde{v}_j(t,x)$ and $\widetilde{\rho}_j(t)$ to obtain non negative generalized subsolutions.

Since, $\mathrm{Arg}(V(x))=O\left(\xi^{\frac{1}{2p}}\right)$ for all $x\in[0,1]$, we introduce the following rescaled function: 
\bqs
\Phi(x):=\frac{\mathrm{Arg}(V(x))}{\Im(\widetilde{z}_+^0(\xi))},
\eqs
and define
\bqs
\underline{\varphi}:=\min_{x\in[0,1]}\Phi(x) \quad \text{ and } \quad \overline{\varphi}:=\max_{x\in[0,1]}\Phi(x).
\eqs
Finally, we set $\Phi_{\mathrm{per}}$ the $1$-periodic extension of $\Phi$. For all $t\geq0$, we can then define the set
\bqs
F(t):=\left\{y\in\R~|~ ct-\frac{c\pi}{2\Im(\widetilde{z}_+^0(\xi))}<\lfloor y\rfloor-c\Phi_{\mathrm{per}}(y)<ct+\frac{c\pi}{2\Im(\widetilde{z}_+^0(\xi))}\right\},
\eqs
which may have several connected components. Upon denoting
\bqs
\underline{\xi}(t):=ct-\frac{c\pi}{2 \Im(\widetilde{z}_+^0(\xi))}+c\overline{\varphi} \quad \text{ and } \quad \overline{\xi}(t):= ct+\frac{c\pi}{2 \Im(\widetilde{z}_+^0(\xi))}+c\underline{\varphi},
\eqs
we only retain the connected component of $F(t)$, denoted by $G(t)$, which contains the interval
\bqs
\left\{ y\in\R~|~ \underline{\xi}(t)< \lfloor y\rfloor <\overline{\xi}(t)\right\}\subset G(t).
\eqs
As a consequence, for each $t\geq0$, there exists $\underline{y}(t)\leq \underline{\xi}(t)< \overline{\xi}(t)\leq \overline{y}(t)$, such that $G(t)=(\underline{y}(t),\overline{y}(t))$. Then, we define
\bqs
\underline{j}(t):=\left\lfloor\underline{y}(t)\right\rfloor \in\Z \quad \text{ and }\quad \overline{j}(t):=\left\lfloor\overline{y}(t)\right\rfloor \in\Z,
\eqs
together with
\bqs
\underline{y}(t)-\underline{j}(t)=:\underline{\epsilon}(t)\in[0,1),  \quad \text{ and }  \overline{y}(t)-\overline{j}(t)=: \overline{\epsilon}(t)\in[0,1).
\eqs
With all the previous notations, we can introduce the following generalized subsolution for the $v$ component
\bqs
\forall t\geq0, \quad \underline{v}_j(t,x):=\left\{\begin{split}
\widetilde{v}_j(t,x), & \quad \underline{j}(t)+1 \leq j \leq \overline{j}(t)-1, \quad x\in [0,1],\\
\widetilde{v}_j(t,x), & \quad j=\underline{j}(t) \text{ and } x\in[\underline{\epsilon}(t),1],\\
0, & \quad j=\underline{j}(t) \text{ and } x\in[0,\underline{\epsilon}(t)],\\
\widetilde{v}_j(t,x), & \quad j=\overline{j}(t) \text{ and } x\in[0,\overline{\epsilon}(t)],\\
0, & \quad j=\overline{j}(t) \text{ and } x\in[\overline{\epsilon}(t),1],\\
0, & \quad \text{ otherwise.}
\end{split}
\right.
\eqs
Regarding the $\rho$ component, we simply define
\bqs
\forall t\geq0, \quad \underline{\chi}(t)=\left\lceil ct-\frac{c\pi}{2 \Im(\widetilde{z}_+^0(\xi))} \right\rceil, \quad \text{ and } \quad \overline{\chi}(t)=\left\lfloor ct+\frac{c\pi}{2 \Im(\widetilde{z}_+^0(\xi))}\right\rfloor,
\eqs
and set
\bqs
\underline{\rho}_j(t):=\left\{\begin{split}
\widetilde{\rho}_j(t),  & \quad \underline{\chi}(t) \leq j \leq \overline{\chi}(t),\\
0, & \quad \text{ otherwise.}\end{split}\right.
\eqs
We claim that $\Phi(x)<0$ for all $x\in[0,1]$ such that one has the ordering
\bqs
\forall t\geq0, \quad \underline{j}(t)<\underline{\chi}(t)<\overline{j}(t)<\overline{\chi}(t).
\eqs
Let us now check that $(\underline{v}_j(t,x),\underline{\rho}_j(t))$ provides a generalized subsolution to the linear system \eqref{lindelta}-\eqref{bddelta}.
\begin{itemize}
\item Case $j\leq \underline{j}(t)-1$. In that case, one has $\underline{v}_j(t,x)=\underline{\rho}_j(t)=0$ for all $x\in[0,1]$ with $\underline{v}_j(t,0)=\underline{v}_{j-1}(t,0)=0$. We also remark that $\underline{\rho}_{\underline{j}(t)}(t)=0$ since $\underline{j}(t)<\underline{\chi}(t)$. As a consequence $(\underline{v}_j(t,x),\underline{\rho}_j(t))$ is a subsolution for all $j\leq \underline{j}(t)-1$.
\item Case $j= \underline{j}(t)$. First of all, we have that $\partial_x \underline{v}_{j}(t,\underline{\epsilon}(t)^+)-\partial_x \underline{v}_{j}(t,\underline{\epsilon}(t)^-)=\partial_x \widetilde{v}_{j}(t,\underline{\epsilon}(t)^+)\geq0$, such that $\underline{v}_{j}$ is a generalized subsolution on $(0,1)$. Regarding the boundary conditions, we have
\bqs
\begin{cases}
-d\partial_{x}\underline{v}_{j}(t,0)+\alpha \underline{v}_{j}(t,0)-\beta \underline{\rho}_{j}(t)=0,\\
d\partial_{x}\underline{v}_{j}(t,1)+\alpha \underline{v}_{j}(t,1)- \beta \underline{\rho}_{j+1}(t)=\beta\left(\widetilde{\rho}_{j+1}(t)-\underline{\rho}_{j+1}(t)\right)\leq0.
\end{cases}
\eqs
For the last inequality, it comes from the fact that since $\xi$ can always be chosen small enough, and thus $\Im(\widetilde{z}_+^0(\xi))$, one can always ensure that  $\underline{\rho}_{j}(t)\geq \widetilde{\rho}_{j}(t)$ for all $j\in \llbracket \underline{j}(t), \chi(t)\rrbracket$. Finally, at $j=\underline{j}(t)$, we have
\bqs
\underbrace{\underline{\rho}_j'(t)}_{= 0}+(2\beta-f'(0))\underbrace{\underline{\rho}_j(t)}_{=0}-\alpha\left( \underbrace{\underline{v}_j(t,0)+\underline{v}_{j-1}(t,1)}_{=0}\right)=0.
\eqs
\item Case $\underline{j}(t)+1\leq j \leq \underline{\chi}(t)-1$. In that range, we have $\underline{v}_{j}(t,x)=\widetilde{v}_{j}(t,x)>0$ for all $x\in[0,1]$ and $\underline{\rho}_{j}(t)=0$. As a consequence, we only need to check the boundary conditions
\bqs
\begin{cases}
-d\partial_{x}\underline{v}_{j}(t,0)+\alpha \underline{v}_{j}(t,0)-\beta \underline{\rho}_{j}(t)=\beta\left(\widetilde{\rho}_{j}(t)-\underline{\rho}_{j}(t)\right) \leq0,\\
d\partial_{x}\underline{v}_{j}(t,1)+\alpha \underline{v}_{j}(t,1)- \beta \underline{\rho}_{j+1}(t)=\beta\left(\widetilde{\rho}_{j+1}(t)-\underline{\rho}_{j+1}(t)\right)\leq0,
\end{cases}
\eqs
since $\underline{\rho}_{j}(t)\geq \widetilde{\rho}_{j}(t)$ for all $j\in \llbracket \underline{j}(t), \chi(t)\rrbracket$. Finally, we have
\bqs
\underbrace{\underline{\rho}_j'(t)}_{= 0}+(2\beta-f'(0))\underbrace{\underline{\rho}_j(t)}_{=0}-\alpha\left( \underbrace{\underline{v}_j(t,0)+\underline{v}_{j-1}(t,1)}_{>0}\right)<0.
\eqs
\item Case $\underline{\chi}(t)\leq j \leq \overline{j}(t)-1$. In that range, we still have $\underline{v}_{j}(t,x)=\widetilde{v}_{j}(t,x)>0$ for all $x\in[0,1]$, but this time one also has that $\underline{\rho}_{j}(t)=\widetilde{\rho}_j(t)$, and  $(\underline{v}_j(t,x),\underline{\rho}_j(t))$ exactly solves the linear system \eqref{lindelta}-\eqref{bddelta}.
\item Case $j=\overline{j}(t)$. This time, we check that $\partial_x \underline{v}_{j}(t,\overline{\epsilon}(t)^+)-\partial_x \underline{v}_{j}(t,\overline{\epsilon}(t)^-)=-\partial_x \widetilde{v}_{j}(t,\overline{\epsilon}(t)^-)\geq0$, such that $\underline{v}_{j}$ is a generalized subsolution on $(0,1)$. Regarding the boundary conditions, we have
\bqs
\begin{cases}
-d\partial_{x}\underline{v}_{j}(t,0)+\alpha \underline{v}_{j}(t,0)-\beta \underline{\rho}_{j}(t)=0,\\
d\partial_{x}\underline{v}_{j}(t,1)+\alpha \underline{v}_{j}(t,1)- \beta \underline{\rho}_{j+1}(t)=-\beta\widetilde{\rho}_{j+1}(t)\leq0.
\end{cases}
\eqs
Finally, we have
\begin{align*}
&\underline{\rho}_j'(t)+(2\beta-f'(0))\underline{\rho}_j(t)-\alpha\left( \underline{v}_j(t,0)+\underline{v}_{j-1}(t,1)\right)\\&=\widetilde{\rho}_j'(t)+(2\beta-f'(0))\widetilde{\rho}_j(t)-\alpha\left( \widetilde{v}_j(t,0)+\widetilde{v}_{j-1}(t,1)\right)=0.
\end{align*}
\item Case $\overline{j}(t)+1\leq j \leq \overline{\chi}(t)$. Let us first note that for all $j\geq \overline{j}(t)+1$ one has $\underline{v}_j(t,x)=0$, such that the boundary conditions read
\bqs
\begin{cases}
-d\partial_{x}\underline{v}_{j}(t,0)+\alpha \underline{v}_{j}(t,0)-\beta \underline{\rho}_{j}(t)=-\beta \underline{\rho}_{j}(t)\leq0,\\
d\partial_{x}\underline{v}_{j}(t,1)+\alpha \underline{v}_{j}(t,1)- \beta \underline{\rho}_{j+1}(t)=- \beta \underline{\rho}_{j+1}(t)\leq0.
\end{cases}
\eqs
Let us remark that in the above computations, when $j=\overline{\chi}(t)$, we have used the fact that $\underline{\rho}_{j+1}(t)=0$, while for all $j\leq \overline{\chi}(t)$ we used the fact that $\underline{\rho}_{j}(t)=\widetilde{\rho}(j)(t)>0$. Finally, we also have
\begin{align*}
&\underline{\rho}_j'(t)+(2\beta-f'(0))\underline{\rho}_j(t)-\alpha\left( \underbrace{\underline{v}_j(t,0)+\underline{v}_{j-1}(t,1)}_{=0}\right) \\
&=\widetilde{\rho}_j'(t)+(2\beta-f'(0))\widetilde{\rho}_j(t)\\
&=\alpha\left( \widetilde{v}_j(t,0)+\widetilde{v}_{j-1}(t,1)\right)\leq0.
\end{align*}
The last inequality comes, once again, by the fact that $\xi$ can always be chosen small enough so as to ensure that $\widetilde{v}_j(t,x)\leq0$ for all $\overline{j}(t)+1\leq j \leq \overline{\chi}(t)$ and $x\in[0,1]$.
\item Case $j>\overline{\chi}(t)$. In that region, we have that $\underline{v}_j(t,x)=\underline{\rho}_j(t)=0$ for all $x\in[0,1]$ with $\underline{v}_j(t,0)=\underline{v}_{j-1}(t,0)=0$.
\end{itemize}
This concludes the proof of the proposition.
\end{Proof}

\begin{Proof}[ of the claim.] Let us show that $\Phi(x)<0$ for all $x\in[0,1]$. By definition, we may write $V(x)=P(x)\ee^{\mathbf{i}\Theta(x)}$ where $P(x)=|V(x)|$ and $\Theta(x)=\mathrm{Arg}(V(x))$ satisfy 
\bqs
\forall x\in(0,1), \quad d\left( P''(x)-(\Theta'(x))^2P(x)+\mathbf{i} \Theta''(x)P(x)+2\mathbf{i}\Theta'(x)P'(x)\right)-(\lambda_*+\widetilde{z}_+^0(\xi))P(x)=0,
\eqs
together with
\bqs
\begin{split}
-d\left( P'(0)+\mathbf{i}\Theta'(0)P(0)\right)+\alpha P(0) &= \beta \ee^{-\mathbf{i}\Theta(0)},\\
d\left( P'(1)+\mathbf{i}\Theta'(1)P(1)\right)+\alpha P(1) &= \beta \ee^{-\mu-\mathbf{i}\Theta(1)}.
\end{split}
\eqs
Keeping the imaginary parts of the above equations, we obtain
\bqs
d\Theta''(x)P(x)+2d\Theta'(x)P'(x)-\Im(\widetilde{z}_+^0(\xi))P(x)=0, \quad x\in(0,1),
\eqs
and
\bqs
dP(0)\Theta'(0)=\beta \sin(\Theta(0)),\quad
dP(1)\Theta'(1)=-\beta \ee^{-\frac{\lambda_*}{c}}\sin\left(\frac{\Im(\widetilde{z}_+^0(\xi))}{c}+\Theta(1)\right),
\eqs
while the real parts give
\bqs
dP''(x)-\left[d(\Theta'(x))^2+\lambda_*+\Re(\widetilde{z}_+^0(\xi)) \right]P(x) = 0, \quad x\in(0,1),
\eqs
and
\bqs
-dP'(0)+\alpha P(0) = \beta \cos(\Theta(0)),\quad
dP'(1)+\alpha P(1) = \beta  \ee^{-\frac{\lambda_*}{c}} \cos\left(\frac{\Re(\widetilde{z}_+^0(\xi))}{c}+\Theta(1)\right).
\eqs
Using the rescaling $\Theta(x)=\Im(\widetilde{z}_+^0(\xi))\Phi(x)$, with $\Phi(x)=O(1)$ as $\xi\rightarrow0$, we obtain
\bqs
d\Phi''(x)P(x)+2d\Phi'(x)P'(x)-P(x)=0, \quad x\in(0,1),
\eqs
and
\bqs
dP(0)\Phi'(0)=\beta \frac{\sin(\Im(\widetilde{z}_+^0(\xi))\Phi(0))}{\Im(\widetilde{z}_+^0(\xi))},\quad
dP(1)\Phi'(1)=-\beta \ee^{-\frac{\lambda_*}{c}} \frac{\sin\left(\Im(\widetilde{z}_+^0(\xi))\left[\frac{1}{c}+\Phi(1)\right]\right)}{\Im(\widetilde{z}_+^0(\xi))}.
\eqs
Let us finally remark that $P(x)=V_*(x)+O\left(\xi^{\frac{1}{2p}}\right)>0$ for all $x\in[0,1]$ provided that $\xi$ is taken small enough, since $V_*(x)>0$ for all $x\in[0,1]$. Integrating once the equation for $\Phi$, we obtain
\bqs
\forall x\in[0,1], \quad \Phi'(x)=\Phi'(0)\left(\frac{P(0)}{P(x)}\right)^2+\frac{1}{d}\int_0^x \left(\frac{P(y)}{P(x)}\right)^2\d y.
\eqs
We claim that necessarily $\Phi'(0)<0$. By contradiction, assume that it is not true and that $\Phi'(0)\geq 0$, then one has $\Phi'(x)>0$ for all $x\in(0,1]$. On the other hand, using the boundary condition at $x=0$, we infer that $\Phi(0)\geq 0$ and thus $0 \leq \Phi(0)<\Phi(1) $ which implies that $\Phi'(1)<0$ by the boundary condition at $x=1$, which is a contradiction since $\Phi'(x)>0$ for all $x\in(0,1]$. As a consequence, either $\Phi'(x)<0$ for all $x\in[0,1)$ and $\Phi'(1)\leq0$, and in that case $\Phi(1)<\Phi(0)<0$, and the boundary condition at $x=1$ implies that $-\frac{1}{c}\leq \Phi(1)$ and thus $-\frac{1}{c}\leq \Phi(x)<0$ for all $x\in[0,1]$. Or there exists a unique $x_0\in(0,1)$ such that $\Phi'(x_0)=0$, since at such a point one has $\Phi''(x_0)=1/d>0$. In particular, we deduce that $\Phi'(1)>0$ and thus $\Phi(1)<-\frac{1}{c}$ which implies that $\Phi(x_0)\leq \Phi(x)\leq \max(\Phi(0),\Phi(1))<0$ for all $x\in[0,1]$. \end{Proof}

\begin{Proof}[ of item (ii) of Theorem~\ref{thmSpSp}.] Let $c\in(0,c_*)$ and choose $c'\in(c,c_*)$ very close to $c_*$ such that, from the previous Proposition~\ref{propsubsoldelta}, we get the existence of $\delta>0$ such that the penalized linear system \eqref{lindelta}-\eqref{bddelta} admits a nonnegative, compactly supported, generalized subsolution that we denote $(\underline{\mathbf{v}}^{c',\delta},\underline{\boldsymbol{\rho}}^{c',\delta})\not\equiv(0,0)$. By regularity of the nonlinearity $f$, there exists $\iota>0$ such that
\bqs
(f'(0)-\delta)u\leq f(u), \quad 0\leq u\leq \iota.
\eqs
Then, one can find $\eta>0$, small enough, such that $\eta \underline{\rho}^{c',\delta}_j(t)\leq \iota$ for all $t\geq0$ and $j\in\Z$. As a consequence $(\eta\underline{\mathbf{v}}^{c',\delta},\eta\underline{\boldsymbol{\rho}}^{c',\delta})\not\equiv(0,0)$ is a nonnegative compactly supported subsolution to the full nonlinear system \eqref{eq1}-\eqref{bd}. By positivity of the solution of the nonlinear system \eqref{eq1}-\eqref{bd} ensured by Theorem~\ref{thmcauchy} and upon eventually reducing the size of $\eta>0$, we can always ensure that at time $t=1$ the unique solution $(\mathbf{v},\boldsymbol{\rho})$ of the Cauchy problem  \eqref{eq1}-\eqref{bd}-\eqref{inv} starting from the nontrivial bounded compactly supported initial datum $(0,0)\not\equiv(\mathbf{h},\boldsymbol{\Lambda})\leq \left(\frac{\beta}{\alpha},1\right)$ satisfies
\bqs
(\eta\underline{\mathbf{v}}^{c',\delta}(0),\eta\underline{\boldsymbol{\rho}}^{c',\delta}(0)) \leq (\mathbf{v}(1),\boldsymbol{\rho}(1)).
\eqs
From the comparison principle of Proposition~\ref{propCPcauchy} we obtain that
\bqs
\forall t\geq1, \quad (\eta\underline{\mathbf{v}}^{c',\delta}(t-1),\eta\underline{\boldsymbol{\rho}}^{c',\delta}(t-1)) \leq (\mathbf{v}(t),\boldsymbol{\rho}(t)).
\eqs
As a consequence, there exists $\nu\in(0,1)$ small such that
\bqs
v_{\lfloor c't \rfloor}(t,x) \geq \eta\underline{v}^{c',\delta}_{\lfloor c't \rfloor}(t-1,x)\geq\frac{\beta}{\alpha}\nu \quad \text{ and } \quad v_{\lfloor c't \rfloor+1}(t,x) \geq \eta\underline{v}^{c',\delta}_{\lfloor c't \rfloor+1}(t-1,x)\geq\frac{\beta}{\alpha}\nu,
\eqs
with
\bqs
\rho_{\lfloor c't \rfloor}(t) \geq \eta\underline{\rho}^{c',\delta}_{\lfloor c't \rfloor}(t-1)\geq\nu \quad \text{ and } \quad \rho_{\lfloor c't \rfloor+1}(t) \geq \eta\underline{\rho}^{c',\delta}_{\lfloor c't \rfloor+1}(t-1)\geq\nu,
\eqs
for all $t\geq1$ and $x\in[0,1]$. Here, we denote by $\lfloor x\rfloor$ the integer part of $x\in\R$. By a symmetry argument, we also obtain
\bqs
v_{-\lfloor c't \rfloor}(t,x) \geq \frac{\beta}{\alpha}\nu \quad \text{ and } \quad v_{-\lfloor c't \rfloor-1}(t,x) \geq \frac{\beta}{\alpha}\nu,
\eqs
and
\bqs
\rho_{-\lfloor c't \rfloor}(t) \geq \nu \quad \text{ and } \quad \rho_{-\lfloor c't \rfloor-1}(t) \geq \nu,
\eqs
for all $t\geq1$ and $x\in[0,1]$. Upon eventually reducing the size of $\nu$ and by positivity of the solution $(\mathbf{v},\boldsymbol{\rho})$ we can always ensure that
\bqs
v_j(1,x)\geq \frac{\beta}{\alpha}\nu \quad\text{ and }\quad \rho_j(1)\geq \nu\quad\text{ for all }\quad x\in[0,1] \quad\text{ and }\quad  -c'-1\leq j \leq c'+1.
\eqs
Since $\left(\frac{\beta}{\alpha}\nu,\nu\right)_{j\in\Z}$ is a homogeneous subsolution of \eqref{eq1}-\eqref{bd}, we can apply a variant of the comparison principle, Proposition~\ref{propCPcauchy}, but with two boundaries as stated in Proposition~\ref{propCPapp2} of the Appendix. More precisely, we set $\zeta(t)=-c't$ and $\xi(t)=c't$, from the previous analysis, we have $v_j(t,x)\geq \frac{\beta}{\alpha}\nu$ and $\rho_j(t)\geq\nu$ for all $t\geq1$, $j\in[\zeta(t)-1,\zeta(t))\cup(\xi(t),\xi(t)+1]$ and $x\in[0,1]$. Furthermore, at time $t=1$, we also have  $v_j(1,x)\geq \frac{\beta}{\alpha}\nu$ and $\rho_j(1)\geq\nu$ for all  $j\in[\zeta(1)-1,\xi(1)+1]$. As a consequence, the comparison principle with two boundaries ensures that
\bqs
\forall t\geq1,~x\in[0,1],~|j|\leq c't, \quad v_j(t,x)\geq \frac{\beta}{\alpha}\nu \quad\text{ and }\quad \rho_j(t)\geq\nu,
\eqs
and thus
\bqq\label{e5.2}
\forall x\in[0,1], \quad \underset{t\to+\infty}{\liminf}~\underset{|j|\leq c't }{\inf}(v_{j}(t,x),\rho_{j}(t))\geq \nu \left(\frac{\beta}{\alpha},1\right).
\eqq
Let us deduce from there that we have 
\bqs
\forall x\in[0,1], \quad \underset{t\to+\infty}{\lim}~\underset{|j|\leq ct }{\inf}(v_{j}(t,x),\rho_{j}(t))=  \left(\frac{\beta}{\alpha},1\right).
\eqs
Since we already know that the solution satisfies $(\mathbf{v},\boldsymbol{\rho})\leq \left(\frac{\beta}{\alpha},1\right)$, let us assume, by contradiction, that there exists two sequences $(t_n)_{n\in\N}$ and $(j_n)_{n\in\N}$ with $0<t_n\rightarrow+\infty$ and $j_n\in\Z$ be such that $|j_n|\leq c t_n$ for all $n\in\N$ and
\bqs
\underset{n\rightarrow+\infty}{\limsup}~ (v_{j_n}(t_n,x),\rho_{j_n}(t_n))<\left(\frac{\beta}{\alpha},1\right).
\eqs
We set $\tilde{v}_j^n(t,x):=v_{j+j_n}(t+t_n,x)$ and $\tilde{\rho}_j^n(t):=\rho_{j+j_n}(t+t_n)$ for each $n\in\N$ and let 
$(\tilde{\mathbf{v}}^n,\tilde{\boldsymbol{\rho}}^n)_{n\in\N}$ be  the corresponding sequence.  Fix any $(N,M)\in\Z^2$ with $N<M$ and any time interval $[\tau,T]\subset\R$ with $\tau<T$. Let $n_0\in\N$ large enough such that $\tau+t_n>0$ for all $n\geq n_0$, then by standard parabolic estimates (similar to the one obtained in the proof of Theorem~\ref{thmcauchy}), the sequence $((\tilde{v}^n_j)_{j=N,\cdots,M},(\tilde{\rho_j}^n)_{j=N,\cdots,M})_{n\geq n_0}$ together with its respective time derivative and space derivatives is uniformly bounded in $\mathscr{C}^{0,\nu'}$ norm on $[\tau,T]\times[0,1]$ for some $\nu'\in(0,1)$. As a consequence, up to a subsequence, $(\tilde{\mathbf{v}}^n,\tilde{\boldsymbol{\rho}}^n)$ converges (locally uniformly) to $(\tilde{\mathbf{v}}^\infty,\tilde{\boldsymbol{\rho}}^\infty)$ a time entire solution of \eqref{eq1}-\eqref{bd}. From \eqref{e5.2}, since $c'>c$, we have that $(\tilde{v}^\infty_j(t,x),\tilde{\rho}^\infty_j(t))\geq \nu \left(\frac{\beta}{\alpha},1\right)$ for all $t\in\R$, $j\in\Z$ and $x\in[0,1]$. By a straightforward comparison to the spatially homogeneous system (solutions independent of $j$ and $x$), we infer that necessarily   $(\tilde{v}^\infty_j(t,x),\tilde{\rho}^\infty_j(t))=\left(\frac{\beta}{\alpha},1\right)$ for all $t\in\R$, $j\in\Z$ and $x\in[0,1]$, which contradicts the fact that $(\tilde{v}^\infty_0(0,x),\tilde{\rho}^\infty_0(0))<\left(\frac{\beta}{\alpha},1\right)$ by construction. As a consequence, we have thus proved that
\bqs
\forall x\in[0,1], \quad \underset{t\to+\infty}{\lim}~\underset{|j|\leq ct }{\inf}(v_{j}(t,x),\rho_{j}(t))=  \left(\frac{\beta}{\alpha},1\right),
\eqs  
and this concludes the proof of the theorem.
\end{Proof}

\section{Large diffusion limit}\label{secLDL}

Motivated by our numerical finding (see panel (c) of Figure~\ref{fig:Speed}), in this section, we study the asymptotic regime when $d\rightarrow+\infty$. For that purpose, we first set $\epsilon:=1/d>0$ such that system \eqref{eq1}-\eqref{bd} rewrites
\bqq\label{eqeps}
\forall t>0,~j\in\Z, \quad \left\{
	\begin{split}
		\epsilon\partial_{t} v_{j}(t, x)&=\partial_{x}^2v_{j}(t,x),\quad x\in(0,1),\\
		\rho_{j}^{\prime}(t)&=f(\rho_{j}(t))+\alpha(v_{j}(t,0)+v_{j-1}(t,1))-2\beta\rho_{j}(t),
	\end{split}\right.
\eqq
with associated Robin boundary conditions
\begin{equation}\label{bdeps}
\forall t>0,~j\in\Z, \quad	\begin{cases}
		-\partial_{x}v_{j}(t,0)+\alpha \epsilon v_{j}(t,0)=\beta \epsilon \rho_{j}(t),\\
		\partial_{x}v_{j}(t,1)+\alpha \epsilon v_{j}(t,1)=\beta \epsilon \rho_{j+1}(t).
	\end{cases}
\end{equation}

\subsection{Derivation of the asymptotic limiting system}

Fix $\epsilon_0>0$. For each $\epsilon\in(0,\epsilon_0]$, we consider $\mathbf{h}^\epsilon=(h_j^\epsilon)_{j\in\Z} \in \mathcal{X}^4$ and $\boldsymbol{\Lambda}^\epsilon=(\Lambda_j^\epsilon)_{j\in\Z} \in \ell^\infty(\Z)$, satisfying the compatibility condition \eqref{compatibility}, with
$$
\mathcal{X}^4:=\left\{\mathbf{u}=(u_j)_{j\in\Z} ~|~ \forall j\in\Z,~ u_{j}\in\mathscr{C}^4([0,1],\R)\right\}.
$$
We further suppose that there exists some positive constant $\kappa>0$ such that 
\bqq
\label{ICeps}
0< \left\|{\mathbf{h}^\epsilon}\right\|_{\infty}+\left\| \boldsymbol{\Lambda}^\epsilon\right\|_{\ell^\infty(\Z)}  \leq \kappa, \quad \left\|{\mathbf{h}^\epsilon}'\right\|_{\infty}+\left\|{\mathbf{h}^\epsilon}''\right\|_{\infty}+\left\|{\mathbf{h}^\epsilon}'''\right\|_{\infty} \leq \epsilon \kappa \quad\text{and} \quad \left\|{\mathbf{h}^\epsilon}''''\right\|_{\infty} \leq \epsilon^2 \kappa,
\eqq
for all $\epsilon\in(0,\epsilon_0]$, and we also assume that the sequences $\mathbf{V}^0$ and $\mathbf{P}^0$ defined as the following limits
\bqq
\label{IClim}
\forall j\in\Z, \quad V_j^0=\underset{\epsilon\rightarrow0}{\lim}~  \int_0^1 h_j^\epsilon(x)\d x, \text{ and } P_j^0=\underset{\epsilon\rightarrow0}{\lim}~\Lambda_j^\epsilon,
\eqq
satisfy $(0,0)\not\equiv(\mathbf{V}^0,\mathbf{P}^0)\in\ell^\infty(\Z)\times\ell^\infty(\Z)$.

For each $\epsilon\in(0,\epsilon_0]$, we shall denote by $(\mathbf{v}^\epsilon,\boldsymbol{\rho}^\epsilon)$ the solution of \eqref{eqeps}-\eqref{bdeps} given by Theorem~\ref{thmcauchy} with initial condition
\begin{equation}\label{inveps}
\forall j\in\Z,\quad \begin{cases}
		v_{j}^\epsilon(0,x)=h_{j}^\epsilon(x),\quad x\in(0,1),\\
		\rho_{j}^\epsilon(0)=\Lambda_{j}^\epsilon,
	\end{cases}
\end{equation}
with  $(\mathbf{h}^\epsilon,\boldsymbol{\Lambda}^\epsilon)$ satisfying the above conditions. For all $t>0$, one has
\bqs
\forall t>0, ~ \forall j\in\Z,\quad 0<v_{j}^\epsilon(t,x)\leq \max\left\{\dfrac{\beta}{\alpha},\kappa,\dfrac{\beta}{\alpha}\kappa\right\},~ x\in[0,1], \text{ and } 0< \rho_{j}^\epsilon(t)\leq \max\left\{1,\kappa,\dfrac{\alpha}{\beta}\kappa\right\},
\eqs
from which we deduce that
\bqs
\forall t>0,~ \forall j\in\Z, \quad \left|{\rho_{j}^{\epsilon}}'(t)\right|\leq C_\kappa,
\eqs
for some constant $C_\kappa>0$. Let us set
\bqs
\forall t>0,\quad \forall j\in\Z, \quad \forall x\in[0,1], \quad w_j^\epsilon(t,x):=\partial_tv_j^\epsilon(t,x).
\eqs
It follows from \eqref{eqeps} and \eqref{bdeps} that $\mathbf{w}^\epsilon=(w_{j}^\epsilon)_{j\in\Z}$ is a solution of
\bqs
\forall t>0,~j\in\Z, \quad \epsilon\partial_{t} w_{j}^\epsilon(t, x)=\partial_{x}^2w_{j}^\epsilon(t,x),\quad x\in(0,1),
\eqs
and
\bqs
\forall t>0,~j\in\Z, \quad	\begin{cases}
		-\partial_{x}w_{j}^\epsilon(t,0)+\alpha \epsilon w_{j}^\epsilon(t,0)=\beta \epsilon {\rho_{j}^\epsilon}'(t),\\
		\partial_{x}w_{j}^\epsilon(t,1)+\alpha \epsilon w_{j}^\epsilon(t,1)=\beta \epsilon {\rho_{j+1}^\epsilon}'(t),
	\end{cases}
\eqs
with initial condition given by
\bqs
w_{j}^\epsilon(0, x)=\frac{1}{\epsilon} {h_j^\epsilon}''(x), \quad x\in[0,1].
\eqs
Thanks to our condition on $\mathbf{h}^\epsilon$, we have that
\bqs
\forall \epsilon\in(0,\epsilon_0],\quad \left\|\mathbf{w}^\epsilon(0)\right\|_{\infty} \leq \kappa.
\eqs
Defining the super solution
\bqs
\forall t\geq0, \quad j\in\Z, \quad x\in[0,1], \quad \overline{w}_j(t,x)=\max\left(\dfrac{\beta}{\alpha}C_\kappa,\kappa\right),
\eqs
the parabolic comparison principle implies that for each $\epsilon\in(0,\epsilon_0]$ one has
\bqs
\forall t>0, \quad \forall j\in\Z, \quad |\partial_tv_j^\epsilon(t,x)| \leq \max\left(\dfrac{\beta}{\alpha}C_\kappa,\kappa\right),
\eqs
and
\bqs
\forall t>0, \quad \forall j\in\Z, \quad |\partial_x^2 v_j^\epsilon(t,x)| \leq \epsilon_0\max\left(\dfrac{\beta}{\alpha}C_\kappa,\kappa\right),
\eqs
for all $x\in[0,1]$. On the other hand, we also have
\bqs
\partial_x v_j^\epsilon(t,x)=\partial_x v_j^\epsilon(t,0)+\epsilon \int_0^x w_j(t,y)\d y=\alpha \epsilon v_j^\epsilon(t,0)-\beta\epsilon \rho_j^\epsilon(t)+\epsilon \int_0^x w_j(t,y)\d y
\eqs
such that
\bqs
\left|\partial_x v_j^\epsilon(t,x)\right| \leq C(\epsilon_0,\kappa),
\eqs
for all $t>0$, $j\in\Z$, $x\in[0,1]$ and $\epsilon\in(0,\epsilon_0]$. Furthermore, from the equation on $\rho_j$, we also infer that
\bqs
\forall t>0,\, j\in\Z, \quad \left|{\rho_j^\epsilon}''(t)\right| \leq \widehat{C}_k.
\eqs

Next, we set
\bqs
\forall t>0,\quad \forall j\in\Z, \quad \forall x\in[0,1], \quad w_j^\epsilon(t,x):=\partial_{tx}^2v_j^\epsilon(t,x),
\eqs
which satisfies 
\bqs
\forall t>0,~j\in\Z, \quad \epsilon\partial_{t} w_{j}^\epsilon(t, x)=\partial_{x}^2w_{j}^\epsilon(t,x),\quad x\in(0,1),
\eqs
and
\bqs
\forall t>0,~j\in\Z, \quad	\begin{cases}
		w_{j}^\epsilon(t,0)=\alpha \epsilon \partial_t v_{j}^\epsilon(t,0)-\beta \epsilon {\rho_{j}^\epsilon}'(t),\\
		w_{j}^\epsilon(t,1)=\beta \epsilon {\rho_{j+1}^\epsilon}'(t)-\alpha \epsilon \partial_tv_{j}^\epsilon(t,1),
	\end{cases}
\eqs
with initial condition given by
\bqs
w_{j}^\epsilon(0, x)=\frac{1}{\epsilon} {h_j^\epsilon}'''(x), \quad x\in[0,1].
\eqs
We note that we have already obtained the uniform bounds
\bqs
\forall t>0,~j\in\Z, \quad \left|\alpha \epsilon \partial_t v_{j}^\epsilon(t,0)-\beta \epsilon {\rho_{j}^\epsilon}'(t)\right| \leq \widetilde{C}(\epsilon_0,\kappa), \quad \left|\beta \epsilon {\rho_{j+1}^\epsilon}'(t)-\alpha \epsilon \partial_tv_{j}^\epsilon(t,1)\right| \leq \widetilde{C}(\epsilon_0,\kappa), 
\eqs
for all $\epsilon\in(0,\epsilon_0]$. Thanks to our condition on $\mathbf{h}^\epsilon$, and using the comparison principle, we obtain that
\bqs
\forall t>0,~j\in\Z,\quad x\in[0,1], \quad |\partial_{tx}^2v_j^\epsilon(t,x)|\leq \max\left(\dfrac{\beta}{\alpha}\widetilde{C}(\epsilon_0,\kappa),\kappa \right),
\eqs
for all $\epsilon\in(0,\epsilon_0]$. And since $v_j^\epsilon(t,x)$ satisfies the identity $\epsilon \partial_{tx}^2v_j^\epsilon(t,x)=\partial_x^3v_j^\epsilon(t,x)$, we deduce that
\bqs
\forall t>0,~j\in\Z,\quad x\in[0,1], \quad |\partial_x^3v_j^\epsilon(t,x)|\leq  \epsilon_0 \max\left(\dfrac{\beta}{\alpha}\widetilde{C}(\epsilon_0,\kappa),\kappa \right)
\eqs
for all $\epsilon\in(0,\epsilon_0]$.

Finally, we set 
\bqs
\forall t>0,\quad \forall j\in\Z, \quad \forall x\in[0,1], \quad w_j^\epsilon(t,x):=\partial_{t}^2v_j^\epsilon(t,x),
\eqs
which satisfies 
\bqs
\forall t>0,~j\in\Z, \quad \epsilon\partial_{t} w_{j}^\epsilon(t, x)=\partial_{x}^2w_{j}^\epsilon(t,x),\quad x\in(0,1),
\eqs
and
\bqs
\forall t>0,~j\in\Z, \quad	\begin{cases}
		-\partial_{x}w_{j}^\epsilon(t,0)+\alpha \epsilon w_{j}^\epsilon(t,0)=\beta \epsilon {\rho_{j}^\epsilon}''(t),\\
		\partial_{x}w_{j}^\epsilon(t,1)+\alpha \epsilon w_{j}^\epsilon(t,1)=\beta \epsilon {\rho_{j+1}^\epsilon}''(t)
	\end{cases}
\eqs
with initial condition given by
\bqs
w_{j}^\epsilon(0, x)=\frac{1}{\epsilon^2} {h_j^\epsilon}''''(x), \quad x\in[0,1].
\eqs
Once again, thanks to our condition on $\mathbf{h}^\epsilon$ and our uniform estimate on $\left|{\rho_{j}^\epsilon}''(t)\right|$, the comparison principle yields
\bqs
\forall t>0,~j\in\Z, \,x\in[0,1], \quad \left|\partial_{t}^2v_j^\epsilon(t,x)\right| \leq \widehat{C}_\kappa,
\eqs
for all $\epsilon\in(0,\epsilon_0]$, which also provides a uniform bound to $\partial_{txx}^3v_j^\epsilon$ using the equation.

As a consequence, we have obtained that $v_j^\epsilon$, $\partial_xv_j^\epsilon$, $\partial_tv_j^\epsilon$, $\partial_x^2v_j^\epsilon$, $\rho_j^\epsilon$ and ${\rho_j^\epsilon}'$ are all uniformly bounded in $\mathscr{C}^{0,\nu}$ for any $\nu\in(0,1]$. We can then apply Arzela-Ascoli's theorem, together with a diagonal extraction argument, to obtain, up to a subsequence, the existence of a limit  $(\mathbf{U},\mathbf{P})$ with
\bqs
\forall x\in[0,1], \quad \underset{\epsilon\rightarrow0}{\lim}~ v_j^\epsilon(t,x)=U_j(t,x), \quad \underset{\epsilon\rightarrow0}{\lim}~ \rho_j^\epsilon(t)=P_j(t), 
\eqs
locally uniformly in $(t,j)\in(0,+\infty)\times\Z$. The convergence also holds for the respective time and space derivatives. At the limit, one has
\bqs
\forall t>0,~j\in\Z, \quad \left\{
	\begin{split}
		0&=\partial_{x}^2U_{j}(t,x),\quad x\in(0,1),\\
		P_{j}^{\prime}(t)&=f(P_{j}(t))+\alpha(U_{j}(t,0)+U_{j-1}(t,1))-2\beta P_{j}(t),\\
		\partial_{x}U_{j}(t,0)&=\partial_{x}U_{j}(t,1)=0.
	\end{split}\right.
\eqs
As a consequence, one necessarily has
\bqs
\forall t>0,~j\in\Z, \quad U_j(t,x)=V_j(t).
\eqs
Integrating \eqref{eqeps} from $x=0$ to $x=1$ and using the Robin boundary conditions \eqref{bdeps}, one also finds
\bqs
\epsilon \frac{\d}{\d t} \int_0^1 v_j^\epsilon(t,x)\d x= \partial_x v_j^\epsilon(t,1)-\partial_x v_j^\epsilon(t,0) = \epsilon \beta (\rho_j^\epsilon(t)+\rho_{j+1}^\epsilon(t))-\alpha\epsilon(v_j^\epsilon(t,0)+v_j^\epsilon(t,1)),
\eqs
from which we get
\bqs
\forall t>0,~j\in\Z, \quad  V_j'(t)=-2\alpha V_j(t)+\beta(P_j(t)+P_{j+1}(t)).
\eqs
Let us finally remark that both sequences $\mathbf{v}^\epsilon$ and $\boldsymbol{\rho}^\epsilon$ are actually $\mathscr{C}^{0,\nu}$ for any $\nu\in(0,1]$, and uniformly bounded, on any time interval $[0,\tau]$ with $\tau>0$, as a consequence, the limiting sequences $\mathbf{U}$ and $\mathbf{P}$ are then also continuous on any time interval $[0,\tau]$. Thus, by definition of the sequences $\mathbf{V}^0$ and $\mathbf{P}^0$, we also have
\bqs
\forall j\in\Z, \quad V_j(0)=V_j^0 \quad \text{ and } \quad P_j(0)=P_j^0.
\eqs
As a consequence, we have obtained the following result.

\begin{theorem}\label{thmLD}Let $\epsilon_0>0$. For any initial sequences $(\mathbf{h}^\epsilon)_{0<\epsilon\leq\epsilon_0}$ and $(\boldsymbol{\Lambda}^\epsilon )_{0<\epsilon\leq\epsilon_0}$, with $\mathbf{h}^\epsilon\in\mathcal{X}^4$ and $\mathbf{h}^\epsilon\in\ell^\infty(\Z)$ for all $\epsilon\in(0,\epsilon_0]$ and satisfying the compatibility condition \eqref{compatibility} together with the assumptions \eqref{ICeps} and \eqref{IClim}, the corresponding unique global classical positive solution $(\mathbf{v}^\epsilon,\boldsymbol{\rho}^\epsilon)$ satisfies
\bqs
\forall x\in[0,1], \quad \underset{\epsilon\rightarrow0}{\lim}~ v_j^\epsilon(t,x)=V_j(t), \quad \underset{\epsilon\rightarrow0}{\lim}~ \rho_j^\epsilon(t)=P_j(t), 
\eqs
locally uniformly in $(t,j)\in[0,+\infty)\times\Z$, wherein $(\mathbf{V},\mathbf{P})$ is solution of the asymptotic system
\bqq\label{asympsyst}
\forall t>0,~j\in\Z, \quad \left\{
	\begin{split}
		 V_j'(t)&=-2\alpha V_j(t)+\beta(P_j(t)+P_{j+1}(t)),\\
		P_{j}^{\prime}(t)&=f(P_{j}(t))+\alpha(V_{j}(t)+V_{j-1}(t))-2\beta P_{j}(t),
	\end{split}\right.
\eqq
with initial condition $V_j(0)=V_j^0$ and $P_j(0)=P_j^0$, $j\in\Z$, defined in \eqref{IClim}.
\end{theorem}

\subsection{Spreading properties of the asymptotic limiting system}

In the following, we focus on the study of the long time behavior of system \eqref{asympsyst} and its spreading properties. For that purpose, we first start by giving the notion of super and sub-solutions and prove a comparison principle.

We say that $(\overline{\mathbf{V}},\overline{\mathbf{P}})$ is a supersolution to \eqref{asympsyst} if for all $j\in\Z$ one has $\overline{V}_j,\overline{P}_j\in\mathscr{C}^1([0,+\infty),\R)$ which satisfy
\bqs
\forall t>0,~j\in\Z, \quad \left\{
	\begin{split}
		 \overline{V}_j'(t)& \geq -2\alpha \overline{V}_j(t)+\beta(\overline{P}_j(t)+\overline{P}_{j+1}(t)),\\
		\overline{P}_{j}^{\prime}(t)&\geq f(\overline{P}_{j}(t))+\alpha(\overline{V}_{j}(t)+\overline{V}_{j-1}(t))-2\beta \overline{P}_{j}(t).
	\end{split}\right.
\eqs
We similarly define a subsolution $(\underline{\mathbf{V}},\underline{\mathbf{P}})$ to \eqref{asympsyst} with the same regularity and all the above inequalities being reversed. We can now state a comparison principle for \eqref{asympsyst} whose proof is a direct consequence of Proposition~\ref{propCPappWQ}.

\begin{proposition}\label{propCPasym}
Let $(\underline{\mathbf{V}},\underline{\mathbf{P}})$ and $(\overline{\mathbf{V}},\overline{\mathbf{P}})$ be respectively a subsolution and supersolution to \eqref{asympsyst}. If we assume that $(\underline{\mathbf{V}},\underline{\mathbf{P}})$ and $(\overline{\mathbf{V}},\overline{\mathbf{P}})$ are locally bounded in time and satisfy for all $j\in\Z$ that $\underline{V}_j(0)\leq\overline{V}_j(0)$ and $\underline{P}_j(0)\leq\overline{P}_j(0)$, then we have $\underline{V}_j(t)\leq\overline{V}_j(t)$ and $\underline{P}_j(t)\leq\overline{P}_j(t)$ for all $t>0$ and  $j\in\Z$. Furthermore, if $(\underline{\mathbf{V}}(0),\underline{\mathbf{P}}(0))\not\equiv(\overline{\mathbf{V}}(0),\overline{\mathbf{P}}(0))$, then we have $\underline{V}_j(t)<\overline{V}_j(t)$ and $\underline{P}_j(t)<\overline{P}_j(t)$ for all $t>0$ and  $j\in\Z$.
\end{proposition}

A direct consequence of the above comparison principle is the uniqueness of bounded solutions of system \eqref{asympsyst}. More generally, for each nontrivial nonnegative initial condition $(\mathbf{V}^0,\mathbf{P}^0)\in\ell^\infty(\Z)\times\ell^\infty(\Z)$, there exists a unique classical global solution $(\mathbf{V},\mathbf{P})$ of system \eqref{asympsyst} with $(\mathbf{V}(0),\mathbf{P}(0))=(\mathbf{V}^0,\mathbf{P}^0)$ such that $V_j,P_j\in\mathscr{C}^1([0,+\infty),\R)$ for all $j\in\Z$, together with uniform bounds 
\bqs
\forall t>0, \quad \forall j\in\Z, \quad 0<V_j(t)\leq \max\left( \left\|\mathbf{V}^0\right\|_{\ell^\infty(\Z)}, \frac{\beta}{\alpha} \right) \quad \text{ and } \quad 0<P_j(t)\leq \max\left( \left\|\mathbf{P}^0\right\|_{\ell^\infty(\Z)}, 1 \right).
\eqs

Regarding the long time behavior of the solutions of system \eqref{asympsyst}, we have the following result which mirrors Theorem~\ref{long}.

\begin{proposition}\label{propinfasym}
	Let $(\mathbf{V},\mathbf{P})$ be the unique global classical solution of \eqref{asympsyst} starting from a nontrivial nonnegative bounded initial sequence $(\mathbf{V}^0,\mathbf{P}^0)\in\ell^\infty(\Z)\times\ell^\infty(\Z)$. Then,
	$$
	\lim_{t\to+\infty}(V_{j}(t),P_{j}(t))=\left(\frac{\beta}{\alpha},1\right), 
	$$
	locally uniformly in $j\in\Z$.
\end{proposition}

\begin{Proof}
Stationary solutions of system \eqref{asympsyst} satisfy
\bqs
\forall j\in\Z, \quad \left\{
	\begin{split}
		0&=-2\alpha V_j+\beta(P_j+P_{j+1}),\\
		0&=f(P_{j})+\alpha(V_{j}+V_{j-1})-2\beta P_{j},
		\end{split}\right.
\eqs
from which we deduce that $V_j=\frac{\beta}{2\alpha}(P_j+P_{j+1})$ and
\bqs
 0 = f(P_j)+\frac{\beta}{2}\left(P_{j-1}-2P_j+P_{j+1}\right),
\eqs
for all $j\in\Z$. As a consequence, from the proof of Theorem~\ref{thm-sta}, we deduce that $(V_j,P_j)=(\beta/\alpha,1)$, for all $j\in\Z$, is the only positive stationary solution to \eqref{asympsyst}.

We now let $N_0>1$ be large enough such that 
\bqs
\beta\left(1-\cos\left(\frac{\pi}{N+1}\right)\right)<f'(0),
\eqs
for all $N\geq N_0$. We then define $\underline{\mathbf{P}}=(\underline{P}_j)_{j\in\Z}$ as
\bqs
\underline{P}_j:=\left\{\begin{split}
\sin\left(\frac{j\pi}{N+1}\right), & \quad  j=1,\dots,N,\\
0, & \quad \text{ otherwise,}
\end{split}\right.
\eqs
and set $\underline{\mathbf{V}}=(\underline{V}_j)_{j\in\Z}$ with
\bqs
\underline{V}_j:=\left\{\begin{split}
\frac{\beta}{2\alpha}\left(\sin\left(\frac{j\pi}{N+1}\right)+\sin\left(\frac{(j+1)\pi}{N+1}\right)\right), & \quad  j=1,\dots,N-1,\\
0, & \quad \text{ otherwise.}
\end{split}\right.
\eqs
By definition, one has $\underline{V}_j=\frac{\beta}{2\alpha}(\underline{P}_j+\underline{P}_{j+1})$ for all $j\in\Z$. As a consequence, there exists $\nu_0>0$ such that $(\nu \underline{\mathbf{V}},\nu\underline{\mathbf{P}})$ is a compactly supported stationary subsolution for all $N\geq N_0$ and $\nu\in(0,\nu_0]$. One can also easily check that
\bqs
 \forall j\in\Z, \quad \overline{V}_j= \max\left( \left\|\mathbf{V}^0\right\|_{\ell^\infty(\Z)}, \frac{\beta}{\alpha} \right) \quad \text{ and } \quad \overline{P}_j= \max\left( \left\|\mathbf{P}^0\right\|_{\ell^\infty(\Z)}, 1 \right),
 \eqs
 gives a stationary supersolution. One can then adapt the arguments of the proof of Theorem~\ref{long} to obtain the local uniform asymptotic convergence of the solutions towards the unique positive stationary solution.
\end{Proof}

Linearizing \eqref{asympsyst} around the trivial steady state, we obtain the following linear system
\bqq\label{linVP}
\forall t>0,~j\in\Z, \quad \left\{
	\begin{split}
		 V_j'(t)&=-2\alpha V_j(t)+\beta(P_j(t)+P_{j+1}(t)),\\
		P_{j}^{\prime}(t)&=(f'(0)-2\beta)P_{j}(t)+\alpha(V_{j}(t)+V_{j-1}(t)).
	\end{split}\right.
\eqq
We  look for exponential solutions of the form
\begin{equation*}
	(V_{j}(t),P_{j}(t))=\ee^{-\mu(j-ct)}(v_0,p_0),
\end{equation*}
where $v_0>0, p_0>0$, $c>0$ and $\mu>0$ to be determined later. Substituting this ansatz into the linear system, we obtain that
 \begin{equation*}
 	\begin{cases}
 		( \mu c+2\alpha)v_{0} -\beta(1+\ee^{-\mu})p_0=0,\\
 		(\mu c+2\beta-f^{\prime}(0))p_{0}-\alpha (1+\ee^{\mu})v_0=0,
 	\end{cases}
 \end{equation*}
which implies that
\begin{equation*}
	\begin{pmatrix}
		 \mu c+2\alpha & -\beta(1+\ee^{-\mu})\\
		-\alpha (1+\ee^{\mu})   & \mu c+2\beta-f^{\prime}(0)\\
	\end{pmatrix}
	\begin{pmatrix}
		v_0\\
		p_0
	\end{pmatrix}
	=
	\begin{pmatrix}
		0\\
		0
	\end{pmatrix}
	.
\end{equation*}
Since we are interested in nontrivial solutions, we must have that
\bqs
	\det
	\begin{pmatrix}
		 \mu c+2\alpha & -\beta(1+\ee^{-\mu})\\
		-\alpha(1+\ee^{\mu})   & \mu c+2\beta-f^{\prime}(0)\\
	\end{pmatrix}
	=0,
\eqs
which also reads
\bqs
(\mu c)^2+(2\beta-f^{\prime}(0)+2\alpha)\mu c+2\alpha(2\beta-f^{\prime}(0))-2\alpha\beta(1+\cosh(\mu))=0.
\eqs
It follows from the above equation that
\bqs
c_{\pm}(\mu)=\dfrac{-(2\alpha+2\beta-f^{\prime}(0))\pm\sqrt{\Delta(\mu)}}{2\mu},
\eqs
where 
$$
\Delta(\mu):=(2\alpha-2\beta+f^{\prime}(0))^2+8\alpha\beta(1+\cosh(\mu))>0.
$$
Only retaining the positive root, we define
\begin{equation}\label{c_11}
	c_{*}^\infty:=\min_{\mu>0}~c_{+}(\mu)=\min_{\mu>0}~\dfrac{-(2\alpha+2\beta-f^{\prime}(0))+\sqrt{\Delta(\mu)}}{2\mu}.
\end{equation}

Let us show that $c_{*}^\infty$ is well-defined. Indeed, we consider the following function
\bqs
\Psi(c,\mu):=-(2\alpha+(2\beta-f^{\prime}(0)))+\sqrt{\Delta(\mu)}-2\mu c.
\eqs
By easy calculations, we have
$$
\Psi(c,0)=-(2\alpha+2\beta-f^{\prime}(0))+\sqrt{\Delta(0)}>0,
$$
\bqs
\forall \mu>0,\quad \dfrac{\partial\Psi(c,\mu)}{\partial\mu}\Bigg\vert_{\mu=0}=-2 c<0,\quad\dfrac{\partial\Psi(c,\mu)}{\partial c}=-2\mu<0,
\eqs
and
\begin{align*}
\dfrac{\partial^2\Psi(c,\mu)}{\partial \mu^2}&=\dfrac{4\alpha\beta}{\sqrt{\Delta(\mu)}\Delta(\mu)}\left(\cosh(\mu)\left(2\alpha-2\beta+f^{\prime}(0)\right)^2+8\alpha\beta\cosh(\mu) +4\alpha\beta(1+\cosh^2(\mu))\right)\\
&>0.
\end{align*}
In view of the above properties of the function $\Psi(c,\mu)$, there exists $c_*^\infty>0$ and $\mu_*>0$ such that
$$
\dfrac{\partial\Psi(c,\mu)}{\partial \mu}\Bigg\vert_{(c_{*}^\infty,\mu_{*})}=0\text{  and  }\Psi(c_*^\infty,\mu_*)=0.
$$
Furthermore,
\begin{itemize}
\item [(i)] if $0<c<c_*^\infty$, then $\Psi(c,\mu)>0,\forall\mu>0$,
\item [(ii)] if $c>c_*^\infty$, then the equation $\Psi(c,\mu)=0$ has two positive real roots $\mu_{1}(c),\mu_{2}(c)$ with $0<\mu_{1}(c)<\mu_{*}<\mu_{2}(c)<+\infty$, such that $\Psi(c,\cdot)<0$ in $(\mu_{1}(c),\mu_{2}(c))$ and $\Psi(c,\cdot)>0$ in $(0,\mu_{1}(c))\cup(\mu_{2}(c),+\infty)$.
\end{itemize}
Therefore, $c_{*}^\infty$ is well-defined. It actually characterizes the spreading speed of \eqref{asympsyst} as stated in the theorem below. Let us remark that our numerical evaluation of the linear spreading speed $c_*$,  given by formula $c_*$, as a function of $d$ while all other parameters being kept fixed suggests that
\bqs
c_* \underset{d\rightarrow+\infty}{\longrightarrow}c_*^\infty.
\eqs
We refer to Figure~\ref{fig:Speed}(c) for an illustration. We leave for future work to rigorously prove such a limit.

\begin{theorem}\label{thmASdinf}
	Let $(\mathbf{V},\mathbf{P})$ be the unique global classical solution of \eqref{asympsyst} starting from a nontrivial nonnegative compactly supported initial sequence $(\mathbf{V}^0,\mathbf{P}^0)$ satisfying $(0,0)\lneqq(\mathbf{V}^0,\mathbf{P}^0)\leq (\beta/\alpha,1)$. Then,
\begin{itemize}
\item[(i)] for all $c>c_*^\infty$, we have
\bqs
\underset{t\rightarrow+\infty}{\lim}~\underset{|j|\geq ct}{\sup}~(V_{j}(t),P_{j}(t))=(0,0);
\eqs
\item[(ii)] for all $c\in(0,c_*^\infty)$, we have
 \bqs
\underset{t\rightarrow+\infty}{\lim}~\underset{|j|\leq ct}{\inf}~(V_{j}(t),P_{j}(t))=\left(\frac{\beta}{\alpha},1\right).
\eqs
\end{itemize}
\end{theorem}

\begin{Proof} Let us first prove item (i) of the theorem. We introduce the sequences
\bqs
\forall t\geq0, \quad j\in\Z, \quad \overline{V_j}(t)= \min\left(\vartheta v_0 \ee^{-\mu_*(j-c_*^\infty t)},\frac{\beta}{\alpha}\right) \quad \text{ and } \quad \overline{P_j}(t)= \min\left(\vartheta \ee^{-\mu_*(j-c_*^\infty t)},1\right),
\eqs
with 
\bqs
v_0=\beta\frac{1+\ee^{-\mu_*}}{ \mu_*c_*^\infty+2\alpha}>0. 
\eqs
Here $\vartheta>0$ is chosen large enough such that $V_j^0 \leq \overline{V_j}(0)$ and $P_j^0 \leq \overline{P_j}(0)$ for all $j\in\Z$ which is always possible since $(\mathbf{V}^0,\mathbf{P}^0)$ is assumed to be compactly supported. By construction $(\overline{\mathbf{V}},\overline{\mathbf{P}})$ is a supersolution of system \eqref{asympsyst}. Thus if $(\mathbf{V},\mathbf{P})$ is the unique global classical solution of \eqref{asympsyst} starting from the nontrivial nonnegative compactly supported initial sequence $(\mathbf{V}^0,\mathbf{P}^0)$ then one has
\bqs
\forall t\geq0, \quad j\in\Z, \quad V_j(t)\leq \overline{V_j}(t) \quad \text{ and } \quad P_j(t)\leq \overline{P_j}(t),
\eqs
from which we readily deduce that for all $c>c_*^\infty$
\bqs
\underset{t\rightarrow+\infty}{\lim}~\underset{j \geq ct}{\sup}~ V_j(t)\leq \underset{t\rightarrow+\infty}{\lim}~\underset{j \geq ct}{\sup}~ \overline{V_j}(t) =0,
\eqs
and
\bqs
\underset{t\rightarrow+\infty}{\lim}~\underset{j \geq ct}{\sup}~ P_j(t)\leq \underset{t\rightarrow+\infty}{\lim}~\underset{j \geq ct}{\sup}~ \overline{P_j}(t) =0.
\eqs
By symmetry, we obtain a similar result for all $j\leq -ct$ which concludes the proof of the first part of the theorem.

The second step of the proof is to devise a compactly supported subsolution whose support moves with speed $c$ close to $c_*^\infty$. So let $c\in(0,c_*^\infty)$ and consider  the linear system
\bqs
\forall t>0,~j\in\Z, \quad \left\{
	\begin{split}
		 V_j'(t)&=-2\alpha V_j(t)+\beta(P_j(t)+P_{j+1}(t)),\\
		P_{j}^{\prime}(t)&=(f'(0)-2\beta)P_{j}(t)+\alpha(V_{j}(t)+V_{j-1}(t)).
	\end{split}\right.
\eqs
Looking once again at exponential solutions of the form
\begin{equation*}
	(V_{j}(t),P_{j}(t))=\ee^{-\mu(j-ct)}(v_0,p_0),
\end{equation*}
with $\mu,v_0,p_0\in\C$ and $c>0$, we see, from the above discussion, that $\mu$ and $c$ should satisfy $\Psi(c,\mu)=0$. We recall  that 
$$
\Psi(c_*^\infty,\mu_{*})=0,\quad \partial_{\mu}\Psi(c_{*}^\infty,\mu_{*})=0,\text{ and }2a:=\partial_{\mu\mu}\Psi(c_{*}^\infty,\mu_{*})>0.
$$
In addition, we also have
$$
\partial_c\Psi(c_*^\infty,\mu_{*})=-2\mu_*<0,\quad \partial_{c\mu}\Psi(c_*^\infty,\mu_{*})=-2<0.
$$
We then consider a neighborhood of $(c_*^\infty,\mu_{*})$, thus we set
$$
\xi:=c_*^\infty-c,\quad \tau:=\mu-\mu_{*}.
$$
The equation $\Psi(c,\mu)=0$ becomes, for $(c,\mu)$ in a neighborhood of $(c_*^\infty,\mu_{*})$:
\begin{equation}\label{tau}
a\tau^2+2\xi\tau+2\mu_*\xi=\phi(\tau,\xi),
\end{equation}
where $\phi(\tau,\xi)$ is analytic in $\tau$ in a neighborhood of $0$, vanishing at $(0,0)$ like $|\tau|^3+|\xi|^2$. For small $\xi>0$, the equation $a\tau^2+2\xi\tau+2\mu_*\xi=0$ has two complex roots
$$
\tau_{\pm}(\xi)=\pm i\sqrt{\dfrac{2\mu_*}{a}\xi}+O(\xi).
$$
By applying Rouch\'e's theorem, we find that equation \eqref{tau} has also two complex roots, which are complex conjugates up to order $\xi$, and are denoted by $\widetilde{\tau}_{\pm}$. These roots satisfy $\widetilde{\tau}_{\pm}(\xi)=\pm i\sqrt{\frac{2\mu_*}{a}\xi}+O(\xi).$ Reverting to the full notation, we observe that for $c$ strictly less than and sufficiently close to $c_*^\infty$, the equation $\Psi(c,\mu)=0$ admits a solution $\mu$ with the following properties: its real part is $\xi$-close to $\mu_{*}$, and hence positive; moreover, it has a nonzero imaginary part of order $\xi^{1/2}$. Setting $p_0=1$, we get that
\bqs
v_0= \frac{\beta (1+\ee^{-\mu})}{ \mu c+2\alpha},
\eqs
and by fixing $\mu=\mu_*+\widetilde{\tau}_{+}(\xi)$, we infer that
\bqs
\Re(v_0)>0, \quad \Im(v_0)<0, \quad \text{ and } \quad \left|\Arg(v_0)\right|=\mathcal{O}(\sqrt{\xi}),
\eqs
where we denoted by $\Arg(v_0)\in(-\pi,\pi]$ the principal argument of $v_0$.
Taking the real parts of the constructed exponential solutions, we set
\bqs
\forall t\geq0, \quad j\in\Z, \quad \left\{\begin{split}
\widetilde{V}_j(t)&:=\Re(V_j(t))=|v_0|\ee^{-\Re(\mu)(j-ct)}\cos\left(\Im(\mu)(j-ct)-\Arg(v_0)\right),\\
\widetilde{P}_j(t)&:=\Re(P_j(t))=\ee^{-\Re(\mu)(j-ct)}\cos\left(\Im(\mu)(j-ct)\right).
\end{split}
\right.
\eqs
In order to obtain compactly supported subsolutions, we truncate the above solutions as follows. We define 
\bqs
\Omega_V(t):=\left\{x\in\R~|~ct-\frac{\pi}{2\Im(\mu)}+\frac{\Arg(v_0)}{\Im(\mu)}\leq x \leq  ct+\frac{\pi}{2\Im(\mu)}+\frac{\Arg(v_0)}{\Im(\mu)}\right\},
\eqs
and
\bqs
\Omega_P(t):=\left\{x\in\R~|~ct-\frac{\pi}{2\Im(\mu)}\leq x \leq  ct+\frac{\pi}{2\Im(\mu)}\right\},
\eqs
and set
\bqq\label{subsolVP}
\forall t\geq0, \quad \underline{V}_j(t):=\left\{\begin{split}
\widetilde{V}_j(t),  & \quad j\in \Omega_V(t),\\
0, & \quad \text{ otherwise},
\end{split}\right.
\quad \quad
\underline{P}_j(t):=\left\{\begin{split}
\widetilde{P}_j(t),  & \quad j\in \Omega_P(t),\\
0, & \quad \text{ otherwise}.
\end{split}\right.
\eqq
Since both $\Im(\mu)=O(\sqrt{\xi})$ and $\left|\Arg(v_0)\right|=O(\sqrt{\xi})$, we thus have that $\frac{\left|\Arg(v_0)\right|}{\Im(\mu)}=O(1)$. We readily have that when $j\in \Omega_V(t) \cap \Omega_P(t)$ or $j\in\Z \backslash \Omega_V(t) \cap \Omega_P(t)$, then $(\underline{V}_j(t),\underline{P}_j(t))$ is a solution of the linear system \eqref{linVP}. On the other hand, if $j\in \Omega_V(t) \backslash \Omega_P(t)$, then
\begin{align*}
 \underline{V}_j'(t)+2\alpha \underline{V}_j(t) - \beta (\underline{P}_j(t)+\underline{P}_{j+1}(t))&= \widetilde{V}_j'(t)+2\alpha \widetilde{V}_j(t)=\beta (\widetilde{P}_j(t)+\widetilde{P}_{j+1}(t)) \leq 0,\\
\underline{P}_j'(t)+(2\beta-f'(0))\underline{P}_j(t)-\alpha(\underline{V}_j(t)+\underline{V}_{j-1}(t))&=-\alpha(\widetilde{V}_j(t)+\widetilde{V}_{j-1}(t)) \leq0,
\end{align*}
upon taking $\xi$ small enough to ensure that both $\widetilde{P}_j(t)<0$ and $\widetilde{P}_{j+1}(t)\leq 0$. A similar argument with $j\in \Omega_P(t) \backslash \Omega_V(t)$ shows that $(\underline{V}_j(t),\underline{P}_j(t))$ is a subsolution of \eqref{linVP} for all $t\geq0$ and $j\in\Z$.

At the moment, we have only constructed a compactly supported subsolution for the linear system \eqref{linVP}. It is not difficult to check that all the above arguments naturally perturb if instead we consider the modified linear system 
\bqq\label{linVPdelta}
\forall t>0,~j\in\Z, \quad \left\{
	\begin{split}
		 V_j'(t)&=-2\alpha V_j(t)+\beta(P_j(t)+P_{j+1}(t)),\\
		P_{j}^{\prime}(t)&=(f'(0)-2\beta-\delta)P_{j}(t)+\alpha(V_{j}(t)+V_{j-1}(t)),
	\end{split}\right.
\eqq
for some small $\delta>0$. More precisely, there exists $\delta_0>0$ such that for any $\delta\in(0,\delta_0]$, one can construct a compactly supported subsolution $(\underline{\mathbf{V}}^\delta,\underline{\mathbf{P}}^\delta)=(\underline{V}_j^\delta,\underline{P}_j^\delta)_{j\in\Z}$ in the form of \eqref{subsolVP}. Furthermore, let $\iota_0>0$ be such that
\bqs
(f'(0)-\delta_0)u\leq f(u), \quad 0\leq u \leq \iota_0.
\eqs
Then one can find $\eta_0>0$ small enough such that $\eta_0 \underline{P}_j^{\delta_0}(t)\leq \iota_0$ for all $t>0$ and $j\in\Z$. As a consequence $(\eta_0\underline{\mathbf{V}}^{\delta_0},\eta_0\underline{\mathbf{P}}^{\delta_0})$ is a compactly supported subsolution to the nonlinear system \eqref{asympsyst}.

We can now prove item (ii) of the theorem. Let $c\in(0,c_*^\infty)$ and choose $c'\in(c,c^\infty_*)$ very close to $c_*^\infty$. From the positivity of the solution of the nonlinear system \eqref{asympsyst}, upon eventually decreasing the size of $\eta_0>0$, we can always ensure that at time $t=1$ one has
\bqs
(\mathbf{V}(1),\mathbf{P}(1))\geq (\eta_0\underline{\mathbf{V}}^{\delta_0}(0),\eta_0\underline{\mathbf{P}}^{\delta_0}(0)),
\eqs
where $(\eta_0\underline{\mathbf{V}}^{\delta_0},\eta_0\underline{\mathbf{P}}^{\delta_0})$ is the compactly supported subsolution associated to the speed $c'$ constructed in the previous step. From the comparison principle of Proposition~\ref{propCPasym} we obtain that
\bqs
\forall t\geq1, \quad (\mathbf{V}(t),\mathbf{P}(t))\geq (\eta_0\underline{\mathbf{V}}^{\delta_0}(t-1),\eta_0\underline{\mathbf{P}}^{\delta_0}(t-1)).
\eqs
There exists $\nu\in(0,1)$ small, depending on $c'$, such that
\bqs
V_{\lfloor c't \rfloor}(t)\geq \eta_0 \underline{V}_{\lfloor c't \rfloor}^\delta(t-1) \geq \frac{\beta}{\alpha} \nu \quad \text{ and } \quad V_{\lfloor c't \rfloor+1}(t)\geq \eta_0 \underline{V}_{\lfloor c't \rfloor+1}^\delta(t-1) \geq \frac{\beta}{\alpha} \nu,
\eqs
with
\bqs
P_{\lfloor c't \rfloor}(t)\geq \eta_0 \underline{P}_{\lfloor c't \rfloor}^\delta(t-1) \geq  \nu \quad \text{ and } \quad P_{\lfloor c't \rfloor+1}(t)\geq \eta_0 \underline{P}_{\lfloor c't \rfloor+1}^\delta(t-1) \geq \nu,
\eqs
for all $t\geq1$. By a symmetry argument, we also obtain that
\bqs
V_{-\lfloor c't \rfloor}(t)\geq \frac{\beta}{\alpha} \nu \quad \text{ and } \quad V_{-\lfloor c't \rfloor-1}(t)\geq  \frac{\beta}{\alpha} \nu,
\eqs
with
\bqs
P_{-\lfloor c't \rfloor}(t)\geq  \nu \quad \text{ and } \quad P_{-\lfloor c't \rfloor-1}(t)\geq \nu,
\eqs
for all $t\geq1$. Upon even reducing the size of $\nu$, by positivity of the solution $(\mathbf{V},\mathbf{P})$, we can always ensure that
\bqs
V_j(1)\geq \frac{\beta}{\alpha}\nu\quad\text{ and }\quad P_j(1)\geq\nu\quad\text{ for all}\quad-c'-1\leq j \leq  c'+1
\eqs
Since $\left(\frac{\beta}{\alpha} \nu,\nu\right)_{j\in\Z}$ is a homogeneous subsolution of \eqref{asympsyst}, we are thus in a position to apply the comparison principle of Proposition~\ref{propCPappWQ2b} with two moving boundaries given by $\zeta(t)=-c't$ and $\xi(t)=c't$. It implies that
\bqs
\forall t\geq1, \quad \underset{|j|\leq ct}{\inf}~V_j(t) \geq  \frac{\beta}{\alpha} \nu \quad \text{ and } \quad \underset{|j|\leq ct}{\inf}~P_j(t) \geq \nu,
\eqs
from which we deduce that
\bqs
\underset{t\rightarrow+\infty}{\liminf} \underset{|j|\leq ct}{\inf}\left(V_j(t),P_j(t)\right)\geq \underset{t\rightarrow+\infty}{\liminf} \underset{|j|\leq c't}{\inf}\left(V_j(t),P_j(t)\right)\geq \left(\frac{\beta}{\alpha},1\right).
\eqs
But from Proposition~\ref{propinfasym}, we have
\bqs
\underset{t\rightarrow+\infty}{\limsup} \underset{|j|\leq ct}{\inf}\left(V_j(t),P_j(t)\right)\leq \left(\frac{\beta}{\alpha},1\right).
\eqs
As a conclusion, we have proved that
 \bqs
\underset{t\rightarrow+\infty}{\lim}~\underset{|j|\leq ct}{\inf}~(V_{j}(t),P_{j}(t))=\left(\frac{\beta}{\alpha},1\right),
\eqs
for all $c\in(0,c_*^\infty)$. This concludes the proof of the theorem.
\end{Proof}

\section{Discussion}

\paragraph{Summary of main results.} In this work, we have proposed a new model to describe biological invasions constrained on infinite homogeneous one dimensional metric graphs. Our model consists of an infinite PDE-ODE system where, at each vertex of the one-dimensional lattice $\Z$, we have a standard logistic equation and connections between vertices are given by diffusion equations on the edges supplemented with Robin like boundary conditions at the vertices. Our first main result is the existence and uniqueness of classical, global in time, positive bounded solutions of our PDE-ODE model. Our second main result is the characterization of the long time behavior of the unique solution of our model, where we prove local uniform convergence towards the unique positive bounded stationary solution of the system. Next, we analyzed the linearized problem around the trivial constant state and derived a theoretical formula for the linear spreading speed of our model, defined as the smallest possible speed for which there exist exponential solutions with prescribed form. We then proved that this linear spreading speed is actually the asymptotic spreading speed of the full nonlinear model, which constitutes the key result of our present study. Finally, we investigated the large diffusion limit of the model and established the convergence towards an asymptotic system  for which we also managed to fully characterize its asymptotic spreading properties. We also illustrated our theoretical findings with a selection of numerical simulations.

\begin{figure}[t!]
\centering
\subfigure[$\lambda=0$.]{\includegraphics[width=.31\textwidth]{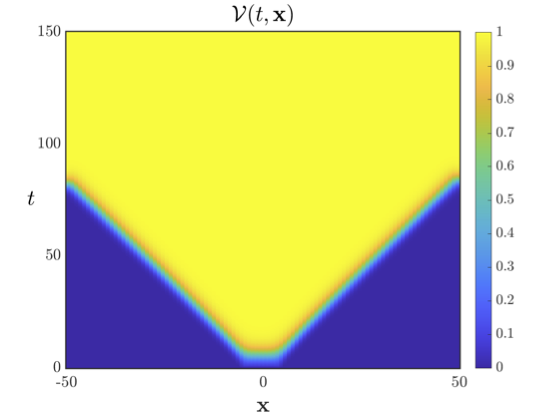}}
\subfigure[$\lambda=2$.]{\includegraphics[width=.31\textwidth]{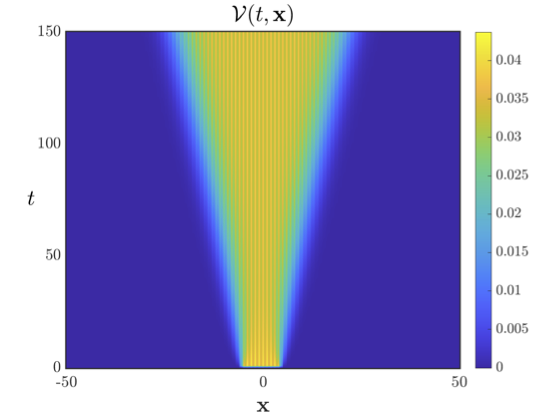}}
\subfigure[$\lambda=4$.]{\includegraphics[width=.32\textwidth]{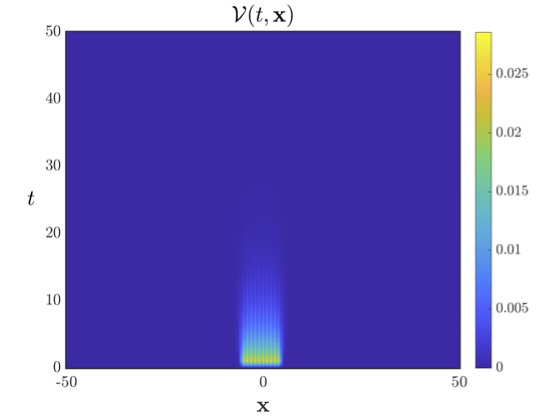}}
\caption{(a) Space-time plot of $\mathcal{V}(t,\mathbf{x})$ solution of \eqref{eq1}-\eqref{bd} starting from a compactly supported centered initial condition. (b)-(c) Space-time plots of the solution in the presence of a hostile environment (see equation \eqref{eqhostile}) starting from the same initial condition. For intermediate values of $\lambda$ (panel (b)), the solution behaves similarly as in the case $\lambda=0$ while for larger values (panel (c)), the solution goes extinct uniformly in space, indicating a possible threshold behavior at some critical value $\lambda_c$. Other values of the parameters are set to $d=\alpha=\beta=1$ and $f(u)=u(1-u)$.}
 \label{fig:STPlambda}
\end{figure}

\paragraph{Natural extensions.} From a biological point of view, it could be interesting to consider several extensions of the model. First of all, roads could be modeled as a hostile environment  such that the diffusion equation of \eqref{eq1} could be replaced by
\bqq\label{eqhostile}
\partial_t v_j=d\partial_x^2 v_j-\lambda v_j,
\eqq
for some $\lambda>0$ representing a death rate on the road. Such a modeling assumption has already been proposed for other reaction-diffusion models \cite{BRR16}. We expect that the presence of a hostile environment will have a direct effect on the stationary solutions of the model and thus on the long time behavior of the solutions. More precisely, we anticipate a threshold effect and the existence of a critical value for $\lambda$ (depending on all other parameters of the model), above which the only stationary solution is the trivial constant steady state, and below which there exists a unique bounded positive stationary solution. For values of $\lambda$ above this critical parameter, solutions of the Cauchy problem are expected to uniformly converge to the trivial constant steady state, and thus go extinct, reflecting the fact that the road is too hostile for the population to survive. This threshold phenomena is numerically confirmed in the case of the logistic nonlinearity as illustrated in Figure~\ref{fig:STPlambda} where we report space-time plots of the solution of \eqref{eq1}-\eqref{bd} with the first equation of \eqref{eq1} modified into \eqref{eqhostile}. For our visualisation, we used the following notation:
\bqs
\text{ for }\quad\mathbf{x}\in [j,j+1], \quad \mathcal{V}(t,\mathbf{x})=v_j(t,x),\quad\text{ with }\quad\mathbf{x} = j+x,
\eqs
with $\mathcal{V}(t,\cdot)$ being defined on $\R$.

\begin{figure}[t!]
\centering
\subfigure[$\nu_+>\nu_->0$.]{\includegraphics[width=.31\textwidth]{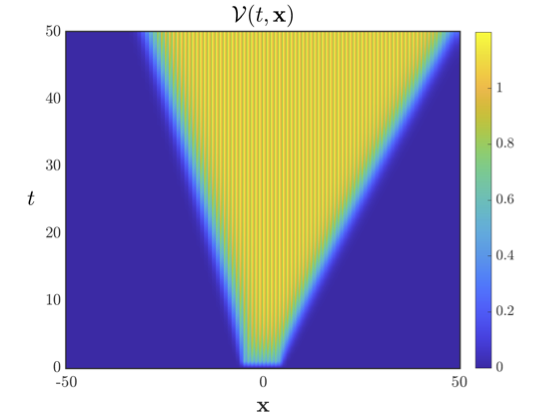}}
\subfigure[$\nu_+=\nu_->0$.]{\includegraphics[width=.31\textwidth]{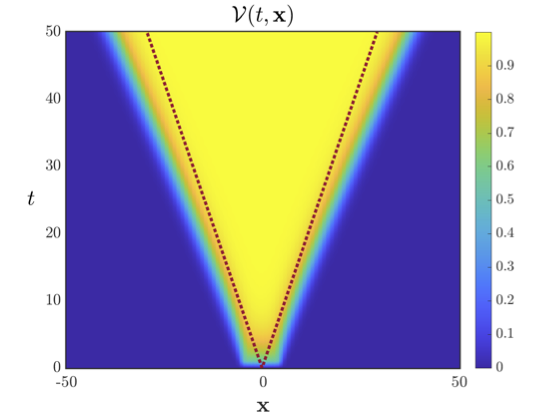}}
\subfigure[$0<\nu_+<\nu_-$.]{\includegraphics[width=.32\textwidth]{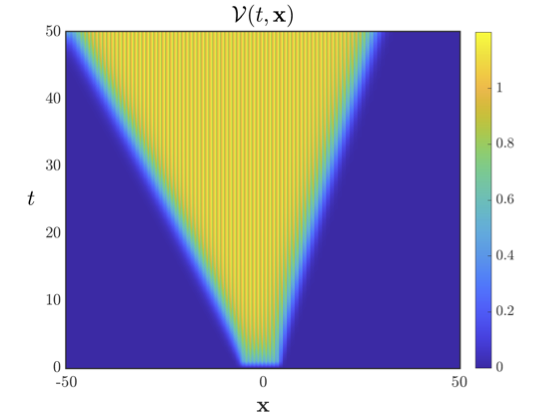}}
\caption{Space-time plot of $\mathcal{V}(t,\mathbf{x})$ solution of \eqref{eq1}-\eqref{bcnu} starting from a compactly supported centered initial condition with logistic nonlinearity $f(u)=u(1-u)$. The dashed dark red lines in panel (b) indicate the position of the level sets $\mathcal{V}(t,\mathbf{x})=1/2$ of the solution when $\nu_+=\nu_-=0$. (a) $(\nu_+,\nu_-)=(3/2,1/2)$. (b) $(\nu_+,\nu_-)=(1,1)$. (c) $(\nu_+,\nu_-)=(1/2,3/2)$. Other values of the parameters are set to $d=\alpha=\beta=1$.}
 \label{fig:STPnu}
\end{figure}

As explained in the introduction, for simplicity, our model neglects the possibility that individuals could pass from one road to an adjacent one. Assuming that such exchanges are homogeneous, in the sense that they are independent of the position $j$, and if $\nu_+>0$ and $\nu_->0$  denote the corresponding exchanges from passing from left to right and right to left respectively, then the Robin boundary conditions \eqref{bd} should be modified according to
\bqq\label{bcnu}
\begin{cases}
-d\partial_{x}v_{j}(t,0)+\alpha v_{j}(t,0)=\beta \rho_{j}(t)+\nu_+v_{j-1}(t,\ell)-\nu_-v_{j}(t,0),\\
d\partial_{x}v_{j}(t,\ell)+\alpha v_{j}(t,\ell)=\beta \rho_{j+1}(t)+\nu_- v_{j+1}(t,0)-\nu_+ v_{j}(t,\ell).
\end{cases}
\eqq
These new exchange terms typically account for the permeability of cell membranes in gap junction models \cite{KS09}. We anticipate that the presence of these new exchange terms will introduce a bias in the direction of propagation. More precisely, we expect that when $\nu_+>\nu_->0$ there should be a faster rightward propagation and a slower leftward propagation while when $0<\nu_+<\nu_-$ there should be a faster leftward propagation and a slower rightward propagation. This is observed numerically in Figure~\ref{fig:STPnu}, panels (a) and (c) for the logistic nonlinearity. In the case of balanced exchanges $\nu_+=\nu_->0$ leftward and rightward propagation occur at the same speed, but with a larger spreading speed than in the original model. This is illustrated in panel (b) of Figure~\ref{fig:STPnu} where the dashed dark red lines indicate the position of the level sets $\mathcal{V}(t,\mathbf{x})=1/2$ of the solution when $\nu_+=\nu_-=0$. The effects of such exchange terms can further be understood in the large diffusion limit $d\rightarrow+\infty$.  Formally, the asymptotic system is expected to read:
\bqs
 \left\{
	\begin{split}
		 V_j'(t)&=-2\alpha V_j(t)+\beta(P_j(t)+P_{j+1}(t))+\nu_-V_{j+1}(t)-(\nu_-+\nu_+)V_j(t)+\nu_+V_{j-1}(t) ,\\
		P_{j}^{\prime}(t)&=f(P_{j}(t))+\alpha(V_{j}(t)+V_{j-1}(t))-2\beta P_{j}(t)
	\end{split}\right.
\eqs
for all $t>0$ and $j\in\Z$. We remark the presence of a new term in the first equation in the form of $\nu_-V_{j+1}-(\nu_-+\nu_+)V_j+\nu_+V_{j-1}$ which acts as an effective advection-diffusion mechanism as can be seen by rewriting it as
\bqs
\nu_-V_{j+1}-(\nu_-+\nu_+)V_j+\nu_+V_{j-1} = \frac{\nu_++\nu_-}{2}\left(V_{j+1}-2V_j+V_{j-1}\right)-\frac{\nu_+-\nu_-}{2}\left(V_{j+1}-V_{j-1}\right).
\eqs
Depending on the strength of advection, we speculate that we may prevent leftward or rightward propagation. We leave the analysis of these natural extensions for future work.

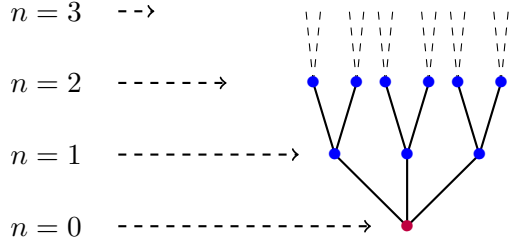
\begin{figure}[t!]
\centering
\begin{tikzpicture}[scale=1.9]
\draw[thick] (0,0) -- (-0.5,0.5);
\draw[thick] (0,0) -- (0.5,0.5);
\draw[thick] (0,0) -- (0,0.5);

\draw[thick] (0.5,0.5) -- (0.35,1);
\draw[thick] (0.5,0.5) -- (0.65,1);
\draw[thick] (-0.5,0.5) -- (-0.35,1);
\draw[thick] (-0.5,0.5) -- (-0.65,1);
\draw[thick] (0,0.5) -- (-0.15,1);
\draw[thick] (0,0.5) -- (0.15,1);

\draw[dashed] (-0.65,1) -- (-0.7,1.5);
\draw[dashed] (-0.65,1) -- (-0.6,1.5);

\draw[dashed] (-0.35,1) -- (-0.4,1.5);
\draw[dashed] (-0.35,1) -- (-0.3,1.5);

\draw[dashed] (-0.15,1) -- (-0.2,1.5);
\draw[dashed] (-0.15,1) -- (-0.1,1.5);

\draw[dashed] (0.65,1) -- (0.7,1.5);
\draw[dashed] (0.65,1) -- (0.6,1.5);

\draw[dashed] (0.35,1) -- (0.4,1.5);
\draw[dashed] (0.35,1) -- (0.3,1.5);

\draw[dashed] (0.15,1) -- (0.2,1.5);
\draw[dashed] (0.15,1) -- (0.1,1.5);

\node[purple] at (0,0) {$\bullet$};
\node[blue] at (-0.5,0.5) {$\bullet$};
\node[blue] at (0.5,0.5) {$\bullet$};
\node[blue] at (0,0.5) {$\bullet$};
\node[blue] at (-0.65,1) {$\bullet$};
\node[blue] at (-0.35,1) {$\bullet$};
\node[blue] at (-0.15,1) {$\bullet$};
\node[blue] at (0.15,1) {$\bullet$};
\node[blue] at (0.35,1) {$\bullet$};
\node[blue] at (0.65,1) {$\bullet$};

\node at (-2.5,0) {$n=0$};
\node at (-2.5,0.5) {$n=1$};
\node at (-2.5,1) {$n=2$};
\node at (-2.5,1.5) {$n=3$};

\draw[-to,thick,dashed] (-2,0) -- (-0.25,0);
\draw[-to,thick,dashed] (-2,0.5) -- (-0.75,0.5);
\draw[-to,thick,dashed] (-2,1) -- (-1.25,1);
\draw[-to,thick,dashed] (-2,1.5) -- (-1.75,1.5);
\end{tikzpicture}
\caption{Homogeneous rooted tree of degree $k=2$. The red node represents the root of the tree. Each generation within the tree is labelled by an integer $n\in\N$.}
\label{fig:tree}
\end{figure}

\paragraph{Beyond the lattice case.} Our model \eqref{eq1}-\eqref{bd} considers the simplest connected metric graph possible: the one dimensional lattice $\Z$. It would be very relevant to extend our framework to other classes of metric graphs. It seems natural to start by considering homogeneous trees, and we already refer to recent developments regarding spreading properties of  reaction-diffusion equations on homogeneous trees \cite{BF21b,Hoffman-Holzer-2019,FHT}. In order to better explain the class of models we have in mind, we introduce some notations. We let $k\in\N$ with $k\geq1$ and shall denote by $\mathbb{T}_k$ a homogeneous tree of degree $k$ with the convention that $\mathbb{T}_1=\Z$. We recall that a homogeneous tree of degree $k$ is an infinite graph where each vertex has precisely $k+1$ adjacent vertices, and we refer to Figure~\ref{fig:tree} for an illustration in the case $k=2$.  As in our original model, we suppose that all edges of the tree have the same length $\ell>0$. To simplify the presentation, we will identify one vertex as being the root of the tree and by convention we will label this vertex as $n=0$ with associated population density $\rho_0(t)$. We will also assume that all populations at some fixed distance away from the root are equal. As a consequence, it will be convenient to denote by $\rho_n(t)$ as a representative population at distance $n\geq1$ from the root. Similarly, we shall also denote by $v_n(t,x)$ a representative population leaving on the edge at distance $n$ from the root. We readily remark that the new model is now indexed by the natural integers $\N$. Only the dynamics for each $\rho_n(t)$ has to be modified according to
\bqs
\left\{
\begin{split}
\rho_0'(t)&=f(\rho_0(t))+(k+1)\left(\alpha v_0(t,0)-\beta \rho_0(t)\right),\\
\rho_n'(t)&=f(\rho_n(t))+\alpha\left(v_{n-1}(t,\ell)+kv_n(t,0)\right)-(k+1)\beta \rho_n(t), \quad n\geq 1.
\end{split}
\right.
\eqs
Let us already remark that the new exchange terms can also be written as
\begin{align*}
\alpha\left(v_{n-1}(t,\ell)+kv_n(t,0)\right)-(k+1)\beta \rho_n(t)&=\alpha\left(v_{n-1}(t,\ell)+v_n(t,0)\right)-2\beta \rho_n(t)\\
&~~~+(k-1)\left( \alpha v_n(t,0)- \beta \rho_n(t)\right),
\end{align*}
where we see the presence of a new term $(k-1)\left( \alpha v_n(t,0)- \beta \rho_n(t)\right)$ which can be interpreted as a drift that may or may not block the propagation within the tree depending on the other parameters of the model. For reaction-diffusion equations set on homogeneous trees \cite{BF21b,Hoffman-Holzer-2019}, the presence of such a term typically prevents propagation within the tree for $k$ large enough and we expect a similar threshold to also happen here. We shall investigate the extension to homogeneous trees in a forthcoming work.

\section*{Acknowledgment}
We would like to thank the two anonymous referees for their valuable comments and suggestions which definitely improved the final version of the manuscript. Authors acknowledge support from ANR project Indyana under grant agreement ANR-21-CE40-0008 and  ANR project ReaCh under grant agreement ANR-23-CE40-0023. 

\section*{Data availability statement} Data sharing is not applicable to this article as no datasets were generated or analyzed during the current study.

\appendix

\section{Representation formula}

We consider the heat equation
\bqq\label{heatapp}
\partial_tv(t,x) = d\partial_x^2 v(t,x),\quad t>0,~ x\in(0,1),\\
\eqq
with inhomogeneous Robin boundary conditions
\bqq\label{robinapp}
\left\{\begin{split}
-d\partial_xv(t,0)+\alpha v(t,0) &= g(t),\\
d\partial_xv(t,1)+\alpha v(t,1) &= h(t),\\
\end{split}\right. \quad t>0,
\eqq
and initial condition
\bqq\label{icapp}
v(0,x)=v_0(x), \quad x\in[0,1].
\eqq
We also define
\bqs
\mathcal{K}(t,x):=\frac{1}{\sqrt{4\pi dt}}\ee^{-\frac{x^2}{4dt}}, \quad t>0, x\in\R.
\eqs
\begin{proposition}\label{proprep}
Assume that $h,g\in \mathscr{C}^0(\R_+,\R)$ and $v_0\in\mathcal{C}^0([0,1],\R)$, then the solution to \eqref{heatapp}-\eqref{robinapp}-\eqref{icapp} can be represented as follows
\bqq\label{repapp}
\begin{split}
v(t,x)&=\int_{0}^1 \K(t,x-y)v_0(y)\d y+\int_0^t \left[ \K(t-s,x-1)h(s)+\K(t-s,x)g(s)\right]\d s \\
&~~~+\int_0^t \left[ -\alpha \K(t-s,x-1) +d \partial_x\K(t-s,x-1)\right] u(s,1)\d s\\
&~~~-\int_0^t \left[ \alpha \K(t-s,x) +d \partial_x\K(t-s,x)\right] u(s,0)\d s,
\end{split}
\eqq
for all $t>0$ and $x\in[0,1]$.
\end{proposition}

\begin{Proof}
Let $w\in \mathscr{C}^2((0,+\infty)\times[0,1],\R)$ and for any $t>0$, from \eqref{heatapp} we have
\bqs
\begin{split}
0&= \int_0^t \int_0^1 w(s,y)\left(\partial_sv(s,y)-d\partial_y^2v(s,y)\right)\d y\d s\\
& = \int_0^1 \left(w(t,y)v(t,y)-w(0,y)v(0,y)\right)\d y -\int_0^t \int_0^1 \left(\partial_sw(s,y)+d\partial_y^2w(s,y)\right)v(s,y)\d y\d s\\
& ~~~ - d \int_0^t \left(w(s,1)\partial_yv(s,1)-w(s,0)\partial_yv(s,0)\right)\d s+d \int_0^t \left( \partial_y w(s,1)v(s,1)-\partial_y w(s,0)v(s,0)\right)\d s,
\end{split}
\eqs
for all $t>0$ and $x\in[0,1]$. We now specify $w$ to
\bqs
w(s,y)=\K(t+\epsilon-s,x-y),
\eqs
for $\epsilon>0$. Note that for all $s\in[0,t]$ and $y\in[0,1]$ it satisfies $\partial_s w+d\partial_y^2w=0$.  We also note that
\bqs
\int_0^1 w(t,y)v(t,y) \d y = \int_0^1 \K(\epsilon,x-y)v(t,y) \d y \underset{\epsilon\rightarrow0}{\longrightarrow} v(t,x),
\eqs
while 
\bqs
\int_0^1 w(0,y)v(0,y) \d y = \int_0^1 \K(t+\epsilon,x-y)v_0(y) \d y \underset{\epsilon\rightarrow0}{\longrightarrow} \int_0^1 \K(t,x-y)v_0(y) \d y.
\eqs
Finally, we simply note that $\partial_y w(s,y)=-\partial_x\K(t+\epsilon-s,x-y)$, and using the Robin boundary condition \eqref{robinapp}, we eventually derive \eqref{repapp}.
\end{Proof}

\section{Comparison principles}

\begin{proposition}\label{propCPapp1}
Let $\mathbf{v}$ and $\boldsymbol{\rho}$ with  $\rho_j\in\mathscr{C}^1([0,+\infty),\R)$ and
\bqs
v_j\in \mathscr{C}^{0}([0,+\infty)\times[0,1],\R),~ \partial_t v_j, \partial_x^2v_j\in\mathscr{C}^{0}((0,+\infty)\times(0,1),\R), \text{ and } \partial_x v_j\in \mathscr{C}^{0}((0,+\infty)\times[0,1],\R),
\eqs
for all $j\in\Z$, which satisfy 
\bqq\label{systvrho}
\begin{cases}
\partial_tv_j(t,x) - d\partial_x^2 v_j(t,x)\geq0,\quad x\in(0,1),\\
\rho_j'(t)-c_j(t)\rho_j(t)  \geq \alpha\left[v_j(t,0)+v_{j-1}(t,1)\right],\\
-d\partial_xv_j(t,0)+\alpha v_j(t,0) \geq \beta \rho_j(t),\\
d\partial_xv_j(t,1)+\alpha v_j(t,1)  \geq \beta \rho_{j+1}(t),
\end{cases}
\eqq
for all $t>0$ and $j\in\Z$ with some $\mathbf{c}=(c_j)_{j\in\Z}\in L^\infty(\R_+,\ell^\infty(\Z))$. Assume that $v_j(0,x)\geq0$ and $\rho_j(0)\geq0$ for all $x\in[0,1]$ and $j\in\Z$, then $v_j(t,x)\geq0$ and $\rho_j(t)\geq 0$ for all $t>0$, $x\in[0,1]$ and $j\in\Z$. If furthermore $\mathbf{v}(0)\not\equiv0$ or $\boldsymbol{\rho}(0)\not\equiv0$, then $v_j(t,x)>0$ and $\rho_j(t)> 0$ for all $t>0$, $x\in[0,1]$ and $j\in\Z$.
\end{proposition}
\begin{Proof}
Fix $T>0$. By assumption on the sequence $\mathbf{c}$, there exists $K>0$ such that 
\bqs
K-c_j(t)>0,\text{ for} ~ t\in(0,T]\text{ and }j\in\Z.
\eqs
For any $\gamma>0$, we define
\bqs\left\{
\begin{split}
w_j(t,x) &:=  \ee^{-\gamma|j|-Kt}v_j(t,x),\\
z_j(t)&:= \ee^{-\gamma|j|-Kt} \rho_j(t).
\end{split}\right.
\eqs
Since $\mathbf{v}$ and $\boldsymbol{\rho}$ are assumed to be locally bounded, we have for each $t\in(0,T]$, $j\in\Z$ and $x\in[0,1]$ that
\bqs\left\{
\begin{split}
w_j(t,x) &\underset{j\rightarrow\pm\infty}{\longrightarrow} 0,\\
z_j(t)&\underset{j\rightarrow\pm\infty}{\longrightarrow} 0.
\end{split}\right.
\eqs
The sequences $\mathbf{w}$ and $\mathbf{z}$ now satisfy
\bqq\label{systwz}
\begin{cases}
\partial_tw_j(t,x) - d\partial_x^2 w_j(t,x)+Kw_j(t,x)\geq0,\quad x\in(0,1),\\
z_j'(t)+(K-c_j(t))z_j(t)  \geq \alpha\left[w_j(t,0)+C^\gamma_j w_{j-1}(t,1)\right],\\
-d\partial_xw_j(t,0)+\alpha w_j(t,0) \geq \beta z_j(t),\\
d\partial_xw_j(t,1)+\alpha w_j(t,1)  \geq \beta C^{-\gamma}_{j+1} z_{j+1}(t),
\end{cases}
\eqq
where the sequence $C^\gamma_j$ is defined as follows
\bqq\label{Cgamma}
C^\gamma_j=\left\{
\begin{split}
\ee^{-\gamma}, & \quad j\geq 1,\\
\ee^{\gamma}, & \quad  j\leq 0.
\end{split}
\right.
\eqq
We now let $\epsilon>0$ and define
\bqs\left\{
\begin{split}
w_j^\epsilon(t,x) &:=w_j(t,x)+\epsilon \ee^{\varrho t +\delta \left(x-\frac{1}{2}\right)^2},\\
z_j^\epsilon (t)&:= z_j(t)+\epsilon \ee^{\varrho t +\frac{\delta}{4}},
\end{split}\right.
\eqs
for two constants $\varrho>0$ and $\delta>0$ that will be fixed later in the proof.
Elementary computations give
\begin{align*}
&\partial_tw_j^\epsilon(t,x) - d\partial_x^2 w_j^\epsilon(t,x)+Kw_j^\epsilon(t,x)\\
&=\partial_tw_j(t,x) - d\partial_x^2 w_j(t,x)+Kw_j(t,x)+\epsilon \left(\varrho+K-2d\delta -4d\delta^2\left(x-\frac{1}{2}\right)^2\right)\ee^{\varrho t +\delta \left(x-\frac{1}{2}\right)^2},
\end{align*}
and 
\begin{align*}
&{z_j^{\epsilon}}'(t)+(K-c_j(t))z_j^\epsilon(t)-\alpha\left[w_j^\epsilon(t,0)+C^\gamma_j w_{j-1}^\epsilon(t,1)\right]\\
&=z_j'(t)+(K-c_j(t))z_j(t)-\alpha\left[w_j(t,0)+C^\gamma_j w_{j-1}(t,1)\right]+\epsilon\left( \varrho+K-c_j(t)-\alpha(1+C^\gamma_j) \right)\ee^{\varrho t +\frac{\delta}{4}},
\end{align*}
together with
\begin{align*}
-d\partial_xw_j^\epsilon(t,0)+\alpha w_j^\epsilon(t,0) - \beta z_j^\epsilon(t)&=-d\partial_xw_j(t,0)+\alpha w_j(t,0) - \beta z_j(t) +\epsilon(d\delta  +\alpha-\beta)\ee^{\varrho t +\frac{\delta}{4}},\\
d\partial_xw_j^\epsilon(t,1)+\alpha w_j^\epsilon(t,1) - C^{-\gamma}_{j+1} \beta z_{j+1}^\epsilon(t)&=d\partial_xw_j(t,1)+\alpha w_j(t,1) - C^{-\gamma}_{j+1} \beta z_{j+1}(t)\\
&~~~ +\epsilon(d\delta  +\alpha-\beta C^{-\gamma}_{j+1})\ee^{\varrho t +\frac{\delta}{4}}.
\end{align*}
As a consequence, we first fix $\delta>0$ such that
\bqs
\delta > \frac{\beta \ee^{\gamma}-\alpha}{d},
\eqs
and then select $\varrho>0$ large enough such that
\bqs
\varrho+K-2d\delta -d\delta^2>0 \text{ and } \varrho-\alpha(1+\ee^{\gamma})>0,
\eqs
that is
\bqs
\varrho>\max\left( 2d\delta +d\delta^2-K, \alpha(1+\ee^{\gamma})\right).
\eqs
With such a choice, the sequences $\mathbf{w}^\epsilon$ and $\mathbf{z}^\epsilon$ now satisfy
\bqq\label{systwzeps}
\begin{cases}
\partial_tw_j^\epsilon(t,x) - d\partial_x^2 w_j^\epsilon(t,x)+Kw_j^\epsilon(t,x)>0,\quad x\in(0,1),\\
{z_j^\epsilon}'(t)+(K-c_j(t))z_j^\epsilon(t)  - \alpha\left[w_j^\epsilon(t,0)+C^\gamma_j w_{j-1}^\epsilon(t,1)\right] >0,\\
-d\partial_xw_j^\epsilon(t,0)+\alpha w_j^\epsilon(t,0)- \beta z_j^\epsilon(t)>0,\\
d\partial_xw_j^\epsilon(t,1)+\alpha w_j^\epsilon(t,1) - \beta C^{-\gamma}_{j+1} z_{j+1}^\epsilon(t)>0,
\end{cases}
\eqq
for all $t\in(0,T]$, $x\in[0,1]$ and $j\in\Z$ with 
\bqs\left\{
\begin{split}
w_j^\epsilon(0,x) &:=w_j(0,x)+\epsilon \ee^{\delta \left(x-\frac{1}{2}\right)^2}>0,\\
z_j^\epsilon (0)&:= z_j(0)+\epsilon \ee^{\frac{\delta}{4}}>0,
\end{split}\right.
\eqs
and for each $t\in(0,T]$, $x\in[0,1]$  and $j\in\Z$
\bqs\left\{
\begin{split}
w_j^\epsilon(t,x) &\underset{j\rightarrow\pm\infty}{\longrightarrow}\epsilon \ee^{\varrho t +\delta \left(x-\frac{1}{2}\right)^2}>0 ,\\
z_j^\epsilon(t)&\underset{j\rightarrow\pm\infty}{\longrightarrow} \epsilon \ee^{\varrho t +\frac{\delta}{4}}>0.
\end{split}\right.
\eqs
As a consequence, there exists some $J>0$ such that $w_j^\epsilon(t,x)>0$ and $z_j^\epsilon(t)>0$ for all $t\in(0,T]$, $x\in[0,1]$ and $|j|\geq J$. Our aim is to show that this is also true for all $|j|\leq J$. By contradiction, assume that there exists $t_0\in(0,T]$, $j_0\in\llbracket-J,J\rrbracket$ and $x_0\in[0,1]$ such that $w^\epsilon_{j_0}(t_0,x_0)=0$ while $w^\epsilon_{j}(t,x)\geq0$ and $z^\epsilon_{j}(t)\geq0$ for all $t\in(0,t_0]$, $x\in[0,1]$ and $j\in\llbracket-J,J\rrbracket$. If $x_0\in(0,1)$, then by definition we have $\partial_tw^\epsilon_{j_0}(t_0,x_0)\leq0$ and $\partial_x^2w^\epsilon_{j_0}(t_0,x_0)\geq0$ such that
\bqs
0\geq \partial_tw_{j_0}^\epsilon(t_0,x_0) - d\partial_x^2 w_{j_0}^\epsilon(t_0,x_0)+Kw_{j_0}^\epsilon(t_0,x_0)>0,
\eqs
which is a contradiction. If $x_0=0$ then one has $\partial_xw_{j_0}^\epsilon(t_0,0)\geq0$ and the boundary condition gives
\bqs
0\geq-d\partial_xw_{j_0}^\epsilon(t_0,0)+\alpha \underbrace{w_{j_0}^\epsilon(t_0,0)}_{=0}>  \beta z_{j_0}^\epsilon(t_0) \geq 0,
\eqs
which is impossible. A similar argument shows that if $x_0=1$ one also reaches a contradiction. Finally, if on the other hand we had assumed that $z_{j_0}(t_0)=0$ while $w^\epsilon_{j}(t,x)\geq0$ and $z^\epsilon_{j}(t)\geq0$ for all $t\in(0,t_0]$, $x\in[0,1]$ and $j\in\llbracket-J,J\rrbracket$, then using the equation satisfied by $z_{j_0}$ we find
\bqs
0\geq \underbrace{{z_{j_0}^\epsilon}'(t_0)}_{\leq0}+(K-c_j(t))\underbrace{z_{j_0}^\epsilon(t_0)}_{=0} > \alpha\left[w_{j_0}^\epsilon(t_0,0)+C^\gamma_j w_{j_0-1}^\epsilon(t_0,1)\right] \geq 0,
\eqs
which is a contradiction. Let us remark that in the above inequality we have used that $w_{j_0-1}^\epsilon(t_0,1)\geq 0$. This holds by definition of $(t_0,x_0,j_0)$ if $j_0\in\llbracket-J+1,J\rrbracket$, and if $j_0=-J$, the fact that $w_{-J-1}^\epsilon(t_0,1)\geq 0$ holds thanks to the definition of $J$ and the fact that $w_j^\epsilon(t,x)>0$ for all $t\in(0,T]$, $x\in[0,1]$ and $|j|\geq J$.

As a conclusion, we have proved that $w_j^\epsilon(t,x)>0$ and $z_j^\epsilon(t)>0$ for all $t\in(0,T]$, $x\in[0,1]$ and $j\in\Z$. Since $\epsilon>0$ was left arbitrary by passing to the limit $\epsilon\rightarrow0$ we obtain that $w_j(t,x)\geq0$ and $z_j(t)\geq 0$ for all $t\in(0,T]$, $x\in[0,1]$ and $j\in\Z$, which concludes the first part of the proof.

In order to prove the last part of the proposition, we shall instead prove that if there exists $t_0\in(0,T]$ and $j_0\in\Z$ such that $\rho_{j_0}(t_0)=0$ or if there exists $t_0\in(0,T]$, $j_0\in\Z$ and $x_0\in[0,1]$ such that $v_{j_0}(t_0,x_0)=0$ then $\rho_j(t)=0$ and $v_j(t,x)=0$ for all $t\in[0,t_0]$, $x\in[0,1]$ and $j\in\Z$. If $\rho_{j_0}(t_0)=0$, then integrating the equation for $\rho_{j_0}$ from $t=0$ to $t_0$, we obtain that
\bqs
0=\rho_{j_0}(t_0)\geq \ee^{\int_0^{t_0}c_{j_0}(s)\ds}\rho_{j_0}(0)+\alpha\int_{0}^{t_0} \ee^{\int_s^{t_0}c_{j_0}(\tau)\d\tau}\left( v_{j_0}(s,0)+v_{j_0-1}(s,1)\right)\ds\geq0.
\eqs
As a consequence, we deduce that $\rho_{j_0}(0)=0$, $v_{j_0}(t,0)=v_{j_0-1}(t,1)=0$ for all $t\in[0,t_0]$ and thus $\rho_{j_0}(t)=0$ for all $t\in[0,t_0]$. But then, the strong maximum principle applied to $v_{j_0}$ and $v_{j_0-1}$ gives that $v_{j_0}(t,x)=v_{j_0-1}(t,x)=0$ for all $t\in[0,t_0]$ and $x\in[0,1]$. Now, by contradiction, if $\mathbf{v}\neq0$ or $\boldsymbol{\rho}\neq0$ on $[0,t_0]$, without loss of generality, we may assume that there exists $p\in\Z$ with $p>j_0$ and $x_*\in[0,1]$ such that $v_p(0,x_*)>0$. By continuity of $v_p$, there exists $r>0$ such that $v_p(0,x)>0$ for all $x\in[0,1]\cap B_r(x_*)$. Recalling that $v_p$ satisfies
\bqs
\begin{cases}
\partial_tv_p(t,x) - d\partial_x^2 v_p(t,x)\geq0,\quad x\in(0,1),\\
-d\partial_xv_p(t,0)+\alpha v_p(t,0) \geq 0,\\
d\partial_xv_p(t,1)+\alpha v_p(t,1)  \geq 0,
\end{cases}
\eqs
since $\rho_j(t)\geq0$ for all $j\in\Z$, the strong maximum principle implies that $v_p(t,x)>0$ for all $t\in(0,t_0]$ and $x\in[0,1]$. This, in turn, also implies that
\bqs
\rho_{p}(t)\geq \ee^{\int_0^tc_p(s)\ds}\rho_{p}(0)+\alpha\int_{0}^t \ee^{\int_s^tc_p(\tau)\d\tau}\left( v_{p}(s,0)+v_{p-1}(s,1)\right)\ds>0,
\eqs
for all $t\in(0,t_0]$. Now, inspecting the equation satisfied by $v_{p-1}$, we have for all $t\in(0,t_0]$ that
\bqs
\begin{cases}
\partial_tv_{p-1}(t,x) - d\partial_x^2 v_{p-1}(t,x)\geq0,\quad x\in(0,1),\\
-d\partial_xv_{p-1}(t,0)+\alpha v_{p-1}(t,0) \geq 0,\\
d\partial_xv_{p-1}(t,1)+\alpha v_{p-1}(t,1) > 0, \\
v_{p-1}(0,x) \geq 0, \quad x\in[0,1].
\end{cases}
\eqs
Once again, the strong maximum principle implies that $v_{p-1}(t,x)>0$ for all $t\in(0,t_0]$ and $x\in[0,1]$. By induction, we reach a contradiction since we eventually end up proving that $v_{j_0}(t,x)>0$ for $t\in(0,t_0]$ and $x\in[0,1]$, which is impossible.
\end{Proof}

\begin{proposition}\label{propCPapp2}
Let $\mathbf{v}$ and $\boldsymbol{\rho}$ with  $\rho_j\in\mathscr{C}^1([0,+\infty),\R)$ and
\bqs
v_j\in \mathscr{C}^{0}([0,+\infty)\times[0,1],\R),~ \partial_t v_j, \partial_x^2v_j\in\mathscr{C}^{0}((0,+\infty)\times(0,1),\R), \text{ and } \partial_x v_j\in \mathscr{C}^{0}((0,+\infty)\times[0,1],\R),
\eqs
for all $j\in\Z$, which satisfy 
\bqq\label{systvrho2b}
\forall t>0,~ \zeta(t)\leq j \leq \xi(t), \quad \begin{cases}
\partial_tv_j(t,x) - d\partial_x^2 v_j(t,x)\geq0,\quad x\in(0,1),\\
\rho_j'(t)-c_j(t)\rho_j(t)  \geq \alpha\left[v_j(t,0)+v_{j-1}(t,1)\right],\\
-d\partial_xv_j(t,0)+\alpha v_j(t,0) \geq \beta \rho_j(t),\\
d\partial_xv_j(t,1)+\alpha v_j(t,1)  \geq \beta \rho_{j+1}(t),
\end{cases}
\eqq
for some $\mathbf{c}=(c_j)_{j\in\Z}\in L^\infty(\R_+,\ell^\infty(\Z))$ and continuous functions $\zeta:\R_+\to\R$ and $\xi:\R_+\to\R$. Assume that $v_j(0,x)\geq0$ and $\rho_j(0)\geq0$ for all $x\in[0,1]$ and $\zeta(0)-1\leq j \leq \xi(0)+1$ together with $v_j(t,x)\geq0$ and $\rho_j(t)\geq 0$ for all $t>0$, $x\in[0,1]$ and $j\in[\zeta(t)-1,\zeta(t))\cup(\xi(t),\xi(t)+1]$,  then $v_j(t,x)\geq0$ and $\rho_j(t)\geq 0$ for all $t>0$, $x\in[0,1]$ and $\zeta(t)\leq j\leq \xi(t)$.
\end{proposition}

\begin{Proof}
The proof is a direct adaptation of the proof of the previous proposition. Fix $T>0$. By assumption on the sequence $\mathbf{c}$, there exists $K>0$ such that 
\bqs
K-c_j(t)>0,\text{ for} ~ t\in(0,T]\text{ and }j\in\Z.
\eqs
For any $\epsilon>0$, we define
\bqs\left\{
\begin{split}
w_j^\epsilon(t,x) &:=  \ee^{-Kt}v_j(t,x)+\epsilon \ee^{\varrho t +\delta \left(x-\frac{1}{2}\right)^2},\\
z_j^\epsilon(t)&:= \ee^{-Kt} \rho_j(t)+\epsilon \ee^{\varrho t +\frac{\delta}{4}},
\end{split}\right.
\eqs
where $\varrho>0$ and $\delta>0$ are taken large enough to ensure that $w_j^\epsilon(t,x)$ and $z_j^\epsilon(t)$ satisfy
\bqs
\forall t>0,~ \zeta(t)\leq j \leq \xi(t), \quad\begin{cases}
\partial_tw_j^\epsilon(t,x) - d\partial_x^2 w_j^\epsilon(t,x)+Kw_j^\epsilon(t,x)>0,\quad x\in(0,1),\\
{z_j^\epsilon}'(t)+(K-c_j(t))z_j^\epsilon(t)  - \alpha\left[w_j^\epsilon(t,0)+ w_{j-1}^\epsilon(t,1)\right] >0,\\
-d\partial_xw_j^\epsilon(t,0)+\alpha w_j^\epsilon(t,0)- \beta z_j^\epsilon(t)>0,\\
d\partial_xw_j^\epsilon(t,1)+\alpha w_j^\epsilon(t,1) - \beta  z_{j+1}^\epsilon(t)>0,
\end{cases}
\eqs
Furthermore, one also has $w_j^\epsilon(0,x)>0$ and $z_j^\epsilon(0)>0$ for all $x\in[0,1]$ and $\zeta(0)-1\leq j \leq \xi(0)+1$ together with $w_j^\epsilon(t,x)>0$ and $z_j^\epsilon(t)> 0$ for all $t>0$, $x\in[0,1]$ and $j\in[\zeta(t)-1,\zeta(t))\cup(\xi(t),\xi(t)+1]$.

By contradiction, assume that there exists $t_0\in(0,T]$, $j_0\in\llbracket\zeta(t_0),\xi(t_0)\rrbracket$ and $x_0\in[0,1]$ such that $w^\epsilon_{j_0}(t_0,x_0)=0$ while $w^\epsilon_{j}(t,x)\geq0$ and $z^\epsilon_{j}(t)\geq0$ for all $t\in(0,t_0]$, $x\in[0,1]$ and $j\in\llbracket\zeta(t),\xi(t)\rrbracket$.  If $x_0\in(0,1)$, then by definition we have $\partial_tw^\epsilon_{j_0}(t_0,x_0)\leq0$ and $\partial_x^2w^\epsilon_{j_0}(t_0,x_0)\geq0$ such that
\bqs
0\geq \partial_tw_{j_0}^\epsilon(t_0,x_0) - d\partial_x^2 w_{j_0}^\epsilon(t_0,x_0)+Kw_{j_0}^\epsilon(t_0,x_0)>0,
\eqs
which is a contradiction. If $x_0=0$ then one has $\partial_xw_{j_0}^\epsilon(t_0,0)\geq0$ and the boundary condition gives
\bqs
0\geq-d\partial_xw_{j_0}^\epsilon(t_0,0)+\alpha \underbrace{w_{j_0}^\epsilon(t_0,0)}_{=0}>  \beta z_{j_0}^\epsilon(t_0) \geq 0,
\eqs
which is impossible. On the other hand, if $x_0=1$, the boundary condition gives
\bqs
0\geq d\underbrace{\partial_xw_{j_0}^\epsilon(t_0,1)}_{\leq0}+\alpha w_{j_0}^\epsilon(t_0,0)>  \beta z_{j_0+1}^\epsilon(t_0) \geq 0,
\eqs
which is also impossible. The fact that $z_{j_0+1}^\epsilon(t_0) \geq 0$ even if $j_0=\xi(t_0)$ is ensured by the assumption that $z_j^\epsilon(t)> 0$ for $j\in[\zeta(t)-1,\zeta(t))\cup(\xi(t),\xi(t)+1]$. Finally, if  we had assumed that $z_{j_0}(t_0)=0$ while $w^\epsilon_{j}(t,x)\geq0$ and $z^\epsilon_{j}(t)\geq0$ for all $t\in(0,t_0]$, $x\in[0,1]$ and $j\in\llbracket\zeta(t),\xi(t)\rrbracket$, then using the equation satisfied by $z_{j_0}$ we find
\bqs
0\geq \underbrace{{z_{j_0}^\epsilon}'(t_0)}_{\leq0}+(K-c_j(t))\underbrace{z_{j_0}^\epsilon(t_0)}_{=0} > \alpha\left[w_{j_0}^\epsilon(t_0,0)+w_{j_0-1}^\epsilon(t_0,1)\right] \geq 0.
\eqs
Once again, the fact that $w_{j_0-1}^\epsilon(t_0,1)\geq0$ even if $j_0=\zeta(t_0)$ comes from the assumption that $w_j^\epsilon(t,x)> 0$ for $j\in[\zeta(t)-1,\zeta(t))\cup(\xi(t),\xi(t)+1]$ and all $x\in[0,1]$.

To conclude the proof one just passes to the limit $\epsilon\rightarrow0$.
\end{Proof}

\begin{proposition}\label{propSMP0}
Let $\lambda>0$, $\kappa>0$ and $(N,M)\in\Z^2$ such that $N<M$. Consider a sequence $\mathbf{w}=(w_j)_{j=N-1,\dots,M+1}$ satisfying 
\bqs
\lambda ( w_{j+1}-2w_j+w_{j-1})-c_j w_j \geq 0, \quad j=N,\dots,M,
\eqs
for some sequence $\mathbf{c}=(c_j)_{j=N,\dots,M}$ satisfying $0\leq c_j \leq \kappa$ for all $j\in\llbracket N,M\rrbracket$. If $w_j\leq0$ for $j\in\llbracket N-1,M+1\rrbracket$ and there exists $j_*\in\llbracket N,M\rrbracket$ such that $w_{j_*}=0$, then $w_{j}=0$ for all $j\in\llbracket N-1,M+1\rrbracket$.
\end{proposition}

\begin{Proof}
By contradiction, let us first assume that there exists $j_0\neq j_*$ with $j_0\in\llbracket N,M\rrbracket$ such that $w_{j_0}<0$. Without loss of generality, let us assume that $j_*<j_0\leq M$ (otherwise we would assume $N\leq j_0<j_*$ and the argument works similarly). Let us define $z_j=\mathrm{e}^{\alpha (j-j_*)}-1$ for all $j\in\llbracket N-1,M+1\rrbracket$, for some $\alpha>0$ to be fixed later. We first note that
\bqs
z_{j+1}-2z_j+z_{j-1}=2\left( \cosh(\alpha)-1\right)\mathrm{e}^{\alpha (j-j_*)}, \quad j=N,\dots,M.
\eqs
As a consequence, we have
\bqs
\lambda ( z_{j+1}-2z_j+z_{j-1})-c_j z_j = \left[2\lambda\left( \cosh(\alpha)-1\right)-c_j\right]\mathrm{e}^{\alpha (j-j_*)}+c_j, \quad j=N,\dots,M,
\eqs
from which we choose $\alpha>0$ such that
\bqs
2\lambda\left( \cosh(\alpha)-1\right)>\kappa.
\eqs
Next, we set $w^\epsilon_j=w_j+\epsilon z_j$ for $j\in\llbracket N-1,M+1\rrbracket$ and fix $\epsilon>0$ according to
\bqs
0<\epsilon<\frac{-w_{j_0}}{z_{j_0}}.
\eqs
Thus, we have that $w^\epsilon_j<0$ for all $j=N-1,\dots,j_*-1$, $w_{j_*}^\epsilon=0$ and $w_{j_0}^\epsilon<0$ by definition of $\epsilon$, which implies the existence of $j_m\in\llbracket j_*,j_0-1\rrbracket$ where $0\leq w_{j_m}^\epsilon=\max_{j\in\llbracket j_*,j_0-1\rrbracket}w_j^\epsilon$, together with
\bqs
\lambda ( w_{j+1}^\epsilon-2w_j^\epsilon+w_{j-1}^\epsilon)-c_j w_j^\epsilon >0, \quad j=N,\dots,M.
\eqs
This yields a contradiction since at $j=j_m$
\bqs
0\geq \lambda ( w_{j_m+1}^\epsilon-2w_{j_m}^\epsilon+w_{j_m-1}^\epsilon)-c_{j_m} w_{j_m}^\epsilon >0.
\eqs
As a consequence, we have proved that $w_j=0$ for all $j\in\llbracket N,M\rrbracket$, but then we have
\bqs
 \lambda w_{N-1} = \lambda ( w_{N+1}-2w_N+w_{N-1})-c_N w_N \geq 0,
\eqs
and thus $w_{N-1}=0$ since by assumption we assumed that $w_{N-1}\leq0$. A similar argument shows that $w_{M+1}=0$, and this concludes the proof.
\end{Proof}

\begin{proposition}\label{propSMPdiscret}
Let $\lambda>0$ and $(N,M)\in\Z^2$ such that $N<M$. Consider two sequences $\underline{\boldsymbol{\rho}}=(\underline{\rho}_j)_{j=N-1,\dots,M+1}$ and $\overline{\boldsymbol{\rho}}=(\overline{\rho}_j)_{j=N-1,\dots,M+1}$ which are respectively lower bounded and upper bounded, and satisfying 
\bqq\label{ineqSMP}
\lambda ( \underline{\rho}_{j+1}-2\underline{\rho}_j+\underline{\rho}_{j-1})+f(\underline{\rho}_j) \geq \lambda ( \overline{\rho}_{j+1}-2\overline{\rho}_j+\overline{\rho}_{j-1})+f(\overline{\rho}_j),
\eqq
for each $j=N,\dots,M$. If $\overline{\rho}_j\geq\underline{\rho}_j$ for $j\in\llbracket N-1,M+1\rrbracket$ and there exists $j_*\in\llbracket N,M\rrbracket$ such that $\overline{\rho}_{j_*}=\underline{\rho}_{j_*}$, then $\overline{\rho}_j=\underline{\rho}_j$ for all $j=N-1,\dots,M+1$.
\end{proposition}

\begin{Proof}
Since $\underline{\boldsymbol{\rho}}=(\underline{\rho}_j)_{j=N-1,\dots,M+1}$ and $\overline{\boldsymbol{\rho}}=(\overline{\rho}_j)_{j=N-1,\dots,M+1}$  are respectively lower bounded and upper bounded, there exist two constants $\underline{\beta}<\overline{\beta}$ such that
\bqs
\underline{\beta} \leq \underline{\rho}_j \leq \overline{\rho}_j \leq \overline{\beta}, \quad j= N-1,\dots,M+1.
\eqs
Next, we let $K>0$ be the Lipschitz constant of $f$ on the interval $[\underline{\beta},\overline{\beta}]$. We define $\tilde{f}(u):=f(u)+Ku$ which is nondecreasing on $[\underline{\beta},\overline{\beta}]$ by construction. We can then rewrite \eqref{ineqSMP} as
\bqs
\lambda ( w_{j+1}-2w_j+w_{j-1})-Kw_j \geq \tilde{f}(\overline{\rho}_j)-\tilde{f}(\underline{\rho}_j)\geq0, \quad j=N,\dots,M,
\eqs
where we have set $w_j:=\underline{\rho}_j-\overline{\rho}_j$, and the last inequality is a consequence of the monotony of $\tilde{f}$. By assumption, we have $w_j\leq0$ for $j\in\llbracket N-1,M+1\rrbracket$ and there exists $j_*\in\llbracket N,M\rrbracket$ such that $w_{j_*}=0$. From Proposition~\ref{propSMP0}, we infer that $w_j=0$ for all $j=N-1,\dots,M+1$ which concludes the proof.
\end{Proof}

\begin{proposition}\label{propCPappWQ}
Let $\mathbf{W}=(W_j)_{j\in\Z}$ and $\mathbf{Q}=(Q_j)_{j\in\Z}$ with $W_j,Q_j\in\mathscr{C}^1([0,+\infty),\R)$ for all $j\in\Z$ which satisfy
\bqs
\forall t>0,~j\in\Z, \quad \left\{
	\begin{split}
		 W_j'(t)& \geq -2\alpha W_j(t)+\beta(Q_j(t)+Q_{j+1}(t)),\\
		Q'_j(t)&\geq c_j(t)Q_j(t)+\alpha(W_{j}(t)+W_{j-1}(t)),
	\end{split}\right.
\eqs
with some $\mathbf{c}=(c_j)_{j\in\Z}\in L^\infty(\R_+,\ell^\infty(\Z))$. Assume that $W_j(0)\geq0$ and $Q_j(0)\geq0$ for all $j\in\Z$, then $W_j(t)\geq0$ and $Q_j(t)\geq 0$ for all $t>0$ and $j\in\Z$. If furthermore $\mathbf{W}(0)\not\equiv0$ or $\mathbf{Q}(0)\not\equiv0$, then $W_j(t)>0$ and $Q_j(t)> 0$ for all $t>0$ and $j\in\Z$.
\end{proposition}

\begin{Proof}
Fix $T>0$. By assumption on the sequence $\mathbf{c}$, there exists $K>0$ such that 
\bqs
K-c_j(t)>0,\text{ for} ~ t\in(0,T]\text{ and }j\in\Z.
\eqs
For any $\gamma>0$, we define
\bqs\left\{
\begin{split}
w_j(t) &:=  \ee^{-\gamma|j|-Kt}W_j(t),\\
q_j(t)&:= \ee^{-\gamma|j|-Kt} Q_j(t).
\end{split}\right.
\eqs
Since $\mathbf{W}$ and $\mathbf{Q}$ are assumed to be locally bounded, we have for each $t\in(0,T]$ and $j\in\Z$ that
\bqs\left\{
\begin{split}
w_j(t) &\underset{j\rightarrow\pm\infty}{\longrightarrow} 0,\\
q_j(t)&\underset{j\rightarrow\pm\infty}{\longrightarrow} 0
\end{split}\right.
\eqs
The sequences $\mathbf{w}$ and $\mathbf{q}$ now satisfy
\bqq\label{systwq}
\forall t\in(0,T], \quad j\in\Z, \quad\left\{
\begin{split}
 w_j'(t)& \geq -(K+2\alpha) w_j(t)+\beta(q_j(t)+C_{j+1}^{-\gamma}q_{j+1}(t)),\\
		q'_j(t)&\geq -(K-c_j(t))q_j(t)+ \alpha(w_{j}(t)+C_j^\gamma w_{j-1}(t)),
\end{split}\right.
\eqq
with $C^\gamma_j$ defined in \eqref{Cgamma}. As in the proof of Proposition~\ref{propCPapp1}, we define 
\bqs\left\{
\begin{split}
w_j^\epsilon(t) &:=w_j(t)+\epsilon \ee^{\rho t},\\
q_j^\epsilon (t)&:= q_j(t)+\epsilon \ee^{\rho t },
\end{split}\right.
\eqs
with $\epsilon>0$ and $\rho>0$ chosen such that
\bqs
\rho>\max\left( \alpha(1+\ee^\gamma), \beta(1+\ee^\gamma)-K-2\alpha\right).
\eqs
With such a choice, we readily have that
\bqq\label{systwqeps}
\forall t\in(0,T], \quad j\in\Z, \quad \left\{
\begin{split}
 {w_j^\epsilon}'(t) &>  -(K+2\alpha) w_j^\epsilon(t)+\beta(q_j^\epsilon(t)+C_{j+1}^{-\gamma}q_{j+1}^\epsilon(t)),\\
		{q_j^\epsilon}'(t)&> -(K-c_j(t))q_j^\epsilon(t)+ \alpha(w_{j}^\epsilon(t)+C_j^\gamma w_{j-1}^\epsilon(t)),
\end{split}\right.
\eqq
with 
\bqs\left\{
\begin{split}
w_j^\epsilon(0) &:=w_j(0)+\epsilon >0,\\
q_j^\epsilon (0)&:= q_j(0)+\epsilon>0,
\end{split}\right.
\eqs
and for each $t\in(0,T]$  and $j\in\Z$
\bqs\left\{
\begin{split}
w_j^\epsilon(t) &\underset{j\rightarrow\pm\infty}{\longrightarrow}\epsilon \ee^{\varrho t}>0 ,\\
q_j^\epsilon(t)&\underset{j\rightarrow\pm\infty}{\longrightarrow} \epsilon \ee^{\varrho t}>0.
\end{split}\right.
\eqs
As a consequence, there exists some $J>0$ such that $w_j^\epsilon(t)>0$ and $q_j^\epsilon(t)>0$ for all $t\in(0,T]$ and $|j|\geq J$. Our aim is to show that this is also true for all $|j|\leq J$. By contradiction, without loss of generality, assume that there exists $t_0\in(0,T]$ and $j_0\in\llbracket-J,J\rrbracket$ such that $w^\epsilon_{j_0}(t_0)=0$ while $w^\epsilon_{j}(t)\geq0$ and $q^\epsilon_{j}(t)\geq0$ for all $t\in(0,t_0]$ and $j\in\llbracket-J,J\rrbracket$. Then by definition, we also have that ${w^\epsilon_{j_0}}'(t_0)\leq0$ such that
\bqs
0\geq  {w_{j_0}^\epsilon}'(t_0) + (K+2\alpha) w_{j_0}^\epsilon(t_0)-\beta(q_{j_0}^\epsilon(t_0)+C_{j_0+1}^{-\gamma}q_{j_0+1}^\epsilon(t_0)) >0,
\eqs
which is a contradiction. As a consequence, one has $w_j^\epsilon(t)>0$ and $q_j^\epsilon(t)>0$ for all $t\in(0,T]$ and $j\in\Z$. Since $\epsilon>0$ was left arbitrary by passing to the limit $\epsilon\rightarrow0$, we obtain that $w_j(t)\geq0$ and $q_j(t)\geq 0$ for all $t\in(0,T]$ and $j\in\Z$, which concludes the first part of the proof.

In order to prove the last part of the proposition, we shall instead prove that if there exists $t_0\in(0,T]$ and $j_0\in\Z$ such that $W_{j_0}(t_0)=0$ or $Q_{j_0}(t_0)=0$ then $W_j(t)=0$ and $Q_j(t)=0$ for all $t\in[0,t_0]$ and $j\in\Z$. Without loss of generality, suppose that $W_{j_0}(t_0)=0$, then we obtain
\bqs
0 = W_{j_0}(t_0) \geq \ee^{-(K+2\alpha) t_0}W_{j_0}(0)+\beta\int_0^{t_0} \ee^{-(K+2\alpha)(t_0-s)}\left(Q_{j_0}(s)+Q_{j_0+1}(s)\right)\ds\geq0,
\eqs
from which we infer that $W_{j_0}(0)=0$ and $Q_{j_0}(s)=Q_{j_0+1}(s)=0$ for all $s\in[0,t_0]$. As a consequence, we deduce that $W_{j_0}(s)=0$ for all $s\in[0,t_0]$. Now, using the equation satisfied by $Q_{j_0}$ and $Q_{j_0+1}$, we infer that $W_{j_0-1}(s)=0$ and $W_{j_0+1}(s)=0$ for all $s\in[0,t_0]$. As a consequence, by induction, we get that $W_j(t)=0$ and $Q_j(t)=0$ for all $t\in[0,t_0]$ and $j\in\Z$, which concludes the proof.
\end{Proof}

\begin{proposition}
Let $\mathbf{W}=(W_j)_{j\in\Z}$ and $\mathbf{Q}=(Q_j)_{j\in\Z}$ with $W_j,Q_j\in\mathscr{C}^1([0,+\infty),\R)$ for all $j\in\Z$ which satisfy
\bqs
\forall t>0,~\zeta(t)\leq j \leq \xi(t), \quad \left\{
	\begin{split}
		 W_j'(t)& \geq -2\alpha W_j(t)+\beta(Q_j(t)+Q_{j+1}(t)),\\
		Q'_j(t)&\geq c_j(t)Q_j(t)+\alpha(W_{j}(t)+W_{j-1}(t)),
	\end{split}\right.
\eqs
with some $\mathbf{c}=(c_j)_{j\in\Z}\in L^\infty(\R_+,\ell^\infty(\Z))$ and continuous functions $\zeta:\R_+\to\R$ and $\xi:\R_+\to\R$. Assume that $W_j(0)\geq0$ and $Q_j(0)\geq0$ for all $\zeta(0)-1\leq j \leq \xi(0)+1$ together with $W_j(t)\geq0$ and $Q_j(t)\geq 0$ for all $t>0$ and $j\in[\zeta(t)-1,\zeta(t))\cup(\xi(t),\xi(t)+1]$, then we have $W_j(t)\geq0$ and $Q_j(t)\geq 0$ for all $t>0$ and $\zeta(t)\leq j \leq \xi(t)$.
\end{proposition}

\begin{Proof}
Fix $T>0$. By assumption on the sequence $\mathbf{c}$, there exists $K>0$ such that 
\bqs
K-c_j(t)>0,\text{ for} ~ t\in(0,T]\text{ and }j\in\Z.
\eqs
Next, let us define
\bqs\left\{
\begin{split}
w_j^\epsilon(t) &:=\ee^{-Kt}W_j(t)+\epsilon \ee^{\rho t},\\
q_j^\epsilon (t)&:= \ee^{-Kt}Q_j(t)+\epsilon \ee^{\rho t },
\end{split}\right.
\eqs
for $\epsilon>0$ and 
\bqs
\rho>\max\left(2\alpha, 2\beta-K-2\alpha\right).
\eqs
As a consequence $(w_j^\epsilon,q_j^\epsilon)$ satisfies all the assumptions with strict inequalities and we can argue as in the previous proof. Without loss of generality, we assume by contradiction that there exists $t_0\in(0,T]$ and $j_0\in [\zeta(t_0),\xi(t_0)]\cap\Z$ such that $w_{j_0}^\epsilon(t_0)=0$ while $w^\epsilon_{j}(t)\geq0$ and $q^\epsilon_{j}(t)\geq0$ for all $t\in(0,t_0]$ and $\zeta(t)\leq j \leq \xi(t)$. Then by definition, we also have that ${w^\epsilon_{j_0}}'(t_0)\leq0$ such that
\bqs
0\geq  {w_{j_0}^\epsilon}'(t_0) + (K+2\alpha) w_{j_0}^\epsilon(t_0)-\beta(q_{j_0}^\epsilon(t_0)+q_{j_0+1}^\epsilon(t_0)) >0,
\eqs
which is a contradiction. In the above inequality, we crucially used our assumption that $w_j^\epsilon(t)>0$ for $j\in(\xi(t),\xi(t)+1]$. We can then pass to the limit $\epsilon\rightarrow0$ and conclude the proof.
\end{Proof}

\begin{proposition}\label{propCPappWQ2b}
Let $(\underline{\mathbf{W}},\underline{\mathbf{Q}})$ and $(\overline{\mathbf{W}},\overline{\mathbf{Q}})$ with $\underline{W}_j,\underline{Q}_j, \overline{W}_j,\overline{Q}_j\in\mathscr{C}^1([0,+\infty),\R)$ for all $j\in\Z$ be respectively subsolution and supersolution of the asymptotic system \eqref{asympsyst} for $t>0$ and $j\in[\zeta(t),\xi(t)]\cap\Z$
for some continuous functions $\zeta:\R_+\to\R$ and $\xi:\R_+\to\R$. Assume that $(\underline{W}_j(0),\underline{Q}_j(0))\leq (\overline{W}_j(0),\overline{Q}_j(0))$  for all $\zeta(0)-1\leq j \leq \xi(0)+1$ together with $(\underline{W}_j(t),\underline{Q}_j(t))\leq (\overline{W}_j(t),\overline{Q}_j(t))$ for all $t>0$ and $j\in[\zeta(t)-1,\zeta(t))\cup(\xi(t),\xi(t)+1]$, then we have $(\underline{W}_j(t),\underline{Q}_j(t))\leq (\overline{W}_j(t),\overline{Q}_j(t))$ for all $t>0$ and $\zeta(t)\leq j \leq \xi(t)$.
\end{proposition}
\begin{Proof}
The proof is a direct consequence of the previous proposition.
\end{Proof}


\begin{thebibliography}{99}	

\bibitem{Berestycki-Roquejoffre-Rossi-2013}
H. Berestycki, J.-M. Roquejoffre and L. Rossi.
 The influence of a line of fast diffusion in Fisher-KPP propagation.
 \newblock \emph{J Math Biol}, 66 (2013), 743--766.
 
\bibitem{BRR16}
H.  Berestycki, J. M. Roquejoffre and L. Rossi.
\newblock Travelling waves, spreading and extinction for Fisher-KPP propagation driven by a line with fast diffusion.
\newblock \emph{ Nonlinear Analysis}, 137, (2016) 171--189.

\bibitem{BF21a} 
C. Besse and G. Faye.
\newblock Dynamics of epidemic spreading on connected graphs.
\newblock \emph{J. Math. Biol}, 82 (2021), 1--52.

\bibitem{BF21b} 
C. Besse and G. Faye.
\newblock Spreading properties for SIR models on homogeneous trees.
\newblock \emph{Bull. Math. Biol}, 83:114 (2021) , pp. 1--27 .


\bibitem{Besse-Faye-Roquejoffre-Zhang-2023}
 C. Besse, G. Faye, J.-M. Roquejoffre and M. Zhang.
  The logarithmic Bramson correction for Fisher-KPP equations on the lattice $\Z$.
  \newblock \emph{Trans. Am. Math. Soc.} 376 (2023), pp. 8553-8619.

\bibitem{PCB22}
P.C. Bressloff. 
\newblock Local accumulation time for diffusion in cells with gap junction coupling. 
\newblock \emph{Phys. Rev. E}, 105(3), (2022) 034404.


\bibitem{DLPZ}
Y. Du, B. Lou, R. Peng and M. Zhou. 
\newblock The Fisher-KPP equation over simple graphs: varied persistence states in river networks.
\newblock \emph{ J. Math. Biol.}, 80(5), (2020) 1559-1616.

\bibitem{FHT}
W. T. L. Fan, W. Hu and G. Terlov. 
\newblock Wave propagation for reaction-diffusion equations on infinite random trees.
\newblock \emph{Commun. Math. Phys.}, 384(1), (2021), 109--163.

 
\bibitem{Faria-Rambaut-2014}
N. Faria, A. Rambaut, M. Suchard, G. Baele, T. Bedford, M. Ward, A. Tatem, J. Sousa, N. Arinaminpathy, J. Pepin, D. Posada, M. Peeters, O. Pybus and P. Lemey.
\newblock HIV epidemiology. the early spread and epidemic ignition of hiv-1 in human populations.
\newblock  \emph{Science}, 346 (2014), 56--61.

\bibitem{Gatto-Bertuzzo-2020}
M. Gatto, E. Bertuzzo, L. Mari, S. Miccoli, L. Carraro, R. Casagrandi and A. Rinaldo.
\newblock Spread and dynamics of the COVID-19 epidemic in Italy: effects of emergency containment measures.
\newblock \emph{ Proc. Nat. Acad. Sci.}, (2020), 10484--10491.

\bibitem{Gou-Ward-2016}
J. Gou and M.J. Ward.
\newblock Oscillatory dynamics for a coupled membrane-bulk diffusion model with Fitzhugh-Nagumo membrane kinetics.
\newblock \emph{SIAM J Appl Math}, 76(2) (2016), 776--804.


\bibitem{Gou-Li-Nagata-Ward-2015}
J. Gou, Y.X. Li, W. Nagata and M.J. Ward.
\newblock Synchronized oscillatory dynamics for a 1-Dmodel ofmembrane kinetics coupled by linear bulk diffusion.
\newblock  \emph{SIAM J Appl Dyn Syst}, 14(4) (2015), 2096--2137.


\bibitem{Hoffman-Holzer-2019} 
A. Hoffman and M. Holzer.
\newblock Invasion fronts on graphs: the Fisher-KPP equation on homogeneous trees and Erd\H{o}s--R\'eyni graphs.
\newblock \emph{Discrete Contin Dyn Syst B}, 24(2) (2019), 671.


\bibitem{Hale-Thomas-2019}
C.M. Hale, C. Thomas, et al..
\newblock Spatiotemporal heterogeneity in the distribution of chikungunya and Zika virus case incidences during their 2014 to 2016 epidemics in Barranquilla, Colombia.
\newblock \emph{Int. J. Environ. Res. Public Health}, 16 (2019), 1759.

\bibitem{JPS}
Y. Jin, R. Peng and J. Shi.
\newblock Population dynamics in river networks. 
\newblock \emph{J. Nonlinear Sci.}, 29(6), (2019), 2501-2545.

\bibitem{KS09}
J. P. Keener and J. Sneyd.
\newblock Mathematical Physiology I: Cellular Physiology.
\newblock \emph{Springer, New York,} (2009) 2nd edition.

\bibitem{KDB}
H. Kravitz, C. Dur\'on and M. Brio. 
\newblock A Coupled Spatial-Network Model: A Mathematical Framework for Applications in Epidemiology.
\newblock \emph{ Bull. Math. Biol.} 86.11 (2024): 132.


\bibitem{Paquin-Lefebvre-Nagata-Ward-2020} F. Paquin-Lefebvre, W. Nagata and M.J. Ward.
\newblock  Weakly nonlinear theory for oscillatory dynamics in a one-dimensional PDE-ODE model of membrane dynamics coupled by a bulk diffusion field.
\newblock \emph{SIAM J Appl Math}, 80(3) (2020), 1520--1545.

\bibitem{LSU}
O. A Ladyzhenskaja, V.A. Solonnikov and N.N. Ural'ceva.
\newblock Linear and Quasi-linear Equations of Parabolic Type.
\newblock \emph{American Mathematical Soc.}, vol 23, (1968).

\bibitem{PW}
M. H. Protter and H.F. Weinberge.
\newblock Maximum principles in differential equations.
\newblock \emph{ Springer Science \& Business Media} (2012).

\bibitem{RB90}
S. V. Ramanan and P. R. Brink.
\newblock Exact solution of a model of diffusion in an infinite chain or monlolayer of cells coupled by gap junctions.
\newblock \emph{ Biophys. J.} 58, 631 (1990).

\bibitem{Weinberger-1982} 
H. Weinberger.
\newblock Long-time behavior of a class of biological models.
\newblock \emph{SIAM J Math Anal}, 13(3) (1982), 353--396.



\end{thebibliography}
\end{document}